\documentclass[sn-mathphys-num]{sn-jnl}
\usepackage[ps]{skak}
\usepackage{graphicx}%
\usepackage{multirow}%
\usepackage{amsmath,amssymb,amsfonts}%
\usepackage{amsthm}%
\usepackage{mathrsfs}%
\usepackage[title]{appendix}%
\usepackage{xcolor}%
\usepackage{textcomp}%
\usepackage{manyfoot}%
\usepackage{booktabs}%
\usepackage{algorithm}%
\usepackage{algorithmicx}%
\usepackage{algpseudocode}%
\usepackage{listings}%
\usepackage{dsfont}
\usepackage{amssymb}   % symboles maths (mathbb, etc.)
\usepackage{amsmath}   % align environment
\usepackage{amsthm}    % theorems environment
\usepackage{stmaryrd}  % pour les doubles crochets-\rrbracket-\rbr
\usepackage{dsfont}    % pour l'indicatrice         \mathds{1}
\usepackage{enumitem}

\usepackage{tikz}      % faire de jolis dessins
\usetikzlibrary {arrows.meta} 
\usetikzlibrary{patterns} %hachures
\usepackage{pgfplots}
\pgfplotsset{width=6.5cm,compat=1.9}

\usepackage{hyperref}   
\hypersetup{
    colorlinks=true,
    linkcolor=purple,
    filecolor=magenta,      
    urlcolor=cyan,
    citecolor=teal
    }

\newcommand{\N}{\mathbb{N}}
\newcommand{\Sf}{\mathbb{S}}
\newcommand{\E}{\mathbb{E}} 
\newcommand{\T}{\mathbb{T}}

\newcommand{\Z}{\mathbb{Z}} 
\newcommand{\R}{\mathbb{R}} 

\newcommand{\Lp}{\mathbb{L}}
\newcommand{\Pro}{\mathbb{P}}

\newcommand{\lbr}{\llbracket}
\newcommand{\rbr}{\rrbracket}

\newcommand{\ind}{\mathds{1}}

\newcommand{\dd}{\mathrm{d}}

\newcommand{\CQFD}{\hfill $\square$}

\newcommand{\sif}{\mbox{ if }}

\newcommand{\und}{\mbox{ and }}

\theoremstyle{plain}
\newtheorem{lemm}{Lemma}[subsection]
\newtheorem{theo}{Theorem}[section]
\newtheorem{prop}{Proposition}[subsection]

\theoremstyle{remark}
\newtheorem*{prof}{Proof}

\newtheorem{deaf}{Definition}[section]

\newcommand{\supess}{\mathop{\smash{\mathrm{esssup}}}}

\everymath{\displaystyle} 

\usepackage{anyfontsize}

\usepackage{pgfplots}
\pgfplotsset{width=6.5cm,compat=1.9}
\usepackage{hyperref}   % hyperliens /!\ en dernier /!\
\hypersetup{
    colorlinks=true,
    linkcolor=purple,
    filecolor=magenta,      
    urlcolor=cyan,
    citecolor = blue
    }

\newcommand{\mZ}{\mathcal{Z}}
\newcommand{\mD}{\mathcal{D}}
\newcommand{\mX}{\mathcal{X}}
\newcommand{\mF}{\mathcal{F}}
\newcommand{\mN}{\mathcal{N}}
\newcommand{\mC}{\mathcal{C}}
\newcommand{\mB}{\mathcal{B}}
\newcommand{\mG}{\mathcal{G}}
\newcommand{\mP}{\mathcal{P}}
\newcommand{\mT}{\mathcal{T}}
\newcommand{\mS}{\mathcal{S}}
\newcommand{\mK}{\mathcal{K}}

\newcommand{\mM}{\mathcal{M}}
\newcommand{\mH}{\mathcal{H}}

\newcommand{\mL}{\mathcal{L}}
\newcommand{\mI}{\mathcal{I}}
\newcommand{\mJ}{\mathcal{J}}

\newcommand{\bV}{\mathbf{V}}
\newcommand{\bO}{\mathbf{O}}
\newcommand{\bF}{\mathbf{F}}
\newcommand{\fm}{\mathfrak{m}}

\newcommand{\ux}{\underline{x}}
\newcommand{\us}{\underline{s}}

\newcommand{\uv}{\underline{v}}
\newcommand{\uz}{\underline{z}}
\newcommand{\uze}{\underline{\zeta}}

\newcommand{\ut}{\underline{t}}
\newcommand{\uom}{\underline{\omega}}
\newcommand{\usig}{\underline{\sigma}}
\newcommand{\utau}{\underline{\tau}}
\newcommand{\um}{\underline{m}}

\newcommand{\ul}{\underline{\ell}}

\newcommand{\ukap}{\underline{\kappa}}
\newcommand{\uchi}{\underline{\smash{\chi}}}

\newcommand{\col}{\mathrm{c}}
\newcommand{\coll}{\mathrm{col}}
\newcommand{\ov}{\mathrm{ov}}
\newcommand{\tree}{\mathrm{tree}}
\newcommand{\cycle}{\mathrm{[cycle]}}

\newcommand{\cy}{\mathrm{c}}

\newcommand{\cls}{\mathrm{agg}}

\numberwithin{equation}{section}

\newcommand{\triple}[1]{{\left\vert\kern-0.25ex\left\vert\kern-0.25ex\left\vert #1 
    \right\vert\kern-0.25ex\right\vert\kern-0.25ex\right\vert}}

\everymath{\displaystyle} 

\begin{document}

\title[Article Title]{Cumulants of the Rayleigh gas mixture model:\\ statistical results}

\subtitle{Fluctuations and large deviations of the nonideal Rayleigh gas}


\author[1]{\fnm{Florent} \sur{Foug\`{e}res}}\email{florent.fougeres@ens.fr}

\date{2026}

\affil[1]{\orgdiv{DMA}, \orgname{\'{E}cole normale sup\'{e}rieure}, \orgaddress{\street{45 rue d'Ulm}, \city{Paris}, \postcode{75005},  \country{France}}}

\abstract{In this paper, we explore the statistical subtleties of the nonideal Rayleigh gas, in a grand canonical mixture framework. This model allows to consider a large amount of tagged particles close to equilibrium, and their empirical measure, whose first-order convergence has been shown to converge to the solution of the linear Rayleigh--Boltzmann equation~\cite{2016brownian}. Thanks to the study of the cumulants of the system, we analyze the asymptotic behaviour of the fluctuations and large deviations of this empirical measure, hence refining the previous statistical results in the same vein as~\cite{2023grandev}. This way, we exhibit the trivial limit behaviour of the fluctuations in any overdilute regime, proving the exact relevance at any statistical scale of the low density limit. In the case of large deviations, we present the linear Boltzmann--Hamilton--Jacobi system driving their asymptotic behaviour.

Eventually, we optimize the geometrical estimates on the billiards dynamics~\cite{2023longcor} to finally achieve a full convergence rate for the cumulants.}


\ackno{I would like to thank warmly my PhD directors Isabelle Gallagher and Sergio Simonella, for their tireless help and support. Thanks also to all the persons who partake in making the department a living, stimulating and cosy place.}

\maketitle

\section{Introduction: model and main results} \label{sec:intro}

\subsection{Mathematical context} 
The Boltzmann equation statistically describes the dynamics of rarefied gas, and may be seen as an intermediate mesoscopic scale between microscopic and macroscopic points of view: it can thus be used to derive fluid mechanics equations from Newton classical equations~(see for example~\cite{1991BGL,2001ns,2008euler}). This model allows on the one hand quantitative simulations~\cite{2022lbm, 2023lbm}, and on the other hand qualitative comprehension of the intrinsic statistical behaviour of fluid matter. In fact, introducing the \emph{entropy} of the system as one of its Lyapunov functional, Boltzmann proved that the density of particles irreversibly converges towards an equilibrium on large time scales. This equilibrium is called the Maxwellian state, where the gas occupies uniformly all the available space, and where the velocities of the particles are distributed according to the Maxwell Gaussian distribution~(\ref{def:maxwellian}), which only depends on the temperature of the system.

The mathematical derivation of this equation from microscopic Newton equations has only been proved in 1975 by Oscar Erasmus Lanford III~\cite{1975lanford, 1976lanford}, in the case of hard sphere interactions, yet on very small times. The main hinderance to long time scales is the emergence of correlations between colliding particles: the quantitative estimates on the system's chaoticity deteriorate over time, preventing from dealing efficiently with recollisions of particles.

Nonetheless, in a framework close to thermodynamic equilibrium (i.e. the Maxwellian state presented above), the statistical stability of the dynamics guarantees the preservation of a certain amount of chaos over large times, allowing a stronger control of correlations. In a \emph{Rayleigh~gas setting}, that describes the behaviour of a small fraction of tagged particles near equilibrium~\cite{1988spohnkineq}, Henk van Beijeren, Lanford, Joel Louis Lebowitz and Herbert Spohn proved in 1980 the derivation \emph{on large time scales} of a linear version of the Boltzmann equation, called the linear Rayleigh--Boltzmann equation~\cite{1980linear,1982lebspohn}.

In 2013, Isabelle Gallagher, Laure Saint-Raymond and Benjamin Texier returned to the work of Lanford and his former student, Francis Gordon King~\cite{1975king}, to generalize their proof in the case of compactly supported potentials, hence providing precise estimates on the convergence rate in Lanford's theorem. These quantitative estimates allowed in 2016 Thierry Bodineau, Gallagher and Saint-Raymond to study the convergence rate in the linear case, along with its dependence on the long time scaling, so as to infer Brownian hydrodynamic limits~\cite{2016brownian}. These estimates on the convergence rates have been improved in~2024 in our previous paper~\cite{fou24}.

Joined by Sergio Simonella, the authors also studied the symmetric linearization of the Boltzmann equation on long time scales~\cite{2023longcor, 2024longder}, precising their study so as to capture rarer dynamics events that are involved in fine correlations between particles. Eventually, they completed their analysis of the fine scales of the dynamics in 2023~\cite{2023grandev}, studying the finest correlations of the system through its \emph{cumulants}~\eqref{def:cumulants}: they exhibited the equations governing the fluctuations of the empirical measure and its large deviations.

The present work is dedicated to performing this study in the Rayleigh~gas setting, where the expected limit of the empirical measure follows the linear Rayleigh--Boltzmann equation. The relevant model to introduce the empirical measure of tagged particles is a gas mixture, introduced in the companion paper~\cite{fou26Rayleigh}, in which most of the particles are initially distributed according to the equilibrium state, and in the low-density regime. Meanwhile, some other tagged particles are initially perturbed from this equilibrium, their number going to infinity, yet necessarily in a more dilute regime, to remain in a Rayleigh setting. We choose the grand canonical ensemble to get symmetry and simplify the study of cumulants. This model allows to study the statistical behaviour of the perturbed particles, and permits in particular to explore the phase transition around the low-density scaling. Indeed, looking at the fluctuations of the empirical measure, we show that there is no threshold ahead of this scaling, which is thus the exact regime to consider, to observe correlations at the limit.

On the other hand, we exhibit a linear version of the Hamilton--Jacobi equation that describes the limit cumulants of the tagged particles, and derive a large deviation principle in a tighter functional framework than~\cite{2023grandev}. These two results are nevertheless still restricted to small time scales, because of the absence of strong \emph{a priori} bounds on the cumulants, due to their non-linear combinatorial structure. Recently, Yu Deng, Zaher Hani and Xiao Ma~\cite{2025DHM} showed that assuming uniform bounds on the regular solutions to the Boltzmann equation, and computing a very fine combinatorial analysis of the dynamics trees, the derivation was valid for long times. Their method might be applied to statistical results like the ones presented in this paper, to extend their time validity, especially since the uniform bounds of the limit are known in the linear framework.

Note that the way we deal with the geometry of recollisions in this paper is different from how it is done in the companion paper~\cite{fou26Rayleigh}: we improve here the analysis performed in the most recent geometric method to get a convergence rate in a full power of~$\varepsilon$ instead of~$\varepsilon |\log \varepsilon|$ (see Section~\ref{app:geom_recoll}).

Section~\ref{sec:intro} recalls the definition of the model we introduced in~\cite{fou26Rayleigh}, and states the main fluctuation and large deviation results. Section~\ref{sec:cumstat} is dedicated to the study of the cumulants of the nonideal Rayleigh gas, providing the equations that drive their evolution in terms of dynamics trees and interaction clusters, short-time bounds and the study of their convergence and limit. Our optimized geometrical computation, necessary to prove the latter results, is postponed to Section~\ref{app:geom_recoll}. Eventually, Section~\ref{sec:fluct_LD_proof} explains how the results from Section~\ref{sec:cumstat} lead to the fluctuation and large deviation results.

Two appendixes follow to expose technical results, about the change of parametrization of the hard sphere scattering rule (Appendix~\ref{app:scattering}), and the functional analysis of the linear Rayleigh--Boltzmann and Boltzmann--Hamilton--Jacobi equations (Appendix~\ref{app:RB}).

\subsection{Microscopic hard sphere model} \label{sec:model2}

The model that we present in the following sections has been introduced in the companion paper~\cite{fou26Rayleigh}. At the microscopical scale, we consider the hard sphere model: the state of a gas of $N$ particles is completely determined by the positions (in the $d$-dimensional torus $\T^d$) and the velocities of every particle, represented by the vector
\begin{equation*}
\uz_N = (z_1, \dots, z_N) \doteq (\ux_N, \uv_N) \in \mD^N \doteq (\T^{d} \times \R^{d})^N.
\end{equation*}
The \emph{hard sphere} model consists in the following exclusion condition, which asserts that two particles cannot get closer than a certain diameter~$\varepsilon > 0$: the positions must belong to the hard sphere exclusion set
\begin{equation} \label{def:dom}
\mathcal{X}^\varepsilon_N \doteq \{ \ux_N \in \T^{dN};\ \forall\ i\neq j,\ d(x_i,x_j) > \varepsilon \},
\end{equation}
where $d(\cdot, \cdot)$ denotes the distance on the torus. Hence, the state of the gas~$\uz_N$ has to remain in the open domain~$\mD^\varepsilon_N \doteq \mX^\varepsilon_N \times \R^{dN}$.

 Within this set, the particles' dynamics is given by the Newton equations for uniform line movement. On the contrary, on the boundary of~$\mathcal{D}^\varepsilon_N$, at least two particles (let us call them $i$ and $j$) collide, by the condition~$|x_i - x_j| = \varepsilon$. Then, if the scalar product~$(x_i - x_j) \cdot (v_i - v_j)$ is positive, it means that the particles are exiting the collision in uniform line movement; but otherwise they are entering the collision and must scatter according to the following system giving the post-collisional velocities~$({v_i}', {v_j}')$:
\begin{equation} \label{eq:coll}
\left\{ \begin{array}{l}
{v_i}' = v_i - \left\langle v_i - v_j, \frac{x_i - x_j}{\varepsilon}\right \rangle\frac{x_i - x_j}{\varepsilon} \\
\ \\
{v_j}' = v_j + \left\langle v_i - v_j, \frac{x_i - x_j}{\varepsilon}\right \rangle\frac{x_i - x_j}{\varepsilon}\cdotp
\end{array} \right.
\end{equation}
In the hard sphere case, the interaction is taken instantaneous and elastic: the system~(\ref{eq:coll}) stems from the preservation of momentum and kinetic energy. 
Appendix~\ref{app:scattering} is dedicated to the parametrization of this scattering system by either one of the following angles,
\begin{equation} \label{eq:scattering_angles}
\omega = \frac{x_i - x_j}{\varepsilon} \mbox{ \ \ \ or \ \ \  } \sigma = \frac{{v_i}' - {v_j}'}{|{v_i}' - {v_j}'|}\cdotp
\end{equation} 
This dynamics is well-defined up to a zero measure set of initial configurations, in which infinite amounts of collisions might happen in finite times, along with collisions between more than two particles at a time. This result was proved by~Roger Keith Alexander~\cite{1975alexander}, and might also be found in~\cite{2013newton}.

\subsection{Statistical grand canonical description}
\label{sec:stat_descr}
As the number of particles~$N$ gets large, this hard sphere dynamics becomes very difficult to compute, especially because it is very chaotic, so that we choose to describe the gas statistically. This paper is dedicated to the study of a gas of identical particles, yet divided in two distinguishable parts, represented by tags~$\ul_N = (\ell_1, \dots, \ell_N) \in \{0,1\}^N$: the tag~$\ell = 0$ will be attributed to `tagged'~particles initially distributed at thermodynamic equilibrium, and the tag~$\ell = 1$ to particles initially perturbed away from that equilibrium.

Our grand canonical model, introduced in~\cite{fou26Rayleigh}, randomizes the number~$N$ of particles to a random variable~$\mN$, the expectancy of which is tuned by a parameter~$\mu$ called the \emph{chemical potential}.
The kinetic limit that we consider is called the \emph{low density limit}, or Boltzmann--Grad limit, and consists in letting this chemical potential~$\mu$ go to infinity while keeping a constant mean free path $\mu^{-1} \varepsilon^{1-d} = 1$, so that the particles' diameter~$\varepsilon$ goes to 0.
In this limit, assuming initial chaos, the first marginal of the density usually converges to the solution to the Boltzmann equation~\cite{1975lanford, 2013newton, 2025DHM}.

This solution of the Boltzmann equation, studied by Boltzmann and Maxwell, when well defined, relaxes in large times to an equilibrium called the \emph{Maxwell state} and defined~\cite{1872boltzmann} as
\begin{equation} \label{def:maxwellian}
M_\beta(x,v) \doteq \left(\frac{\beta}{2\pi} \right)^{d/2} \exp\left( -\frac{\beta}{2} |v|^2 \right).
\end{equation}
The parameter~$\beta$ stands for the inverse temperature of the system, tuning its intensive (kinetic) energy. 
It appears that the density $M_\beta^{\otimes N} \ind_{\mD^\varepsilon_N}$ is an equilibrium of the microscopic hard sphere dynamics, and in this paper we consider specific initial conditions that are close to this thermodynamic equilibrium, to retrieve a \emph{linear version} of the Boltzmann equation, whose theory is much simpler and hence might be derived for long time scales. More precisely, we consider the \emph{nonideal Rayleigh gas mixture model} introduced in the companion paper~\cite{fou26Rayleigh}: we work in the grand canonical ensemble and we add the additional tagging variable~$\ul_N$ indicating to which set each particle belongs. 

The particles at equilibrium are taken in the usual low density (Boltzamnn--Grad) limit, whereas the tagged perturbed particles only occupy a tiny fraction of the gas, smaller than the Boltzmann--Grad density: otherwise indeed they would behave like a classical Boltzmann dilute gas, satisfying the non-linear Boltzmann equation. This all boils down to the following scaling,
\begin{equation}
  \tag{$S_{\varepsilon, \mu, \lambda}$}
  \mu \varepsilon^{d-1} = 1 \mbox{\ \ and \ \ } 1 \ll \lambda \ll \mu,
\label{eq:scaling}
\end{equation}
where~$\mu > 0$ corresponds to the chemical potential of the particles at equilibrium, and $\lambda > 0$ to that of tagged particles. The notation~$\lambda \ll \mu$ simply means that $\lambda \mu^{-1}$ goes to 0 as $\lambda$ and $\mu$ go to infinity. Formally, \textbf{the particles at initial equilibrium will be tagged with a~0, and the initially perturbed `tagged' particles will be tagged with a~1}. The tags of all particles hence form a vector $\ul_n \in \Lambda_n \doteq \{0,1\}^n$, identified to the corresponding subset $ \ul_n \subset \lbr 1, n \rbr$, with the following notation 
\begin{equation*}
|\ul_{n}| = \| \ul_n \|_1 = \left\lvert \{ i \leq n \ , \ \ell_i = 1 \} \right\rvert  \mbox{\ \ and \ \ } \varphi_0^{\otimes \ul_{n}} (\uz_{\ul_{n}}) = \prod_{\substack{i \leqslant n \\ \ell_i = 1}} \varphi_0(z_i).
\end{equation*}
Indeed, while the majority of particles are initially distributed according to equilibrium, the tagged particles will be perturbated by an initial perturbation~$\varphi_0 \in \Lp^\infty_x \Lp^\infty_v(M_{\beta/2})$, happening in the phase space (contrary to~\cite{2016brownian, fou24} where it only happened in positions).
Now, the mixed grand canonical ensemble consists in relaxing the number of particles, weighting it with a mixed Poisson law depending on the number of tagged particles. 

More precisely, at fixed~$\lambda$, and $\mu = \varepsilon^{1-d}$, we defined in~\cite{fou26Rayleigh} the \emph{correlation functions}~$(F_n^\varepsilon)$ of the system, describing its evolution, based on the marginals of the canonical densities of the gas. Denoting $\tilde{\ul}_{n+p} \doteq (\ul_n, \ul^*_p)$ and $\tilde{\uz}_{n+p} \doteq (\uz_n, \uz^*_p)$, these correlation functions initially satisfy~\cite{fou26Rayleigh}
\begin{align}  \label{eq:initial_corr_funct}
F_n^\varepsilon(0, \uz_n, \ul_n)&  = \sum_{p\geqslant 0} \sum_{\ul^*_p \in \Lambda_p } \frac{\lambda^{|\ul^*_{p}|} \mu^{p - |\ul^*_{p}|} }{p!} \int \dd \uz^*_{p}\ M_\beta^{\otimes n+p}(\tilde{\uv}_{n+p}) \varphi_0^{\otimes \tilde{\ul}_{n+p}}(\tilde{\uz}_{\tilde{\ul}_{n+p}})  \ind_{\mX_{n+p}^\varepsilon}(\tilde{\ux}_{n+p}),
\end{align}
where the normalizing (grand canonical) partition function is defined as follows 
\begin{align} \label{def:partition_function_Z}
\mZ_{\mu} & \doteq \sum_{p\geqslant 0} \sum_{\ul_p \in \Lambda_p } \frac{\lambda^{|\ul_{p}|} \mu^{p - |\ul_{p}|} }{p!} \int M_\beta^{\otimes p}(\uv_{p}) \varphi_0^{\otimes \ul_{p}}(\uz_{\ul_{p}})  \ind_{\mX_{p}^\varepsilon}(\ux_{p}) \dd \uz_{p}.
\end{align}

Our probabilistic study is based on the random variables $(Z_{\varepsilon, i}^{[t]}, L_i)_{1\leqslant i \leqslant \mN}$, giving the states at time~$t$, and the tags of the particles. The probability density of the initial state $(Z_{\varepsilon, i}^{[0]}, L_i)$ is given by the initial correlation functions above~\eqref{eq:initial_corr_funct}, and the evolution at time~$t$ is a deterministic piecewise affine function of the initial state. 
The correlation functions are in fact defined such that for any observable $H_n \in \mC^\infty_{\mathrm{c}}(\mD^n \times \Lambda_n)$, for any~$t \geqslant 0$, we have  (see once again~\cite{fou26Rayleigh})
\begin{align}
\E\left[ \sum_{1 \leq i_k \neq i_j \leq \mathcal{N}} H_n(Z_{\varepsilon,i_1}^{[t]}, L_{i_1}, \dots, Z_{\varepsilon, i_n}^{[t]}, L_{i_n})  \right] 
= \sum_{\ul_n  \in \Lambda_n} \mu^{n - |\ul_n|} \lambda^{|\ul_n|}  \int_{\mathcal{D}^\varepsilon_n} F_n^\varepsilon(t,\uz_n, \ul_n)  H_n(\uz_n, \ul_n) \dd \uz_n. \label{eq:empmes}
\end{align} 
This formula between the correlation functions and observables will be  used in the combinatorial computations of Section~\ref{sec:cumulant_generating_function}, when expanding the cumulant generating function.

By the usual BBGKY arguments~\cite{fou26Rayleigh, fouthese}, the correlation functions satisfy the following hierarchy in our mixed scaling~\eqref{eq:scaling}
\begin{align} \label{eq:BBBGKY}
\partial_t F_n^\varepsilon + \uv_n \cdot \nabla_{\ux_n} F_n^\varepsilon  = \mC^{\langle 0 \rangle}_n  F^\varepsilon_{n+1} + \frac{\lambda}{\mu} \mC^{\langle 1 \rangle}_n F^\varepsilon_{n+1},
\end{align}
where the \emph{collision operators} are defined as
\begin{align} 
\mC_n^\ell F_{n+1}^{\varepsilon} = \sum_{i=1}^n  \int &\dd \omega \dd v_{n+1} \langle \omega, v_{n+1} - v_i \rangle_+ \times \nonumber \\ 
&  \Bigl[F_{n+1}^{\varepsilon}(\uz_n', \ul_n, x_i + \varepsilon \omega, v_{n+1}', \ell) - F_{n+1}^{\varepsilon}(\uz_n, \ul_n, x_i - \varepsilon \omega, v_{n+1}, \ell) \Bigr] . \label{def:collision_op2}
\end{align}
Morally, we look at the influence of a $(n+1)$-th particle---with tag~$\ell$---on the dynamics, colliding with one of the $n$ existing ones with angle~$\omega$ and velocity~$v_{n+1}$, whence the name \emph{collision} operators.
 The \emph{cross section} $\langle \omega, v_{n+1} - v_i\rangle_+$ weights the likelihood of such a collision. 

\subsection{Linear Rayleigh--Boltzmann equation} \label{sec:RB}

In the case of a single tagged particle, it has been shown~\cite{2016brownian} that the first correlation function~$F_1^\varepsilon$ converges in the low-density limit to the solution $g\doteq M_\beta \varphi$ of the linear \emph{Rayleigh--Boltzmann equation} with initial condition~$\varphi_0$:
\begin{equation} \label{eq:phi} \def\arraystretch{1.7}
\left\{\begin{array}{rcl}
\partial_t \varphi + v\cdot \nabla_x \varphi & =&  \int_{\Sf^{d-1}} \int_{\R^d} [\varphi(v') - \varphi(v)] M_\beta(v_c) \langle \omega, v_c - v \rangle_+ \dd v_c \dd \omega, \\
\varphi(0,x,v) &=& \varphi_0. \end{array} \right.
\end{equation}
This linear equation is globally well-posed in the velocity-weighted space $\Lp^\infty_x \Lp^\infty_v(M_{\beta/2})$ (see Appendix~\ref{app:RB}), and allows to derive the linear heat equation in the hydrodynamic limit~\cite{2016brownian}.

The formal limit of the hierarchy on the correlation functions~\eqref{eq:BBBGKY} with initial conditions~\eqref{eq:initial_corr_funct}, is satisfied by the family~$(M_\beta^{\otimes n} \varphi^{\otimes \ell_n})_{n \geqslant 1}$, with~$\varphi$ the solution to the Rayleigh--Boltzmann equation~\eqref{eq:phi} with initial data~$\varphi_0$. This result is formalized in the companion paper~\cite[Theorem~1]{fou26Rayleigh}, where we show the long-time convergence of the correlation functions with a quantitative convergence~rate.

\subsection{Main results: fluctuations and large deviations} \label{sec:stat_refinements}
% numérotation corollary ??§§
As a consequence to the convergence of the correlation functions, we also prove in~\cite[Corollary~1.1]{fou26Rayleigh} a law of large numbers for the empirical measure defined as follows.
For any observable~$H \in \mC^\infty_{\mathrm{c}}(\mD \times \Lambda_1)$ we define the following random variables: the \emph{empirical measure of all particles}
\begin{equation} \label{def:emp_meas}
\pi^\varepsilon_t[H] \doteq \frac{1}{\mu} \sum_{i=1}^\mN H(Z^{[t]}_{\varepsilon,i}	, L_i),
\end{equation}
and the \emph{empirical measure of tagged particles}
\begin{equation} \label{def:emp_meas_tag}
\tilde{\pi}^\varepsilon_t[H] \doteq \frac{1}{\lambda} \sum_{i=1}^\mN H(Z^{[t]}_{\varepsilon,i}, L_i) \ind_{L_i = 1}.
\end{equation}

To study further than this first order convergence, we will introduce its fluctuation field around its expected value, and present its limit. Eventually, we state a large deviation principle for the empirical measure, that we prove in the following sections.

Indeed, the empirical measure of the tagged particles, defined in~\eqref{def:emp_meas_tag} for an observable~$H$, can be seen as the observation of a measure~$\tilde{\pi}^\varepsilon_t \in \mM(\mD)$ on the domain~$\mD$, writing 
\begin{equation} \label{def:emp_meas_tag2}
\tilde{\pi}^\varepsilon_t[H] = \int H(z) \dd \tilde{\pi}^\varepsilon_t(z).
\end{equation}

The family~$(\tilde{\pi}^\varepsilon_s)_{0 \leqslant s \leqslant t}$ defines a measure on the trajectories of $\mD^{[0,t]}$. The set $\mathrm{Traj}([0,t], \mM(\mD))$ of such measures is endowed with the Skorokhod topology.
On the other hand, one may define the \emph{fluctuation field}
\begin{equation} \label{def:fluctuation_field}
\zeta^\varepsilon_t = \sqrt{\lambda} \bigl( \tilde{\pi}^\varepsilon_t - \E\left[ \tilde{\pi}^\varepsilon_t\right] \bigr),
\end{equation}
to capture the next small order after the law of large numbers~\cite[Corollary~1.1]{fou26Rayleigh}.

\begin{theo}[Convergence of the fluctuation field] \label{theo:fluct}There exists a time~$T>0$ such that, in the mixed scaling~\eqref{eq:scaling}, the fluctuation field defined above converges in law on~$[0,T]$ to a Gaussian process~$(\zeta_t)$, whose equal-time covariance is given for any $t < T$ by
\begin{equation}\label{eq:cov_fluct}
\E\Bigl[\zeta_t[g]\zeta_t[h]\Bigr] = \int_\mD M(v) \varphi(t,z) g(z) h(z) \dd z,
\end{equation}
where~$\varphi$ denotes the solution to the linear Rayleigh--Boltzmann equation~\eqref{eq:phi}.
\end{theo}
Indeed, unlike in~\cite{2023grandev} where the limit fluctuation field satisfies a linear stochastic equation of the form 
\begin{equation*} 
\dd \tilde{\zeta}_t = \mL \tilde{\zeta}_t\ \dd t + \dd \tilde{\eta}_t,
\end{equation*}
here the limit fluctuation field is trivial and does not depend on the second cumulant (see next Section~\ref{sec:cumulant_generating_function}); there is no interference between the tagged particles, on the mere condition that $\lambda \ll \varepsilon^{1-d}$, without any phase transition between this scaling and the nonlinear scaling~$\lambda \sim \alpha \varepsilon^{1-d}$.
The proof of this theorem is given in Section~\ref{sec:fluct_proof}, thanks to finer objects than the correlation functions, called~\emph{cumulants}, that capture finer scales of the dynamics and allow to rescale its rare events to characterize them.

To study the \emph{large deviations} of the dynamics, we will harness a weaker topology than the Skorokhod one. For a measure~$\fm = (\fm_s)_{s \in [0,t]} \in \mathrm{Traj}([0,t], \mM(\mD))$, and an observable~$h \in \mC^\infty_c([0,t] \times \mD)$, we define the \emph{filtered mean}
\begin{equation} \label{def:crochet}
\{ h, \fm \}_t = \int_{\mD} h(t,z) \dd \fm_t(z) - \int_0^t \int_{\mD} (\partial_s + v\cdot \nabla_x ) h(s,z) \dd \fm_s(z),
\end{equation}
that filters the transported part of the considered observables. The first quantity that will be relevant for the large deviation principle will be the limit object~$\mI(t,h)$, that is defined in~\eqref{def:Itg} as the limit of an expectancy, and which can also be defined as the mild solution (Proposition~\ref{prop:identifHJ}) to the Hamilton--Jacobi system
\begin{equation} \label{syst:HJ1}
\left\{ \begin{array}{lll}
(\partial_s - v\cdot \nabla_x) q^{[t]} = \frac{\partial \mH}{\partial p} (q^{[t]}, p^{[t]})  & , & q^{[t]}(0) = M \varphi_0, \\
(\partial_s - v\cdot \nabla_x)(p^{[t]} - h) = -\frac{\partial \mH}{\partial q}(q^{[t]}, p^{[t]}) & , & p^{[t]}(t) = h(t),
\end{array} \right.
\end{equation}
with the Hamiltonian
\begin{equation*}
 \mH\left( q , p \right) \doteq \int \dd v_2 \dd \omega \langle v_{1}  - v_{2},  \omega \rangle_+ M_\beta(v_2) q(z_1) (e^{p(z_1) - p(z_1')} - 1),
\end{equation*}
in the sense of
\begin{align} 
\mI(t,h) = \mI(0,h) + \int_0^t \dd s \int_{\mD} q^{[t]}(s) (\partial_s - v\cdot \nabla_x)(p^{[t]}(s) - h(s))  + \int_0^t  \mH(q^{[t]}(s), p^{[t]}(s)) \dd s.
\end{align}
More precisely, we are interested in its Legendre transform, defined for~$ \mathfrak{v} \in \mathrm{Traj}([0,t], \mM(\mD))$ as
\begin{equation} \label{def:Legendre}
\mathbf{\Lambda}(t, \mathfrak{v}) \doteq \sup_{h \in \mathbb{B}_{t, \beta}} \Bigl[ \{h, \mathfrak{v} \} - \mI(t,h) - 1 \Bigr],  
\end{equation} 
where the supremum is taken over observables in~$\mathbb{B}_{t, \beta}$, defined in \eqref{def:gset} as the set of observables with bounded transport, and such that~$e^h$ is uniformly dominated by the inverse $\beta/4$--Gaussian. 

Eventually, for the lower bound of the large deviation principle, we need to consider the set~$\mathbf{S}_t$ of strong solutions, on~$[0,t]$, of a biased linear Boltzmann equation of the form
\begin{equation} \label{eq:biasedBoltz} 
(\partial_s - v\cdot \nabla_x) \mathfrak{v} = \int \dd v_\cy \dd \omega \langle v  - v_\cy,  \omega \rangle_+   M_\beta(v_\cy) \left(\mathfrak{v}(v')  e^{p(z) -  p(z')} - \mathfrak{v}(v) e^{p(z') - p(z)} \right),
\end{equation}
for some~$p \in  \mathbb{B}_{t, \beta}$. The large deviation principle might then be formulated as follows.
\begin{theo}[Large deviations of the empirical measure] \label{theo:LD}
In our mixed scaling~\eqref{eq:scaling}, considering the tagged empirical measure~\eqref{def:emp_meas_tag2}, there exists a time~$T > 0$ such that, for any~$t \in (0,T]$, we have the following large deviation upper bound 
\begin{align} \label{ineq:upperbound}
\limsup_{\varepsilon \to 0} \frac{1}{\lambda} \log \Pro \bigl(\tilde{\pi}^\varepsilon_t \in \bF\bigr) \leq - \inf_{\mathfrak{v} \in \bF}  \mathbf{\Lambda}(t, \mathfrak{v}),
\end{align}
when~$\bF$ is a closed set in the Skorokhod topology. Additionally, when~$\bO$ is an open set in this topology, one has the large deviation lower bound
\begin{align}  \label{ineq:lowerbound}
\liminf_{\varepsilon \to 0} \frac{1}{\lambda} \log \Pro \bigl(\tilde{\pi}^\varepsilon_t \in \bO\bigr) \geq - \inf_{\mathfrak{v} \in \bO \cap \mathbf{S}_t} \mathbf{\Lambda}(t, \mathfrak{v}).
\end{align}
\end{theo}
Note that the most useful result is the upper bound on the probability of deviation, which has no restriction on its infimum. Nevertheless, the lower bound, which precises that the upper bound is optimal, is here restricted (like in~\cite{2023grandev}) to solutions not that far from the Boltzmann linear equation, since they must be solutions to the biased equation~\eqref{eq:biasedBoltz}. The proof of this theorem is the subject of Section~\ref{sec:LD_proof}, once again based on the convergence of the cumulants. For this reason, since without a priori bounds on these objects our method does not allow to show their convergence on large times, the theorems above are restricted as in~\cite{2023grandev} to short times.

\section{Cumulants of the nonideal Rayleigh gas} \label{sec:cumstat}

\subsection{Cumulant generating function} \label{sec:cumulant_generating_function}

The cumulant generating function is the functional $\log \E\left[ \exp\left( \mu \pi^\varepsilon_t[H] \right) \right]$,
containing all the information on the moments of the empirical measure.
Using the identity~\eqref{eq:empmes} between the correlation functions and observables, the formal expansion of this functional leads naturally~\cite{2023grandev} to objects called \emph{cumulants}, which we define below. We then prove rigourously the said formal expansion of the cumulant generating function~\eqref{def:cumulant_generating_function2}.

\subsubsection{Cumulants} \label{sec:cum_def}

The definition of the cumulants is based on a decomposition into partitions, so that for $\sigma \in \mP_n$ a partition of $\lbr 1, n \rbr$, we denote $|\sigma|$ the number of subsets $(\sigma_i)_{1 \leqslant i \leqslant |\sigma|}$ that compose this partition, which is not to be confused with the cardinal~$|\sigma_i|$ of one of these subsets. Eventually, we denote~$\mP_n^k \subset \mP_n$ the set of partitions~$\sigma \in \mP_n$ that contain exactly $k = |\sigma|$ subsets.
\begin{deaf}[Cumulants] The cumulants associated to a family~$(G_n)_{n\geqslant 1}$ are defined as
\begin{equation} \label{def:cumulants}
g_n(\uz_n, \ul_n) \doteq \sum_{\sigma \in \mP_n} (-1)^{|\sigma|-1} (|\sigma|-1)! G_{[\sigma]}(\uz_n, \ul_n),
\end{equation}
where for any partition~$\sigma \in \mP_n$, we denote
\label{page:G_sig}
\begin{equation*}
G_{[\sigma]}(\uz_n, \ul_n) = \prod_{i = 1}^{|\sigma|} G_{|\sigma_i|}(\uz_{\sigma_i}, \ul_{\sigma_i}).
\end{equation*}
Note that the $n$-th cumulant~$g_n$ is constructed from all the correlation functions~$(G_i)$ for $i \leq n$. Indeed, it decomposes the interactions between $n$~particles into products of interactions within the subsets of every possible partition, to measure the defects of independence. One may easily check that for a tensorized family $(G_n = G_1^{\otimes n})$, all the cumulants $g_n$ vanish for~$n\geq 2$. 

We denote~$(f^\varepsilon_n)$ the cumulants associated to the hierarchy~$(F_n^\varepsilon)$. Note for example that the second cumulant is $f_2^\varepsilon = F_2^\varepsilon - {F_1^\varepsilon}^{\otimes 2}$, encoding the defect of independence between pairs of particles. One will see in Section~\ref{sec:cum_expansion} that when expanding the pseudo-trajectory formula, this corresponds to the rare dynamics in which a distinguished couple of particles interact together.
The cumulants of the exclusion indicators~$(\ind_{\mX^\varepsilon_n})$ are denoted~$\phi_n$, studied in further detail in the companion paper~\cite{fou26Rayleigh}.
\end{deaf}
\begin{prop}[Inversion formula] \label{prop:inversion_formula} With the definition above, the injectivity of cumulants is a consequence of the following inversion formula
\begin{equation} \label{eq:inversion_cumulants}
G_n = \sum_{\sigma \in \mP_n} g_{[\sigma]}.
\end{equation}
\end{prop}
The proof is written for example in~\cite[Proposition~2.2.1]{2023grandev} or our PhD thesis~\cite[Proposition~5.1.1]{fouthese}.

\subsubsection{Expanding the cumulant generating function} \label{sec:exp_cum_gen_fun}

We go now to the proof of an expansion for the cumulant generating function. 
By linearity, and using the fact that before having been assigned a tag, the particles are exchangeable, one has
\begin{align*}
\E\left[ \exp\left( \sum_{i = 1}^\mN H(Z_{\varepsilon, i}^{[t]}, L_i) \right) \right] & = 1 + \sum_{k \geqslant 1} \frac{1}{k!} \E\left[ \left( \sum_{i = 1}^\mN H(Z_{\varepsilon, i}^{[t]}, L_i)\right)^k \right] \\
& = 1 + \sum_{k \geqslant 1} \frac{1}{k!} \E\left[ \sum_{n=1}^k  \frac{1}{n!} \sum_{\substack{k_1 + \dots + k_n = k \\ (k_j \geqslant 1)_{j\leq n}}} \frac{k!}{k_1!\dots k_n!} \sum_{\substack{ (i_j \leqslant \mN)_{j\leq n} \\ i_j \neq i_{j'}}} H_{i_1}^{k_1}\dots  H_{i_n}^{k_n}   \right],
\end{align*}
denoting $H_{i} \doteq H(Z_{\varepsilon, i}^{[t]}, L_i)$, and expanding the power $k$ by partitioning the resulting sum according to the number of different particles implied in each product. 
Hence, using the relation~\eqref{eq:empmes} between observables and correlation functions,  and permuting the sums to get rid of the sum over $k$, we get (once again extending the correlation functions by~0 outside of their domain)
\begin{align}
\E\left[ \exp\left( \sum_{i = 1}^\mN H(Z_i^{t}, L_i) \right) \right] & = 1 + \sum_{n \geqslant 1} \frac{1}{n!} \sum_{\ul_n \in \Lambda_n} \lambda^{|\ul_n|} \mu^{n - |\ul_n|} \int F_n^{\varepsilon}(t, \ul_n) \prod_{i=1}^n \left( \sum_{k_i \geqslant 1} \frac{1}{k_i!} H(z_i,\ell_i)^{k_i} \right)  \nonumber \\
& = 1 + \sum_{n \geqslant 1} \frac{1}{n!} \sum_{\ul_n \in \Lambda_n} \lambda^{|\ul_n|} \mu^{n - |\ul_n|}   \int F_n^\varepsilon(t, \ul_n) (e^H - 1)^{\otimes n}(\ul_n) . \label{eq:gen_corr_func}
\end{align}
Now that we dispose of an expansion for this expectancy in terms of correlation functions, let us see how the cumulants appear along with a logarithm. We start from the formula~\eqref{eq:gen_corr_func} above and use the inversion formula~\eqref{eq:inversion_cumulants}. Hence, using the exchangeability of particles to reduce the partitions to the number of elements in each of their subsets, we get
\begin{align*}
&\E\left[ \exp\left( \sum_{i = 1}^\mN H(Z_i^{[t]}, L_i) \right) \right] \\
& \hspace{0mm} = 1 + \sum_{n \geqslant 1} \frac{1}{n!} \sum_{\ul_n \in \Lambda_n} \lambda^{|\ul_n|} \mu^{n - |\ul_n|} \sum_{\sigma \in \mP_n} \prod_{i=1}^{|\sigma|} \int_{\mD^{|\sigma_i|}} f^\varepsilon_{\sigma_i}(t, \ul_{\sigma_i}) (e^H - 1)^{\otimes \sigma_i}(\ul_{\sigma_i}) \\
&  \hspace{0mm} = 1 +  \sum_{n \geqslant 1} \frac{1}{n!} \sum_{s=1}^n \frac{1}{s!}\sum_{\substack{p_1 + ... + p_s = n \\ (p_i \geqslant 1)_{i\leqslant s}}} \frac{n!}{p_1!\dots p_s!} \prod_{i=1}^s \sum_{\ul^{(i)} \in \Lambda_{p_i}} \lambda^{|\ul^{(i)}|} \mu^{p_i - |\ul^{(i)}|} \int_{\mD^{p_i}} f^\varepsilon_{p_i}(t, \ul^{(i)}) (e^H - 1)^{\otimes p_i}(\ul^{(i)}), \nonumber
\end{align*}
where the denominator~$s!$ stems from the arbitrary order that we impose on the partitions' subsets. We have also split the labels~$\ul_n \in \Lambda_n$ into the labels~$\ul^{(i)} \in \Lambda_{p_i}$ on each subset. Eventually, we sum over~$n$ to relax the condition on the subset cardinals~$p_i$, and factorize everything as
\begin{align*}
\E\left[ \exp\left( \sum_{i = 1}^\mN H(Z_i^{[t]}, L_i) \right) \right] 
&  = 1 +  \sum_{s \geqslant 1} \frac{1}{s!} \left( \sum_{p \geqslant 1}  \frac{1}{p!} \sum_{\ul_p \in \Lambda_p } \lambda^{|\ul_p|} \mu^{p - |\ul_p|} \int_{\mD^p} f^\varepsilon_{p}(t, \ul_{p}) (e^H - 1)^{\otimes p}(\ul_{p}) \right)^s, \nonumber
\end{align*}
which makes appear the exponential of the quantity defined below.
\begin{deaf}[Cumulant generating function]
Thanks to the computation above, the two following definitions are equivalent, defining the \emph{cumulant generating function}:
\begin{align} \label{def:cumulant_generating_function}
 \mathfrak{G}_\varepsilon^{[t]}[H]& \doteq \log \E\left[ \exp\left( \sum_{i = 1}^\mN H(Z_i^{[t]}, L_i) \right) \right] \\[1em]
  & \doteq \sum_{p \geqslant 1}  \frac{1}{p!} \sum_{\ul_p \in \Lambda_p } \lambda^{|\ul_p|} \mu^{p - |\ul_p|} \int_{\mD^p} f^\varepsilon_{p}(t, \ul_{p}) (e^H - 1)^{\otimes p}(\ul_{p}). \label{def:cumulant_generating_function2}
\end{align}
Note that it is not directly a generating function as one can be used to, since it is not an expansion in powers of the observable~$H$, but in powers of its exponential, which makes appear combinations of cumulants when deriving along~$H$, not directly cumulants.
\end{deaf}
Let us finally observe that if the observable is of the form $\hat{H}(z, \ell) = h(z) \ind_{\ell = 1}$, (i.e. counting only the tagged particles) then the cumulant generating function writes
\begin{equation}
\mathfrak{G}_\varepsilon^{[t]}[\hat{H}] = \sum_{p \geqslant 1}  \frac{\lambda^p}{p!} \int_{\mD^p} f^\varepsilon_{p}(t, \underline{1}_p) (e^{\hat{H}} - 1)^{\otimes p}
. \label{def:cumulant_generating_function_marq}
\end{equation}
We do not renormalize yet the cumulant generating function, since according to whether the observable~$H$ weights all the particles or only the tagged ones, the suitable scale will be~$\mu$ or~$\lambda$.

\subsection{Pseudo-trajectories} \label{sec:pseudo-traj}
The first step to understand the cumulants is to find an equation on them implying only other cumulants, which we do by computing an expansion of the dynamics, based on the interactions between the particles~(Lemma~\ref{lem:cum_pseudo_traj}). We start by defining a pseudo-trajectory formulation of the hierarchy on the correlation functions.

We henceforth choose to denote $p_\mu \doteq \frac{\lambda}{\mu}$ the fraction of initially perturbed particles, so that iterating Duhamel formula as in~\cite{2013newton} or~\cite{2016brownian}, we can write the Dyson~expansion
\begin{equation} \label{eq:Dyson}
F_n^\varepsilon(t) = \sum_{k\geqslant 0} \sum_{\ul^*_k \in \Lambda_k } p_\mu^{|\ul^*_k|} Q_{n, \ul^*_k}(t) F^\varepsilon_{n+k}(0),
\end{equation}
developing the choice of the encountered tags~$\ul^*_k \doteq (\tilde{\ell}_{n+1}, \dots, \tilde{\ell}_{n+k})$, with the successive-collision operators defined as
\begin{equation} \label{def:successive_collision_operator}
Q_{n, \ul^*_k}(t) \doteq \int_{T_k(t)} \Theta_n(t-t_1)\mathcal{C}^{\tilde{\ell}_{n+1}}_n \Theta_{n+1}(t_1 - t_2)\dots \mC^{\tilde{\ell}_{n+k}}_{n+k-1} \Theta_{n+k}(t_k)  \dd \ut_{k},
\end{equation}
where $\Theta_n(\tau)$ denotes the transport semi-group operator in $\mD^\varepsilon_n$ with specular reflections, for a time~$\tau$. The collision times are integrated over
\begin{equation} \label{def:time_ensemble}
T_k(t) \doteq \Bigl\{\ \ut_k \ \Bigl | \ 0 \doteq t_{k+1} \leq t_k \leq \dots \leq t_1 \leq  t_0 \doteq t \Bigr\}.
\end{equation}
The main idea of the proof, coming from Lanford's original paper~\cite{1975lanford}, is to use a coupling between this expansion and its limit version, implying imaginary histories of the particles, that eventually lead to the state $\uz_n$ at time~$t$. These histories, called \emph{pseudo-trajectories}, are non-physical trajectories that---in a way---allow to extend the method of characteristics for the successive-collision operators. 

Indeed, the transport operators appearing in~\eqref{def:successive_collision_operator} correspond to following the characteristics of free transport, with specular reflections: taking the first operator $\Theta_n(t-t_1)$ of a functional amounts to consider this functional at time $t_1$,  in a state~$\uz_n^{[t_1]}$ given by the hard sphere dynamics. 

Then, the first collision operator~\eqref{def:collision_op2} writes
\begin{align*} 
&\mC_n^\ell F_{n+1}^\varepsilon = \sum_{i=1}^n \sum_{s_{1} = \pm 1} s_{1} \int \dd \omega_1 \dd v_{n+1} \langle \omega_1, v_{n+1} - v_i \rangle_+  F_{n+1}^\varepsilon(\uz_n^{\langle s_{1}\rangle}, \ul_n, x_i + s_1 \varepsilon \omega_1, v_{n+1}^{\langle s_{1}\rangle}, \ell),
\end{align*}
where $\uz_n^{\langle +1 \rangle} = \uz_n'$ and $\uz_n^{\langle -1 \rangle} = \uz_n$, scattered for the gain term, and let unchanged for the loss term, so that the collision is always incoming, allowing to pursue the backwards method of characteristics with the next transport operator. Hence, for given collision parameters~$(i, s_1, \omega_1, v_{n+1})$, this operator can be seen as a weighted adjunction of a particle to the characteristics---or pseudo-trajectory---which scatters (or not, according to $s_1$) with particle~$i$, creating a new state~$\uz_{n+1}^{[t_1]} \doteq (\uz_{n}^{[t_1]}, x_i+s_1 \varepsilon \omega_1, v_{n+1}^{\langle s_1 \rangle})$. The integration and sum over these collision parameters will yield an integral over pseudo-trajectories. Iterating this extended method of characteristics and tracking the pseudo-trajectories $(\uz_{n+j}^{[t_j]})$ thus constructed, we bring the analysis back to the value of the functional at time~$\tau = 0$, in the state~$\uz_{n+k}^{[0]}$. 

We will have to record the labels of the existing particles meeting the new one, the velocities of the particles that spring up, the angles at which the encounter happens, and whether they scatter or not.
The pseudo-trajectories will also keep track of the tags of the encountered particles.

Here is precisely how we construct the pseudo-trajectories.
The choice of the successive encountered tags is registered in $\ul^*_k = (\tilde{\ell}_{n+1}, \dots, \tilde{\ell}_{n+k})$, and expanding all the sums in all the collision operators~(\ref{def:successive_collision_operator}), we can sum it up to the history~$(m_1, \dots, m_k)$ of which particle encountered the $(n+i)$-th new one. These particles naturally belong to the following set 
\begin{equation*} \label{eq:Mnk}
\mM_{n,k} \doteq \Bigl\{ (m_1, \dots, m_k) \ \Bigl| \ \forall i \leq k,  m_i \leq n+i-1 \Bigr\}.
\end{equation*}
We consider the scattering labels~$(s_1, \dots, s_k) \in \{ \pm 1\}^k$. The fact that some encounters do not scatter, along with the fact that particles are artificially added, is why the pseudo-trajectories are not physical trajectories.
Once the total history 
\begin{equation} \label{def:history}
\uchi_k \doteq (\um_k, \ul^*_k, \us_k) \in \mH_{n,k} \doteq \mM_{n,k} \times \Lambda_k \times \{\pm 1\}^k
\end{equation}
is fixed, for given collision parameters~$(\uom_k, v_{n+1}, \dots, v_{n+k})$, we can construct the pseudo-trajectories for every~$\uz_n =\uz_n^{[t]}$, backwards in time to a configuration~$\uz_{n+m}^{[0]}$, following the inductive procedure below:
\begin{equation} \label{eq:pseudotraj}
\left\{ \def\arraystretch{1.5} \begin{array}{l}
\uz_n^{[t]} \doteq \uz_n  \\
\forall \ i \in \lbr 0, k \rbr,\ \forall \tau \in (t_{i+1}, t_i),\ \uz_{n+i}^{[\tau]} \mbox{ follows (backwards) the physical hard sphere dynamics}\\
\forall \ i \in \lbr 1, k \rbr,\ \uz_{n+i}^{[t_{i}]} = \left(\uz_{n+i-1}^{[t_i^+],\langle s_{i} \rangle }, x_{m_i} + s_i \varepsilon \omega_i, v_{n+1}^{\langle s_{i} \rangle}\right).
\end{array} \right.
\end{equation}
One may observe that the change of velocities in the last step is automatic by the hard sphere dynamics' boundary condition, but it will not be for the limit version of pseudo-trajectories, since the limit particles are formally pointwise.
In the end, one can write the \emph{pseudo-trajectory formulation} of the Dyson expansion
\begin{equation}
 \label{eq:pseudo_traj_form}
F_{n}^\varepsilon(t) = \sum_{k \geqslant 0} \sum_{\uchi_k} p_\mu^{|\ul^*_k|} \int_{T_k(t)} \dd \ut_k \int \dd \underline{\omega}_k \dd v_{n+1} \dots \dd v_{n+k} \prod_{i=1}^{k} s_i \langle \omega_{i},v_{n+i} - v_{m_i}^{[t_i^+]}  \rangle_+ F_{n+k}^\varepsilon\bigl(0, \uz_{n+k}^{[0]}, \tilde{\ul}_{n+k} \bigr).
\end{equation} 
A small technical detail lies in the fact that the added particles must satisfy the exclusion condition. A way to deal with it may be to impose a condition on the domain on integration of the collision angles~\cite{2023grandev}, yet here to simplify we merely change the definition of the pseudo-trajectories: if at any moment the exclusion condition is violated by the adjunction of a particle, then the trajectories are frozen in this state until time~$\tau = 0$, so that the integral formally vanishes thanks to the initial distribution~$F_{n+k}^\varepsilon(0)$ being~0 outside of $\mD^\varepsilon_{n+k}$. 

The formal limits~$(\uze_{n+i}^{[t_i]})$ of the pseudo-trajectories are defined as their hard sphere versions~\eqref{eq:pseudotraj} for $\varepsilon = 0$, with the noticeable difference that in the dynamics on each time interval~$(t_{i+1}, t_i)$, the particles are pointwise and hence follow the free flow without any scattering. They will be used in Section~\ref{sec:conv_cumulants} to identify the limit cumulants.

\subsection{Pseudo-trajectory measure} \label{sec:pseudo-traj-meas}

We start from the pseudo-trajectory formulation~\eqref{eq:pseudo_traj_form}; recall that a pseudo-trajectory is fully determined by its final configuration~$\uz_n = \uz_n^{[t]}$, and the following parameters: number of collisions, a collision history, collision times, collision angles and collision velocities, summed up in the parametrizing vector
\begin{equation*}
\Psi_n \doteq (k, \uchi_k, \ut_k, \uom_k, \uv_k^*) \in \coprod_{k\geqslant 0} \{k\} \times  \mH_{n,k} \times T_k(t) \times \left( \Sf^{d-1} \times \R^d \right)^k,
\end{equation*} 
referring to the definitions of the collision times~\eqref{def:time_ensemble} and collision history~\eqref{def:history} 
\begin{equation*}
\uchi_k \doteq (\um_k, \ul^*_k, \us_k) \in \mH_{n,k} \doteq \mM_{n,k} \times \Lambda_k \times \{\pm 1\}^k.
\end{equation*}
Hence, denoting $v_{n+i} = v_i^*$ the velocities of added particles, and introducing the measure
\begin{equation} \label{def:measure}
\dd \nu_{[t]}(\Psi_n) \doteq p_\mu^{|\ul^*_k|} \dd \ut_k \dd \uom_k \dd \uv_k^* \prod_{i=1}^k s_i \langle \omega_{i}, v_{n+i} - v_{m_i}^{[t_i^+]} \rangle_+ ,
\end{equation}
the pseudo-trajectory formulation rewrites as
\begin{equation} \label{eq:int_meas_pseudo_traj}
F_n^\varepsilon(t) =  \int  F^\varepsilon\left(0, \uz^{[0]}_{\Psi_n}\right)\ \dd \nu_{[t]}(\Psi_n),
\end{equation}
where $\uz^{[0]}_{\Psi_n}$ denotes the initial configuration deduced from the pseudo-trajectory, according to the construction given in Section~\ref{sec:pseudo-traj} (the initial configuration~$\uz^{[0]}_{\Psi_n}$ also includes the \emph{tags} of the corresponding particles). Note finally that the sums over $k$ and $\uchi_k$ in the pseudo-trajectory formulation result from the domain of integration of the measure~$\nu_{[t]}$, and so remain implicit.

The formula giving an expression for the expectation of the empirical measure according to the correlation functions~\eqref{eq:empmes}, may be generalized~\cite[Proposition~3.3.1]{2023grandev} to observables~$H_n(\uz_n^{[0,t]})$ depending on the whole trajectory on~$[0,t]$, as
\begin{align*}
& \E\left[ \sum_{1 \leq i_k \neq i_j \leq \mathcal{N}} H_n(Z_{i_1}^{[0,t]}, L_{i_1}, \dots, Z_{i_n}^{[0,t]}, L_{i_n})  \right] \\
&\hspace{20mm} = \mu^n \sum_{\ul_n  \in \Lambda_n} p_\mu^{|\ul_n|}  \int_{\mathcal{D}^\varepsilon_n} \dd \uz_n \int F^\varepsilon(0, \uz^{[0]}_{\Psi_n})  H_n(\uz_{\Psi_n}^{[0,t]}) \dd \nu_{[t]}(\Psi_n),
\end{align*}
where $H_n(\uz_{\Psi_n}^{[0,t]})$ is projecting the whole trajectory on the trajectories of the~$n$~first particles. The object that we will study, for tensorized observables, is hence the \emph{$H$-weighted correlation function}
\begin{equation} \label{def:Hweighted_F}
F_n^\varepsilon[H](t, \uz_n, \ul_n) \doteq \int F^\varepsilon(0, \uz^{[0]}_{\Psi_n})  H^{\otimes n}(\uz_{\Psi_n}^{[0,t]})\ \dd \nu_{[t]}(\Psi_n).
\end{equation}
Note that for~$H\equiv 1$, we retrieve the correlation functions. Moreover, the cumulant generating function~\eqref{def:cumulant_generating_function} also generalizes to obervables depending on the whole trajectory~$H_n(\uz_n^{[0,t]})$, as follows
\begin{align} 
 \mathfrak{G}_\varepsilon^{[0,t]}[H]& \doteq \log \E\left[ \exp\left( \sum_{i = 1}^\mN H(Z_i^{[0,t]}, L_i) \right) \right] \nonumber \\[1em]
  & = \sum_{p \geqslant 1}  \frac{1}{p!} \sum_{\ul_p \in \Lambda_p } \lambda^{|\ul_p|} \mu^{p - |\ul_p|} \int_{\mD^p} f^\varepsilon_{p}\left[ e^H - 1\right](t, \ul_{p}). \label{def:cumulant_generating_function3}
\end{align}
Indeed, the computation of Section~\ref{sec:exp_cum_gen_fun} are identical for the expansion of this generalized version of the cumulant generating function. The following section gives a formula for the $H$-weighted cumulants~$f^\varepsilon_{p}\left[ H \right]$.

\subsection{Expansion and equation on the cumulants} \label{sec:cum_expansion}

We perform from now on the expansion on the cumulants, that may be found in~\cite[Chapter~3, \emph{Tree expansions of the hard-sphere dynamics}]{2023grandev}, for the general indistinguishable case.
The purpose is to start with the pseudo-trajectory equation~\eqref{eq:int_meas_pseudo_traj} on the correlation functions, and to identify the cumulants in an expansion similar to the one characterizing them~(see the inversion formula in Proposition~\ref{prop:inversion_formula}).
This expansion is a first illustration of the link between cumulants and dynamical interactions.

We start by splitting the pseudo-trajectories into \emph{non-interacting aggregates}, to factorize the pseudo-trajectory measure. 
At a fixed pseudo-trajectory $\Psi_n = (k, \uchi_k, \ut_k, \uom_k, \uv_k) $, for a set of particles $A \subset \lbr 1, n\rbr$, we will denote 
\begin{equation} \label{eq:P(m)}
P(A) \subset A \cup \lbr n+1, n+k \rbr
\end{equation}
 the set of all particles contained in the pseudo-trajectory stemming from $A$ (depending on the history~$\uchi_k$). Then, if $\ind^{[0,t]}_{p \not\sim q}$ is the indicator that particles~$p$ and~$q$ never collided on~$[0,t]$ during the pseudo-trajectory~$(\uz_n,\Psi_n)$, we denote for $A, A'$ two subsets of~$ \lbr 1, n \rbr$, 
\begin{equation*} \ind_{A \not\sim A'} \doteq \prod_{p \in P(A)} \prod_{q \in P(A')} \ind^{[0,t]}_{p \not\sim q}
\end{equation*}
the indicator that no particle stemming from~$A$ encountered any particle stemming from $A'$, during the whole pseudo-dynamics on~$[0,t]$. Under the condition $A \not\sim A'$, the pseudo-trajectories are independent and the measure and observable factorize as
\begin{equation*}
\dd \nu_{[t]}(\Psi_{A \cup A'}) H^{\otimes A \cup A'}(\uz_{\Psi_{A \cup A'}}^{[0,t]}) = \dd \nu_{[t]}(\Psi_{A}) H^{\otimes A}(\uz_{\Psi_{A}}^{[0,t]}) \times \dd \nu_{[t]}(\Psi_{A'}) H^{\otimes A'}(\uz_{\Psi_{ A'}}^{[0,t]}).
\end{equation*}
We keep this factorized writing for the moment, to remember that the pseudo-trajectories are defined and constructed independently on each aggregate. In practice, they behave as if the particles from two different aggregates could not interact: they overlap one another when they meet.

On the other hand, we denote the \emph{aggregating condition}~$\cls(A)$ the indicator that all the particles in the aggregate~$A$ are connected through collisions in their pseudo-trajectories, i.e. that there exists a path $(i_1, \dots, i_{|A|}) \in \lbr 1, |A| \rbr^{|A|}$ such that
\begin{equation} \label{def:cls} \{ i_1, \dots, i_{|A|} \} = \lbr 1, |A| \rbr \ \und \ \forall j \in \lbr 1,\ |A| -1 \rbr,\  \ind_{A_{i_j} \sim A_{i_{j+1}}}^{[0,t]} = 1.
\end{equation}
This eventually leads to the following conditioning according to the partition~$\kappa \in \mP_n$ into aggregates:
\begin{align*}
F_n^\varepsilon[H](t) &  
=  \sum_{\kappa \in \mP_n} \int \left[ \prod_{i=1}^{|\kappa|}\cls(\kappa_i)\ \dd \nu_{[t]}(\Psi_{\kappa_i}) H^{\otimes \kappa_i}(\uz_{\Psi_{\kappa_i}}^{[0,t]}) \right]\left( \prod_{1 \leqslant i < j \leqslant |\kappa|} \ind_{\kappa_{i} \not\sim \kappa_j}\right) F^\varepsilon(0, \uz^{[0]}_{\Psi_n}).
\end{align*}
Henceforth, we denote the measure weighted by the observable
\begin{equation} \label{def:measure_weighted}
\dd \nu^{[H]}_{[0,t]}(\Psi_n) \doteq H^{\otimes n}(\uz_{\Psi_n}^{[0,t]}) \dd \nu_{[t]}(\Psi_n).
\end{equation}
Now, still in the idea of factorizing, we want to decorrelate the condition of exclusion between the aggregates~$\ind_{\kappa_i \not\sim \kappa_j}$, using the cumulants of the exclusion (Definition~\ref{def:cumulants}) thanks to which we write
\begin{equation*}
\prod_{1 \leqslant i < j \leqslant |\kappa|} \ind_{\kappa_{i} \not\sim \kappa_j} = \sum_{\rho \in \mP_{|\kappa|}} \phi_{[\rho]}(\kappa_1, \dots, \kappa_{|\kappa|}) 
\end{equation*}
recalling the notation
\begin{equation*}
\phi_{[\rho]}(\kappa_1, \dots, \kappa_{|\kappa|}) =  \prod_{i = 1}^{|\rho|} \phi_{|\rho_i|}(\ukap_{\rho_i}).
\end{equation*}
 Finally, we still need to expand the initial distribution~$F^\varepsilon(0, \uz^{[0]}_{\Psi_n})$ to completely factorize the formula, to make appear a cumulant expansion. We need to expand it in a coarser way than the partition~$\rho$, to be compatible with the product~$\phi_{[\rho]}$ above, so that we consider the \emph{cluster cumulants}, defined by the expansion
\begin{equation} \label{eq:initial_cluster_cumulants_1}
F^\varepsilon(0, \uz^{[0]}_{\Psi_{\rho_1}}, \dots, \uz^{[0]}_{\Psi_{\rho_{|\rho|}}}) = \sum_{\sigma \in \mP_{|\rho|}} f^{\varepsilon, \rho}_{[\sigma]}(0, \uz^{[0]}_{\Psi_{\rho_1}}, \dots, \uz^{[0]}_{\Psi_{\rho_{|\rho|}}})
\end{equation}
due once again to the inversion formula (Proposition~\ref{prop:inversion_formula}). This expansion is made in clusters, and so depends on the partition of the aggregates~$\rho \in \mP_{|\kappa|}$.
The  following notation $\sigma \vartriangleleft \rho \vartriangleleft \kappa$ means that the partition~$\sigma$ is coarser than~$\rho$, which is coarser then~$\kappa$, and in the end we write
\begin{align*}
F_n^\varepsilon[H](t) &  
=  \sum_{\kappa \in \mP_n} \int \left[ \prod_{i=1}^{|\kappa|}\cls(\kappa_i)\ \dd \nu^{[H]}_{[0,t]}(\Psi_{\kappa_i}) \right] \sum_{\sigma \vartriangleleft \rho \vartriangleleft \kappa} \phi_{[\rho]} f^{\varepsilon,\rho}_{[\sigma]}(0).
\end{align*}
To identify with a cumulant expansion, the last step is to use the same trick as in the proof of the inversion formula to invert the partitions (Proposition~\ref{prop:inversion_formula}), starting by considering a partition~$\sigma$ of~$\lbr 1, n\rbr$, then considering a partition~$\rho$ \emph{on each subset}~$\sigma_j$. This is made possible by to the compatibility condition
\begin{equation*}
f^{\varepsilon,\rho}_{\sigma_j} = f^{\varepsilon,\rho_{|\sigma_j}}_{\sigma_j},
\end{equation*}
and thanks to the factorizing identity
\begin{align*}
\phi_{[\rho]} f^{\varepsilon,\rho}_{[\sigma]} = \prod_{i = 1}^{|\rho|} \phi_{\rho_i} \prod_{j = 1}^{|\sigma|} f^{\varepsilon, \rho}_{\sigma_j} = \prod_{j=1}^{|\sigma|} \left[ f^{\varepsilon, \rho}_{\sigma_j} \prod_{i \in \sigma_j} \phi_{\rho_i} \right].
\end{align*}
Since the initial cluster cumulants~$f_n^{\varepsilon, \rho}$ do not depend on the subdivision~$\kappa$, we can also invert~$\kappa$ and~$\rho$ to finally get the cumulant expansion
\begin{equation*}  
F^\varepsilon_n[H](t) = \sum_{\sigma \in \mP_n} \prod_{j=1}^{|\sigma|} \int	 \sum_{ \rho \in \mP_{\sigma_j}}  f^{\varepsilon,\rho}_{\sigma_j}(0) \prod_{j=1}^{|\rho|} \left[ \sum_{\kappa \in \mP_{\rho_j}}\phi_{\rho_j} \prod_{i=1}^{|\kappa|} \cls(\kappa_i)\ \dd \nu^{[H]}_{[0,t]}(\Psi_{\kappa_i}) \right],
\end{equation*}
whence the following lemma by injectivity.

\begin{lemm} \label{lem:cum_pseudo_traj} The cumulants satisfy the following integral equation, in terms of aggregates and clusters of interaction
\begin{equation}   \label{eq:cum_pseudo_traj}
f_n^\varepsilon[H](t, \uz_n,\ul_n) = \sum_{ \rho \in\mP_n} \int f^{ \varepsilon,\rho}_{n}(0, \uz^{[0]}_{\Psi_{\rho_1}}, \dots, \uz^{[0]}_{\Psi_{\rho_{|\rho|}}}) \prod_{j=1}^{|\rho|} \left[ \sum_{\kappa \in \mP_{\rho_j}} \phi_{\rho_j} \prod_{i=1}^{|\kappa|} \cls(\kappa_i)\ \dd \nu^{[H]}_{[0,t]}(\Psi_{\kappa_i}) \right].
\end{equation}
\end{lemm}
Note that the indicators~$\cls(\kappa_i)$ and the cumulants~$\phi_{[\rho]}$ depend on the whole pseudo-trajectory, and hence depend on the time~$t$.

\subsection{Initial cluster cumulants} \label{sec:initial_cluster_cum}
This section is dedicated to proving Lemma~\ref{lem:clust_cum} below, yielding an explicit formulation of the initial cluster cumulants implicitly defined in~\eqref{eq:initial_cluster_cumulants_1}, and useful to find bounds in Section~\ref{sec:bounds}.
We denote
\begin{equation*}\aleph_{|\rho|} \doteq \{\uz^{[0]}_{\Psi_{\rho_i}} \ , \ 1 \leq i \leq |\rho| \} \ \ \und \ \ \Gamma_p \doteq \{ z_i^* \ , \ 1 \leq i \leq p \}
\end{equation*}
the set of clusters and the set of integrated particles, $N \doteq |\Psi_{\rho_1}| + \dots + |\Psi_{\rho_{|\rho|}}|$ the total  number of particles contained in the clusters, and 
\begin{equation*}
\left[ \ind_{\mX_{(\cdot)}} \right]^{\otimes \aleph_{|\rho|}} \doteq \prod_{\uz \in \aleph_{|\rho|}} \ind_{\mX_{|\uz|}	}(\uz)
\end{equation*}
the indicator that \emph{within} each cluster, the particles exclude themselves.
Finally, we denote as before $\ind_{x \sim y} = 1 - \ind_{x \not\sim y}$ the exclusion condition between two subsets of particles~$x$ and~$y$ (potentially clusters), and for $S \subset \aleph_{|\rho|}$, $k\leq p$ and $l \leq q$, we introduce the \emph{cumulants of the cluster exclusion}
\begin{align*}
\phi_{S, k+l} \doteq \sum_{G \in \mC_{S \cup \lbr k+l \rbr}} \prod_{\{x,y\} \in E_G} (-\ind_{x\sim y}),
\end{align*}
where~$\mC_A$ stands for the connected graphs on a set~$A$.

\begin{lemm}[Explicit formula for the initial cluster cumulants] \label{lem:clust_cum}
For any~$\lambda, \mu > 0$, with the notation above, the initial cluster cumulants can be written as
\begin{align*}
 f_n^{\varepsilon,\rho}(0, \uz_{\Psi_{\rho_1}}, \dots, \uz_{\Psi_{\rho_{|\rho|}}}) = \left[ \ind_{\mX_{(\cdot)}} \right]^{\otimes \aleph_{|\rho|}} M_\beta^{\otimes N} \varphi_0^{\otimes \ul_{N}}  \sum_{p, q \geqslant 0} \frac{\lambda^{p} \mu^{q}}{p!q!}  \int \varphi_0^{\otimes p} M_\beta^{\otimes p + q} \phi_{\aleph_{|\rho|}, p+q}.
\end{align*}
\end{lemm}

\begin{prof}
To obtain an explicit formulation for the initial cluster cumulants, we will identify them using the following characterization 
\begin{equation} \label{eq:clust_cum}
F^\varepsilon(0, \uz^{[0]}_{\Psi_{\rho_1}}, \dots, \uz^{[0]}_{\Psi_{\rho_{|\rho|}}}) = \sum_{\sigma \in \mP_{|\rho|}} f^{\varepsilon, \rho}_{[\sigma]}(\uz^{[0]}_{\Psi_{\rho_1}}, \dots, \uz^{[0]}_{\Psi_{\rho_{|\rho|}}}).
\end{equation}
Using the exchangeability of particles, we have
\begin{align*}
F^\varepsilon(0, \uz^{[0]}_{\Psi_{\rho_1}}, \dots, \uz^{[0]}_{\Psi_{\rho_{|\rho|}}}) & = \frac{1}{\mZ_{\mu}} \sum_{p\geqslant 0}  \sum_{\ul_{p}^* \in \Lambda_p} \frac{\lambda^{|\ul_{p}^* |} \mu^{p - |\ul_{p}^*|}}{p!}  \int \varphi_0^{\otimes \tilde{\ul}_{N+p}} M_\beta^{\otimes N+p} \ind_{\mX^\varepsilon_{N+p}} \dd \uz_{p}^* \\
& =   \frac{1}{\mZ_{\mu}} \sum_{p\geqslant 0}  \sum_{q = 1}^p \frac{\lambda^{q} \mu^{p - q}}{q!(p-q)!} \int \varphi_0^{\otimes \ul_N } \varphi_0^{\otimes q} M_\beta^{\otimes N+p} \ind_{\mX^\varepsilon_{N+p}} \dd \uz_{p}^* 
\end{align*}
To retrieve a cumulant formulation like~\eqref{eq:clust_cum} above, we will expand into cumulants the exclusion condition---which is the only obstruction to independence. To get rid of the partition function~$\mZ_\mu$ in the process, we will not expand the condition according to the integrating variables~$\uz_p^*$, to make appear the partition function at the numerator. More precisely, decomposing the exclusion condition within each cluster, we write
\begin{align*}
\ind_{\mX_{N+p}}(\uz^{[0]}_{\Psi_{\rho_1}}, \dots, \uz^{[0]}_{\Psi_{\rho_{|\rho|}}}, \uz_p^*) & = \left[ \ind_{\mX_{(\cdot)}}\right]^{\otimes \aleph_{|\rho|}}  \times \ind_{\hat{\mX}_{|\rho|+p}}(\uz^{[0]}_{\Psi_{\rho_1}}, \dots, \uz^{[0]}_{\Psi_{\rho_{|\rho|}}}, \uz_p^*),
\end{align*} 
with
\begin{align*}
\ind_{\hat{\mX}_{|\rho|+p}}(\uz^{[0]}_{\Psi_{\rho_1}}, \dots, \uz^{[0]}_{\Psi_{\rho_{|\rho|}}}, \uz_p^*) & = \prod_{x,y \in \aleph_{|\rho|} } \ind_{x \not\sim y} \prod_{x \in \aleph_{|\rho|} ,y \in \Gamma_p } \ind_{x \not\sim y}
 \prod_{x,y \in \Gamma_p} \ind_{x \not\sim y} \\
 & = \sum_{G \in \mG_{(\aleph_{|\rho|} \cup \Gamma_p)}} \prod_{\{x,y\} \in E_G} (- \ind_{x\sim y}).
\end{align*} 
We want now to split the sum above according to the partition induced by the connected components of the graph~$G$, as for the standard cumulants of the exclusion~\cite[Appendix~A]{fou26Rayleigh}. Nevertheless, since we want to make appear the partition function, we will compute a coarser splitting. Indeed, we will merely isolate the part of the graph which is not connected to any cluster of~$\aleph_{|\rho|}$, i.e. the subset~$B \subset \Gamma_{p}$ of integrated particles that are only connected between themselves. We still decompose the rest of the graph into its connected components, yielding a partition of~$\aleph_{|\rho|} \cup B^c$ in which each subset has to contain a cluster (otherwise it would have been kept in~$B$), which we denote~$\sigma \in \tilde{\mP}[\aleph_{|\rho|} ,B^c]$. Hence, one can write
\begin{align*}
\ind_{\hat{\mX}_{|\rho|+p}}(\uz^{[0]}_{\Psi_{\rho_1}}, \dots, \uz^{[0]}_{\Psi_{\rho_{|\rho|}}}, \uz_p^*) & = \sum_{B \subset \Gamma_p}  \sum_{G \in \mG_{B}} \prod_{\{x,y\} \in E_G} (- \ind_{x\sim y}) \sum_{\sigma \in \tilde{\mP}[\aleph_{|\rho|} ,B^c] } \prod_{i=1}^{|\sigma|} \sum_{G \in \mC_{\sigma_i}} \prod_{\{x,y\} \in E_G} (- \ind_{x\sim y}) \\
& = \sum_{B \subset \Gamma_p}  \ind_{\mX_{|B|}}(\uz^*_B) \sum_{\sigma \in \tilde{\mP}[\aleph_{|\rho|} ,B^c]} \prod_{i=1}^{|\sigma|} \phi_{\sigma_i},
\end{align*} 
gathering the first sum over~$\mG_B$ by inverting the usual expansion, and recognizing cumulants of the exclusion~\cite[Appendix~A]{fou26Rayleigh}.
In the end, a partition~$\sigma \in \tilde{\mP}[\aleph_{|\rho|} ,B^c]$ is a partition of $\aleph_{|\rho|}$ on which some particles from~$B^c$ are added, either tagged by~$\ul^*_p$ or not. We will use the symmetry of particles to reduce these added particles to their numbers of 1-tags $k_1, \dots, k_{|\sigma|} \geq 0$ and of 0-tags $l_1, \dots, l_{|\sigma|} \geq 0$, and similarly for the particles in~$B$, in numbers~$k_0, l_0 \geq 0$. Omitting the indicator~$\left[ \ind_{\mX_{(\cdot)}} \right]^{\otimes \aleph_{|\rho|}}$ to gain lisibility (it factorizes and distributes very easily), we thus write
\begin{align*}
&F^\varepsilon\left[0, \uz^{[0]}_{\Psi_{\rho_1}}, \dots, \uz^{[0]}_{\Psi_{\rho_{|\rho|}}}\right] \\
& = \frac{\varphi_0^{\otimes \ul_N }}{\mZ_{\mu}} \sum_{p\geqslant 0}  \sum_{q = 1}^p \frac{\lambda^{q} \mu^{p - q}}{q!(p-q)!} \sum_{\sigma \in \mP_{\aleph_{|\rho|}}} \sum_{\substack{k_0 + \dots + k_{|\sigma|} = q \\ l_0 + \dots + l_{|\sigma|} = p-q}} \frac{q!}{k_0!\dots k_{|\sigma|}!} \frac{(p-q)!}{l_0!\dots l_{|\sigma|}!} \int  \varphi_0^{\otimes q} M_\beta^{\otimes N+p}  \ind_{\mX^\varepsilon_{k_0+l_0}} \prod_{i=1}^{|\sigma|} \phi_{\sigma_i, k_i+l_i}\\
& = \frac{\varphi_0^{\otimes \ul_N} M_\beta^{\otimes N}}{\mZ_{\mu}} \sum_{k_0, l_0 \geqslant 0} \frac{\lambda^{k_0} \mu^{l_0}}{k_0! l_0!} \int  \varphi_0^{\otimes k_0} M_\beta^{\otimes l_0 + k_0} \ind_{\mX^\varepsilon_{k_0+l_0}} \times \sum_{\sigma \in \mP_{\aleph_{|\rho|}}}   \prod_{i=1}^{|\sigma|}  \sum_{k_i, l_i \geqslant 0} \frac{\lambda^{k_i} \mu^{l_i}}{k_i!l_i!}  \int  \varphi_0^{\otimes k_i} M_\beta^{\otimes l_i + k_i} \phi_{\sigma_i, k_i+l_i} \\[1em]
& = \sum_{\sigma \in \mP_{\aleph_{|\rho|}}}  \prod_{i=1}^{|\sigma|} M_\beta^{\otimes \sigma_i} \varphi_0^{\otimes \ul_{\sigma_i}}  \sum_{k_i, l_i \geqslant 0} \frac{\lambda^{k_i} \mu^{l_i}}{k_i!l_i!}  \int  \dd \uz_{l_i+k_i}^* \varphi_0^{\otimes k_i} M_\beta^{\otimes l_i + k_i} \phi_{\sigma_i, k_i+l_i},
\end{align*}
recognizing the partition function~$\mZ_\mu$ at the numerator. The lemma is thus proved by identification once added back the indicator~$ \ind_{\mX_{(\cdot)}} ^{\otimes \aleph_{|\rho|}}$. \CQFD
\end{prof}

\subsection{Dynamics trees} \label{sec:dyn_trees}
Starting from the formulation~\eqref{eq:cum_pseudo_traj} of the cumulants in terms of pseudo-trajectories, we will rewrite this integral to make appear the \emph{correlations}' history between the particles. Indeed, the cumulants encode the small defects of particles' independence, which correspond to fortuitous encounters happening between the pseudo-trajectories of the $n$~considered particles. In fact, we prove that asymptotically, the $n$-th cumulant~$f_n^\varepsilon$ is supported over trajectories implying exactly $(n-1)$~of these fortuitous encounters, thus connecting all the pseudo-trajectories. This way, each of these interactions has to be clustering, to eventually form a minimally connected graph.

We will control these fortuitous interactions between pseudo-trajectories in the integrated form
\begin{equation}   \label{eq:In_rho}
\sum_{\rho \in \mP_n} \int \dd \uz_n \int f^{ \varepsilon,\rho}_{n}(0, \uz^{[0]}_{\Psi_{\rho_1}}, \dots, \uz^{[0]}_{\Psi_{\rho_{|\rho|}}}) \prod_{j=1}^{|\rho|} \left[ \sum_{\kappa \in \mP_{\rho_j}} \phi_{\rho_j} \prod_{i=1}^{|\kappa|} \cls(\kappa_i)\ \dd \nu^{[H]}_{[0,t]}(\Psi_{\kappa_i}) \right],
\end{equation}
which is the integral over~$\uz_n$ of~$f_n^\varepsilon[H](t, \uz_n,\ul_n)$ in its form~\eqref{eq:cum_pseudo_traj}. The integral over the endstate~$\uz_n$ allow to act on the arrival positions~$\ux_n$ of the trajectories to make the fortuitous interactions appear or disappear. In the end, we will perform a change of these variables~$\ux_n$ that harnesses these interactions.

To that end, these fortuitous encounters will be recorded in \emph{dynamics trees}, allowing us to express the conditions of their existence. They stem from two different sources: first, \emph{recollisions} appear within the aggregates $(\kappa_i)$ to make them connected; and on the other hand, on each cluster~$\rho_j$ of aggregates, the extended cumulants of their mutual exclusion $\phi_{\rho_j}$ make appear \emph{overlapping} conditions between the aggregates, recalling that when two different aggregates meet, they overlap one another without interacting.
This formulation in dynamics trees will be the first step to prove bounds on the cumulants, and to compute their limit (Sections~\ref{sec:bounds} and~\ref{sec:conv_cumulants}).

\paragraph{Clustering recollisions}

As explained above, within a typical aggregate~$\kappa_i \subset \lbr 1, n \rbr$, denoting $A \doteq \kappa_i$ and $a \doteq |A|$, each pseudo-trajectory~$\Psi_A$ contains recollisions~\eqref{def:cls} that connect the aggregate altogether, according to the condition~$\cls(A)$. 

\begin{deaf}[Clustering recollisions] In a given pseudo-trajectory, looked at backwards in time starting from time~$t$, we call \emph{clustering recollision} the first recollision connecting the pseudo-trajectories stemming from two different particles. In an aggregate~$A$, we order these clustering recollisions backwards in time from $1$ to~$(a-1)$, which defines a \emph{recollision tree}~$T^\col(\uz_{\Psi_A}^{[0,t]}) \in \mT_A^\prec$ whose ordered edges record these recollisions: two particles $m,m' \in A \subset \lbr 1, n \rbr$ are connected in this tree if and only if a clustering recollision happens between the pseudo-trajectories stemming from~$m$ and~$m'$. Eventually, we define
\begin{equation*}
(\tau^\col_i, \omega^\col_i)_{1\leqslant i \leqslant a-1} \in ([0,t]\times \Sf^{d-1})^{a-1}
\end{equation*}
the ordered collision times and angles at which the recollision occurs.
\end{deaf}
For simplicity, we denote $T^\col_A \doteq T^\col(\uz_{\Psi_A}^{[0,t]})$. We now partition according to the associated tree
\begin{equation*}
\cls(A) = \sum_{\Upsilon \in \mT_{A}^\prec } \ind_{T^\col_A = \Upsilon} .
\end{equation*}
Under the condition~$T^\col_A = \Upsilon$, we will perform a change of variables that harnesses each of the clustering recollision conditions. To do so, we need to know \emph{which particles}, from the pseudo-trajectories stemming from $m, m' \in A$, collided to connect them.
Recall that~$P(m)$ denotes the set of all particles added in the pseudo-trajectory stemming from $m$. Then, for every edge $e = (m,m') \in E_\Upsilon$, we define $\coll_e(A)$ the pair of particles  $(p_e, q_e) \in P(m)\times P(m')$ that collided on~$(0,t)	$ to connect the pseudo-trajectories of~$m$ and~$m'$, when creating the edge $e$. Eventually, we also partition according to which particle collided with which, for each clustering recollision:
\begin{equation*}
\cls(A) = \sum_{\Upsilon \in \mT_{A}^\prec }  \ind_{T^\col_A = \Upsilon}  \sum_{\substack{p_e \in P(m), q_e \in P(m') \\ e=(m,m')\in  \E_\Upsilon}} \ \ \prod_{e \in  E_\Upsilon} \ind_{(p_e, q_e) = \coll_e(A)}.
\end{equation*} 
We introduce the following notation for the indicator that the clustering particles are those we want
\begin{equation*}
\mathfrak{S}_{(p_e, q_e)_{\Upsilon}}[A]  \doteq \prod_{e \in E_\Upsilon} \ind_{(p_e, q_e) = \coll_e(A)}.
\end{equation*}
We can now harness the successive clustering conditions, which correspond to the ordered edges of the tree~$T^\col_A$, enumerated~$(e_1, \dots, e_{a-1}) = \bigl( (m_1, m'_1), \dots, (m_{a-1}, m'_{a-1}) \bigr)$ backwards in time. The relevant variables that we want to make appear are the times and angles of the collisions, when the current integration variables are the arrival positions $\ux_A$ at time $t$.

As we are interested in the relative positions of the particles implied in the clustering recollisions, we start with the following change of variable of Jacobian~1, denoting $\overline{x}_A$ the \emph{barycenter of the positions}~$\ux_A$,
\begin{equation*}
\ux_A \mapsto \left( \overline{x}_A,  \left( \vec{x}_e \doteq x_{p_e} - x_{q_e} \right)_{e \in E_{\Upsilon}} \right),
\end{equation*}
where we take the convention~$p_e > q_e$, to be consistent in the following with the measure~$\nu$.
 Thus, once recursively determined the clustering conditions associated to all the edges preceding an edge~$e$ such that $\coll_e(A) = (p_e, q_e)$, the condition associated with this edge states that at time~$\tau_{e}^\col$, the relative position $x_{p_e}(\tau_e^\col) - x_{q_e}(\tau_e^\col)$ must belong to the sphere of radius $\varepsilon$, defining the recollision angle~$\omega_e^\col \in \Sf^{d-1}$. Hence, denoting~$\tau_\dd$ the time of the closest deflection that underwent $p_e$ or $q_e$ between~$t$ and $\tau_e^\col$, and $\Delta x_\dd$ the relative distance travelled between $t$ and $\tau_\dd$, the relative positions at time~$t$ might be retrieved from
\begin{equation} \label{eq:pre-changement-de-var}
x_{p_e} - x_{q_e} = \varepsilon \omega_e^\col - (\tau_e^\col - \tau_\dd)(v_{p_e}({\tau_e^\col}^+) - v_{q_e}({\tau_e^\col}^+)) - \Delta x_\dd. 
\end{equation}
The parameters $(\tau_e^\col, \omega_e^\col)$ are unique because we consider the very first time (backwards) at which $p_e$ and $q_e$ meet, so that the mapping~\eqref{eq:pre-changement-de-var} above is invertible: we consider the local change of variable 
\begin{equation*} 
x_{p_e} - x_{q_e} = \vec{x}_{e} \mapsto (\tau_e^\col, \omega_e^\col) \in [0, \tau_{e'}^\col]\times \Sf^{d-1},
\end{equation*} 
where $e'$ is the edge ordered right before~$e$. Since~$\Delta x_\dd$ is a piecewise affine function recording the part of the dynamics that does not depend on~$\tau_e^\col$ (but depends on the previous recursively determined clustering conditions, the dependency in~$\tau_e^\col$ being explicit in~\eqref{eq:pre-changement-de-var}), the Jacobian of this change of variables is
\begin{equation}
\label{eq:jacob}
 \varepsilon^{d-1} \bigl \langle v_{p_e}({\tau_e^\col}^+ )  - v_{q_e}({\tau_e^\col}^+),  \omega_e^\col \bigr\rangle_+ = \mu^{-1} \bigl\langle v_{p_e}({\tau_e^\col}^+ )  - v_{q_e}({\tau_e^\col}^+),  \omega_e^\col \bigr\rangle_+.
\end{equation}
Iterating this computation, we change the variables~$(\vec{x}_e)_{e \in E_\Upsilon}$ to the ordered times~$\utau^\col_{a-1}$ and angles~$\uom^\col_{a-1}$ of clustering recollisions (also indexed by the ordered edges of~$E_\Upsilon$), finally leading to
\begin{align*}
& \int \dd \ux_A \Bigl[\ind_{T^\col_A = \Upsilon}\mathfrak{S}_{(p_e, q_e)_{\Upsilon}}[A] \Bigr] \dd \nu^{[H]}_{[0,t]}(\Psi_{A}) \\
& \hspace{8mm} = \int \frac{\dd \overline{x}_A}{\mu^{a-1}} \dd \uom_{a-1}^\col \dd \utau_{a-1}^\col \Bigl[\ind_{T^\col_A = \Upsilon}\mathfrak{S}_{(p_e, q_e)_{\Upsilon}}[A] \Bigr] \dd \nu^{[H]}_{[0,t]}(\Psi_{A})  \prod_{e \in E_{\Upsilon}} \langle v_{p_e}({\tau_e^\col}^+ )  - v_{q_e}({\tau_e^\col}^+),  \omega_e^\col \rangle_+.
\end{align*}

\paragraph{Clustering overlaps}
Now for a general cluster $R = \rho_j$, of size $r = |R|$, we expand the cumulant
\begin{align*}
\phi_{R}   =  \sum_{G \in \mC_{|\kappa|}}  \prod_{\{i,j\} \in E_G} (-\ind_{\kappa_i \sim \kappa_j})
\end{align*}
making appear \emph{overlap conditions} $\ind_{\kappa_i \sim \kappa_j}$. At fixed pseudo-trajectory parameters~$\Psi_R$, changing the barycenter~$\overline{x}_{\kappa_i}$ of an aggregate moves rigidly the whole aggregate pseudo-trajectory, since by construction the trajectories are defined independently on each aggregate. That way, we can act on the barycenters to make them match the overlaps given by a given graph~$G \in \mC_{|\kappa|}$. In the same fashion as for the recollision trees, we define the following objects. 

\begin{deaf}[Clustering overlaps] In a given pseudo-trajectory, we call~\emph{clustering overlap} the first overlap (backwards in time) appearing between independently defined aggregates, connecting them according to the overlap condition discussed above. As for recollisions, clustering overlaps define an \emph{overlap tree}~$T^\ov(\uz_{\Psi_R^{[0,t]}}) \in \mT_{|\kappa|}^\prec$, and for each edge~$e = \{m,m'\}$ associated to an overlap, we define
\begin{equation*}
\tau^\ov_e = \sup\left\{ 0 \leq s \leq t\ , \ \ind_{\kappa_m \sim \kappa_{m'}}^{[s,t]} = 1  \right\} 
\end{equation*}
its overlap time (arbitrarily taking the latest possible, to make them consistent with recollisions when constructing pseudo-trajectories backwards in time), along with the associated overlap angle~$\omega^\ov_e$.
\end{deaf}
 As for the recollision trees, we denote~$T^\ov_R \doteq T^\ov(\uz_{\Psi_R^{[0,t]}})$ the overlap tree. Now, we partition the expansion of the cumulants according to the associated tree, treating separately when the graph appearing in the sum is the overlap tree, and when it presents some additional cycles. Denoting $\mC_{|\kappa|}(\Upsilon)$ the connected graphs containing~$\Upsilon$, we have
\begin{equation}\def\arraystretch{2} \label{def:tree_cycle}
\begin{array}{rccc}
\phi_{R} = & \sum_{\Upsilon \in\mT^\prec_{|\kappa|}} \ind_{T^\ov_R = \Upsilon} \prod_{\{i,j\} \in E_\Upsilon} (-\ind_{\kappa_i \sim \kappa_j}) & + &\sum_{\Upsilon \in\mT^\prec_{|\kappa|}} \ind_{T^\ov_R = \Upsilon} \sum_{G \in \mC_{|\kappa|}(\Upsilon)\setminus \{\Upsilon\}}  \prod_{\{i,j\} \in E_G} (-\ind_{\kappa_i \sim \kappa_j}) \\
\doteq  &  \phi_R^{[\tree]} & + & \phi_R^\cycle .
\end{array}
\end{equation}
Focusing first on the tree part, one can write
\begin{align*}
\sum_{\kappa \in \mP_R}\phi^{[\tree]}_{R} \prod_{i=1}^{|\kappa|} \cls(\kappa_i) \dd \nu^{[H]}_{[0,t]}(\Psi_{\kappa_i})  = (-1)^{|\kappa|-1} \sum_{\Upsilon \in \mT^\prec_{|\kappa|}} \ind_{T^\ov_R = \Upsilon} \prod_{i = 1}^{|\kappa|} \cls(\kappa_i) \dd \nu^{[H]}_{[0,t]}(\Psi_{\kappa_i}),
\end{align*}
the overlap conditions being included in the overlap tree matching condition. At fixed overlap tree $\Upsilon = T^{\ov}_R$, we act now on the aggregate barycenters that appeared during the change of variable harnessing the recollisions. Recall that moving these barycenters moves rigidly the aggregate pseudo-trajectories. Like for the clustering recollisions, we hence focus on the relative positions of these barycenters through the change of variable, of Jacobian~1,
\begin{equation*}
(\overline{x}_{\kappa_i})_{i \leqslant |\kappa|} \mapsto \left( \overline{x}_R,  \left( \vec{x}_e \doteq x_{p_e} - x_{q_e} \right)_{e \in E_{\Upsilon} } \right),
\end{equation*}
conditionning on the overlapping particles $(p_e, q_e)$ associated to the edge~$e$, as it was done for recollisions.

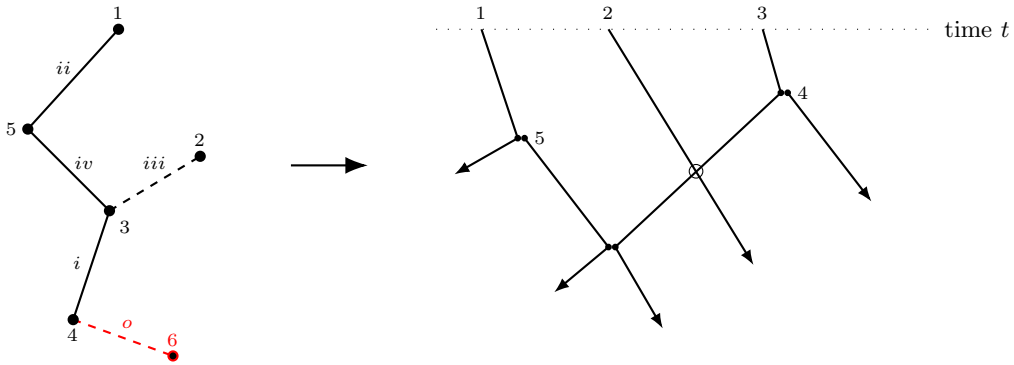
\begin{figure}[h!]\centering 
\begin{tikzpicture}
\begin{scope}[scale = 0.6]
\draw[black, thick, fill = black] (2,0) circle (1 mm) node[above]{$_1$};
\draw[black, thick] (2,0) -- (0,-2.2) -- (1.8, -4) -- (1,-6.4);
\draw[black, thick, dashed] (1.8, -4) -- (3.8,-2.8);
\draw[red, thick, dashed] (1, -6.4) -- (3.2,-7.2);
\draw[black, thick, fill = black] (0,-2.2) circle (1 mm) node[left]{$_5$};
\draw[black, thick, fill = black] (1.8, -4) circle (1 mm) node[below right]{$_3$};
\draw[black, thick, fill = black] (1,-6.4) circle (1 mm) node[below]{$_4$};
\draw[red, thick, fill = black] (3.2,-7.2) circle (1 mm) node[above]{$_6$};
\draw[black, thick, fill = black] (3.8,-2.8) circle (1 mm) node[above]{$_2$};
\draw[black, thick] (1.1,-5.5) node[above]{$_i$};
\draw[black, thick] (0.8,-1.2) node[above]{$_{ii}$};
\draw[black, thick] (2.8,-3.3) node[above]{$_{iii}$};
\draw[black, thick] (1.25,-3.3) node[above]{$_{iv}$};
\draw[red, thick] (2.2,-6.8) node[above]{$_o$};
\draw[black, thick, -{Latex[length=3mm, fill = black]}] (5.8,-3) -- (7.5,-3);
\end{scope}
%%%%%%%%%%%%%%%%%%%%%%%%%%%%%%%%%%%%%%%%%%%%%%%%%%%%%%%%%%%%%%
%                                                            %
%%%%%%%%%%%%%%%%%%%%%%%%%%%%%%%%%%%%%%%%%%%%%%%%%%%%%%%%%%%%%%
\begin{scope}[xshift = 6cm, scale = 0.6]
\draw[black, thick, -{Latex[length=2mm, fill = black]}] (0,0) node[above]{$_1$} -- (0.8,-2.4) -- (-0.6, -3.2);
\draw[black, thin, loosely dotted] (-1,0)  -- (10,0) node[right]{\small time~$t$};
\draw[black, thick, fill = black] (0.8,-2.4) circle (0.5 mm);
\draw[black, thick, fill = black] (0.95,-2.4) circle (0.5 mm) node[right]{$_5$};
\draw[black, thick, fill = black] (6.6,-1.4) circle (0.5 mm);
\draw[black, thick, fill = black] (6.75,-1.4) circle (0.5 mm) node[right]{$_4$} ;
\draw[black, thick, fill = black] (2.8,-4.8) circle (0.5 mm);
\draw[black, thick, fill = black] (2.95,-4.8) circle (0.5 mm);
\draw[black, fill = white] (4.73,-3.13) circle (1.5 mm);
\draw[black, thick, -{Latex[length=2mm, fill = black]}] (0.95,-2.4) -- (2.8,-4.8) -- (1.6,-5.8);
\draw[black, thick, -{Latex[length=2mm, fill = black]}] (2.8,0) node[above]{$_2$} -- (6,-5.2);
\draw[black, thick, -{Latex[length=2mm, fill = black]}] (6.2, 0) node[above]{$_3$} -- (6.6,-1.4) -- (2.95,-4.8) -- (4,-6.6) ;
\draw[black, thick, -{Latex[length=2mm, fill = black]}] (6.75,-1.4) -- (8.6,-3.8);
\end{scope}
\end{tikzpicture}
\caption{Example of retrieving the pseudo-trajectory dynamics (right) from the dynamics tree (left)} \label{fig:dyn_tree}
\end{figure}

\paragraph{Clustering trees}
At this point, we observe that the knowledge of a recollision tree on each aggregate, and of an overlap tree between these aggregates, is equivalent to the knowledge of a global decorated \emph{clustering tree}~$T^\star_R \in \mT^{\prec, \star}_r$, recording all the clustering interactions between the $r$~particles of a cluster~$R$, and recording whether it is a recollision or an overlap through a decoration on each edge. This leads to the mapping
\begin{equation*} \def\arraystretch{2}
\begin{array}{ccc}  \mT_{|\kappa|}^\prec \times \prod_{i=1}^{|\kappa|} \mT_{|\kappa_i|}^\prec & \longrightarrow  & \mT_{r}^{\prec, \star} \\
T^\ov_R, \bigl(T^\col_{\kappa_i} \bigr)_{i \leqslant |\kappa|} & \longmapsto & T^\star_R .
\end{array}
\end{equation*}
The decorations are denoted~$\us_{E_\Upsilon}$ by analogy with the scattering labels of the pseudo-trajectory histories~\eqref{def:history}: indeed, this time again, we choose that a sign~$+1$ corresponds to a scattering (i.e. a recollision), and a sign~$(-1)$ to no scattering (i.e. an overlap).
In the decorated clustering tree, the aggregates correspond to the decorated connected components retrieved when removing the overlap edges (we denote~$\mathfrak{CC}^\star(\Upsilon)$ the set of these decorated connected components), so that the tree part of our cluster measure writes
\begin{align*}
\sum_{\kappa \in \mP_R}\phi^{[\tree]}_{R} \prod_{i=1}^{|\kappa|} \cls(\kappa_i) \dd \nu^{[H]}_{[0,t]}(\Psi_{\kappa_i})  = \sum_{\Upsilon \in\mT^{\prec, \star}_R} (-1)^{|\mathfrak{CC}^\star(\Upsilon)|-1}  \int \ind_{T^\star_R = \Upsilon} \prod_{A \in \mathfrak{CC}^\star(\Upsilon)} \dd \nu^{[H]}_{[0,t]}(\Psi_{A}).
\end{align*}
For each global tree~$\Upsilon \in \mT_R^{\prec, \star}$, we thus have~$(r-1)$~encounter times and angles $(\tau^\star_e, \omega^\star_e)_{e\in E_\Upsilon}$, either associated to a recollision or to an overlap, according to the decorations~$(s_e)_{e \in E_\Upsilon}$. Once computed the recollision changes of variable in each aggregate (i.e. decorated connected component of~$\Upsilon$), we compute the same on their barycenters to harness the overlaps, eventually leading to the variables~$(\overline{x}_R, \utau^\star_{r-1}, \uom^\star_{r-1})$. Since the overlap change of variable has the same Jacobian as the recollision one~\eqref{eq:jacob}, we end up with
\begin{align*}
& \int \dd \ux_R \sum_{\kappa \in \mP_R}\phi^{[\tree]}_{R} \prod_{i=1}^{|\kappa|} \cls(\kappa_i) \dd \nu^{[H]}_{[0,t]}(\Psi_{\kappa_i})
 = \sum_{\Upsilon \in \mT_R^{\prec, \star}}  \sum_{\substack{p_e \in P(m), q_e \in P(m') \\ e=(m,m')\in  E_\Upsilon}}   \int \frac{\dd \overline{x}_R}{\mu^{r-1} } \dd \uom^\star_{E_\Upsilon} \dd \utau^\star_{E_\Upsilon} \\ 
 & \hspace{35mm} \times \ind_{T^\star_R = \Upsilon} \mathfrak{S}_{(p_e, q_e)_{\Upsilon}}[R]  \prod_{e \in E_\Upsilon} s_{e} \langle v_{p_e}({\tau^\star_e}^+ )  - v_{q_e}({\tau^\star_e}^+),  \omega^\star_e \rangle_+ \prod_{A \in \mathfrak{CC}^\star(\Upsilon)} \dd \nu^{[H]}_{[0,t]}(\Psi_{A}).
\end{align*}
Now, keeping in mind that the pseudo-trajectories do not interact between decorated connected components of the clustering tree, one may gather back the measures~$\nu^{[H]}_{[0,t]}$ on the whole set~$R$. This means that we consider a unique history
\begin{equation*}% \label{eq:history2}
(k, \um_k, \ul_k^*, \us_k, \ut_k, \uom_k, \uv_k^*) 
\end{equation*}
for all the dynamics, although the pseudo-trajectories are still independently defined on each aggregate. Moreover, we want to keep track of the time ordering between the particle adjunctions and the clustering encounters. To that end, we will add to each edge~$e \in E_{\Upsilon}$ of the clustering tree the information of the time interval~$[t_{j+1}, t_j]$ (between two particle adjunctions) in which the clustering encounter occur. This way, adding this information to the knowledge of which particles are involved in the clustering encounter, each edge~$e=\{m,m'\}$ is associated with the triplet
\begin{equation*}
(p_e, q_e, i_e) \in P(m) \times P(m') \times \lbr 0,k\rbr,
\end{equation*}
such that $\tau_e \in [t_{i_e +1}, t_{i_e}]$. Keeping the notation $\mathfrak{S}_{(p_e, q_e, i_e)}$ for this triplet's compatibility, we write
\begin{align*}
 \int \dd \ux_R \sum_{\kappa \in \mP_R}\phi^{[\tree]}_{R} \prod_{i=1}^{|\kappa|} \cls(\kappa_i) \dd \nu^{[H]}_{[0,t]}(\Psi_{\kappa_i})
 & = \sum_{\Upsilon \in \mT_R^{\prec, \star}}  \sum_{\substack{p_e, q_e, i_e \\ e\in  E_\Upsilon}}  \int \frac{\dd \overline{x}_R}{\mu^{r-1} } \dd \uom^\star_{E_\Upsilon} \dd \utau^\star_{E_\Upsilon} \ind_{T^\star_R = \Upsilon} \mathfrak{S}_{(p_e, q_e, i_e)_{\Upsilon}}[R] \\ 
 & \hspace{10mm}  \times   \prod_{e \in E_\Upsilon} s_{e} \langle v_{p_e}({\tau^\star_e}^+ )  - v_{q_e}({\tau^\star_e}^+),  \omega^\star_e \rangle_+ \ \dd \nu^{[H]}_{[0,t]}(\Psi_{R}).
\end{align*}
With all this in mind, we are now able to gather all the particles adjunctions and clustering encounters into a single tree encoding the whole dynamics. Indeed, at fixed number of added particles~$k$, the choice of the history~$(\us_k, \um_k)$ contained in the measure~$\nu$, along with the clustering information, may be encoded in a dynamics tree~$T^\dd_R \in \mT_{r+k}^{\prec, \star}$, in which two vertices are connected either if a clustering encounter happens between them, or if one of them was added to the other one, i.e. if they form a couple~$(n+j, m_j)$, with the associated signed decoration. The edges are time ordered from the ordered clustering conditions in~$\Upsilon = T^\star_R$, and the intervals~$(i_e)_{e \in E_\Upsilon}$ defined above. This provides the following mapping
\begin{equation*} \def\arraystretch{1.2}
\begin{array}{cccccccccccccc}
\mT_R^{\prec, \star} &\times& \lbr 1, r+k \rbr^{E_\Upsilon} & \times & \lbr 0, k \rbr^{E_\Upsilon} &  \times& \{\pm 1\}^k& \times &\mM_{R,k} & \rightarrow & \mT_{R,k}^\dd \subset \mT_{r+k}^{\prec, \star} \\
T^\star_R &,& (p_e, q_e)_e & , & (i_e)_e & ,& \us_k &,& \um_k & \mapsto & T^\dd_R
\end{array}
\end{equation*}
with the condition $p_e, q_e \in P(m) \times P(m')$ for every $e = \{m,m'\} \in E_\Upsilon$. This mapping is injective, and made bijective by restriction to its image~$\mT_{R,k}^\dd \subset \mT_{r+k}^{\prec, \star}$, containing all the admissible dynamics trees. More precisely, an admissible dynamics tree is such that the added particles from~$\lbr r+1, r+k \rbr$ appear in increasing order, and in particular do not have encounters before appearing. For such a dynamics tree, the global dynamics might be recovered, starting from the highest ordered edge (the closest to time~0) and reconstructing the trajectories up to time~$t$, as in the example given in Figure~\ref{fig:dyn_tree} (note that in this example the edge~0 is not admissible as it would connect the vertices~4 and~6 before they appear). From the dynamics, the parameters $(T^\star_R, (p_e, q_e, i_e)_e, \us_k, \um_k)$ can be retrieved, making the mapping bijective.
Recalling the expression~\eqref{def:measure} of the measure, this lets us with encounter times~$\utau_{E_\Upsilon}$, angles~$\uom_{E_\Upsilon}$, stemming either from the adjunction of a particle in the pseudo-trajectory expansion, or from a clustering encounter (recollision or overlap). Generalizing the notation~$\langle v_{m}({\tau_e}^+ )  - v_{m'}({\tau_e}^+),  \omega_e \rangle_+$ to the adjunction of particles (even if in this case one of the particles does not really exist at time~$\tau_e^+$), one can write the condensed measure
\begin{align*}
&\hspace{0mm}\int \dd \ux_R \sum_{\kappa \in \mP_R}\phi^{[\tree]}_{R} \prod_{i=1}^{|\kappa|} \cls(\kappa_i) \dd \nu^{[H]}_{[0,t]}(\Psi_{\kappa_i}) \\
&  = \sum_{k \geqslant 0} \sum_{\ul_k \in \Lambda_k}  \sum_{\Upsilon \in \mT_{R,k}^{\dd}}  \int \frac{p_{\mu}^{|\ul_k|}}{\mu^{r-1} }\dd \overline{x}_R   \dd \uv^*_k
 H^{\otimes R} \ind_{T^\dd_R = \Upsilon} \prod_{e = \{m,m'\} \in E_\Upsilon} s_{e} \langle v_{m}({\tau_e}^+ )  - v_{m'}({\tau_e}^+),  \omega_e \rangle_+ \dd \omega_e \dd \tau_{e} .
\end{align*}
Although the order in the velocities' difference above do not change the calculus (as we integrate~$\omega_e$ over the whole symmetric sphere), recall anyway that this order has been chosen such that~$m > m'$.
Summing everything up, we get the following proposition.
\begin{prop}[Dynamics tree formula for the integrated cumulants] \label{prop:dyn_trees}
For any time $t>0$, the cumulants decompose into two terms 
\begin{align*}
f^\varepsilon_n[H](t)
& = f^\varepsilon_n[H]^{[\tree]}(t)  + f^\varepsilon_n[H]^\cycle(t), 
\end{align*}
the first one containing only trees in its expansion, and the other including cycles. Indeed, the integral of the first tree term rewrites as the expansion
\begin{align*}
\int f^\varepsilon_n[H]^{[\tree]}(\uz_n) \dd \uz_n
& = \sum_{\rho \in \mP_n} I_n^{[\rho]}[H](t, \ul_n),
\end{align*}
where 
\begin{align}
I_n^{[\rho]}[H](t, \ul_n) & \doteq \int f^{\varepsilon,\rho}_{n}(0, \uz^{[0]}_{\Psi_{\rho_1}}, \dots, \uz^{[0]}_{\Psi_{\rho_{|\rho|}}}) \dd \uv_{n} \prod_{j=1}^{|\rho|} \frac{\dd \overline{x}_{\rho_j}}{\mu^{|\rho_j|-1}} \sum_{k_j \geqslant 0} \sum_{\ul_{k_j}} \sum_{\Upsilon \in \mT_{\rho_j, k_j}^{\dd}} p_{\mu}^{|\ul_{k_j}|} H^{\otimes \rho_j} \ind_{T^\dd_{\rho_j} = \Upsilon}  \nonumber \\ 
& \hspace{35mm}  \times \dd \uv^*_{k_j}  \prod_{e \in E_\Upsilon} s_{e} \langle v_{m}({\tau_e}^+ )  - v_{m'}({\tau_e}^+),  \omega_e \rangle_+ \dd \omega_e \dd \tau_{e}. \label{def:In_rho_H}
\end{align}
The term $f^\varepsilon_n[H]^\cycle$ corresponds to a similar expansion, stemming from the part~$\phi_{\rho_j}^\cycle$ of the exclusion cumulants~\eqref{def:tree_cycle}. Recall that the cluster cumulants $f_n^{\varepsilon, \rho}$ were defined in~\eqref{eq:initial_cluster_cumulants_1}.
\end{prop}

\subsection{Discarding overlap cycles: a tree inequality }

\label{sec:DiscCycl}
We start the estimates on the cumulants by controlling the part of the cumulant expansion containing cycles.
A cycle in the overlap conditions imposes a strong geometric constraint, which will provide smallness (Section~\ref{app:geom_recoll}). But first, we will compute a \emph{tree inequality} to simplify the sum over all the connected graphs (see~\cite{fouthese} for details about the general tree inequality).
 Indeed, recall that the cycle part of the exlusion cumulants writes, for a general cluster~$R$,
\begin{align*}
\phi_{R}^\cycle = \sum_{\Upsilon \in\mT^\prec_{|\kappa|}} \ind_{T^\ov(R) = \Upsilon} \sum_{G \in \mC_{|\kappa|}(\Upsilon)\setminus \{\Upsilon\}}  \prod_{\{i,j\} \in E_G} (-\ind_{\kappa_i \sim \kappa_j}).
\end{align*}
Denoting~$\ind^{\between}_R$ the indicator that the global dynamics contains a cycle, the sum over the cycling graphs rewrites as 
\begin{align*}
\sum_{G \in \mC_{|\kappa|}(\Upsilon)\setminus \{\Upsilon\}}  \prod_{\{i,j\} \in E_G} (-\ind_{\kappa_i \sim \kappa_j}) & = \ind^{\between}_R \prod_{\{i,j\} \in E_\Upsilon} (-\ind_{\kappa_i \sim \kappa_j}) \sum_{E' \subset E_\Upsilon^c} \prod_{\{i,j\} \in E'} ( -\ind_{\kappa_i \sim \kappa_j}) \\
& = \ind^{\between}_R \prod_{\{i,j\} \in E_\Upsilon} (-\ind_{\kappa_i \sim \kappa_j}) \prod_{\{i,j\} \in E_\Upsilon^c} (1 -\ind_{\kappa_i \not\sim \kappa_j})
\end{align*}
computing the inverse expansion of the exclusion cumulants, the latter product being smaller than~1. We are hence brought back to the tree case, and dominate the cycle part by the tree part, with the additional strong cycle condition
\begin{align*}
\int f^\varepsilon_n[H]^{[\mathrm{cycle}]}
& \leq \sum_{\rho \in \mP_n} \int f^{\varepsilon,\rho}_{n}(0)\dd \uv_{n} \prod_{j=1}^{|\rho|} \frac{\dd \overline{x}_{\rho_j}}{\mu^{|\rho_j|-1}} \sum_{k_j \geqslant 0} \sum_{\ul_{k_j}} \sum_{\Upsilon \in \mT_{\rho_j, k_j}^{\dd}} p_{\mu}^{|\ul_{k_j}|} H^{\otimes \rho_j} \ind_{T^\dd_{\rho_j} = \Upsilon}  \nonumber \\ 
& \hspace{25mm}  \times \dd \uv^*_{k_j}  \prod_{e \in E_\Upsilon} s_{e} \langle v_{m}({\tau_e}^+ )  - v_{m'}({\tau_e}^+),  \omega_e \rangle_+ \dd \omega_e \dd \tau_{e} \times \ind^\between_{\rho_j}.
\end{align*}
This term is negligible in front of the tree one, because of the strong geometric cycle condition, as stated in Section~\ref{sec:disc_cycles}.

\subsection{Bounding the cumulants on short times} \label{sec:bounds}
We hence need to study the tree part of the integrated cumulants, written (Proposition~\ref{prop:dyn_trees})
\begin{align}
\int f^\varepsilon_n[H]^{[\tree]}(t)
& = \sum_{\rho \in \mP_n} \int f^{\varepsilon,\rho}_{n}(0) \dd \uv_{n} \prod_{j=1}^{|\rho|} \frac{\dd \overline{x}_{\rho_j}}{\mu^{|\rho_j|-1}} \sum_{k_j \geqslant 0} \sum_{\ul_{k_j}} \sum_{\Upsilon \in \mT_{\rho_j, k_j}^{\dd}} \int \dd \uv^*_{k_j}  \ind_{T^\dd_{\rho_j} = \Upsilon}   \label{eq:cumulant_exp_int} \\ 
& \hspace{55mm} \times  \prod_{e \in E_\Upsilon} s_{e} \langle v_{m}({\tau_e}^+ )  - v_{m'}({\tau_e}^+),  \omega_e \rangle_+ \dd \omega_e \dd \tau_{e}. \nonumber
\end{align}
As in~\cite{2023grandev}, this bound is only valid on short times, considering that we do not have a priori bounds for the cumulants, unlike for the correlation functions~\cite{fou26Rayleigh}.
\begin{prop}[Cumulant bound] \label{prop:cum_bound}
For $\mu$ large enough in the mixed scaling~\eqref{eq:scaling}, there exists a time $T_\beta$ depending on~$\beta$ and~$d$ such that for any $t< T_\beta$, one has for an absolute constant $C > 0$ that
\begin{equation*}
\left\lvert \int f^\varepsilon_n[H]^{[\tree]}(t) \right\rvert \leq \frac{C^n n!}{\mu^{n-1}} \| H \|_\infty^n.
\end{equation*}
More precisely, we have
\begin{equation} \label{ineq:In_1_n}
 \left\lvert  \int f^\varepsilon_n[H]^{[\tree]}(t) - I_n^{ \lbr 1, n \rbr }[H](t) \right\rvert \leq \frac{C^n n!}{\mu^{n-1}} \| H \|_\infty^n \cdot \varepsilon.
\end{equation}
The latter inequality specifies that the leading term in the integral formula for the cumulants is the one corresponding to the trivial partition $\rho = \{ \lbr 1, n \rbr \}$.
\end{prop}
We denote from now on the total kinetic energy, preserved at fixed number of particles by the transport and by elastic collisions,
\begin{equation*}
\| \uv_k\|^2 \doteq \sum_{i = 1}^k |v_i|^2,
\end{equation*} 
where $|v_i|$ is the Euclidean norm of the velocity $v_i \in \R^d$. We also denote as in~\cite{fou26Rayleigh} the initial bound
\begin{equation} \label{def:C0}
C_0 \doteq \max \left[ \| M_\beta \|_{1, \beta/2} \ ; \ \| M_\beta \varphi_0 \|_{1, \beta/4}\ ; \  \left\| M_\beta M_{\beta/2}^{-\frac{1}{2}} \varphi_0 \right\|_{\Lp^\infty(\mD)}  \right].
\end{equation}
\begin{prof}
\textbf{First step: initial cluster cumulants.}
To get bounds on this cumulant, we will start integrating the \emph{initial} cluster cumulants, of which we recall the following expression~(Lemma~\ref{lem:clust_cum}), where~$N \doteq |\Psi_{\rho_1}| + \dots + |\Psi_{\rho_{|\rho|}}| $ denotes the total number of particles contained in the clusters,
\begin{align*}
 f_n^{\varepsilon,\rho}(0, \uz^{[0]}_{\Psi_{\rho_1}}, \dots, \uz^{[0]}_{\Psi_{\rho_{|\rho|}}}) = \left[ \ind_{\mX_{(\cdot)}} \right]^{\otimes \aleph_{|\rho|}} M_\beta^{\otimes N} \varphi_0^{\otimes \ul_{N}}  \sum_{p, q \geqslant 0} \frac{\lambda^{p} \mu^{q}}{p!q!}  \int  \dd \uz_{p+q}^* \varphi_0^{\otimes p} M_\beta^{\otimes p + q} \phi_{\aleph_{|\rho|}, p+q}.
\end{align*}
Using the definition~\eqref{def:C0} of~$C_0$ to control~$M_\beta$ and~$\varphi_0$, we have
\begin{align*}
 \int \left\lvert  f^{\varepsilon,\rho}_{n}(0, \uz^{[0]}_{\Psi_{\rho_1}}, \dots, \uz^{[0]}_{\Psi_{\rho_{|\rho|}}}) \right\rvert \prod_{k=1}^{|\rho|} \dd \overline{x}_{\rho_k} & \leq C_0^N e^{-\beta \| \uv_N\|^2}\sum_{p, q \geqslant 0} \frac{(C_0\lambda)^{p} \mu^{q}}{p!q!} \int  M_\beta^{\otimes p+q} \left\lvert\phi_{\aleph_{|\rho|}, p+q} \right\rvert \prod_{k=1}^{|\rho|} \dd \overline{x}_{\rho_k}  \dd \uz_{p+q}^*,
\end{align*}
and we apply the tree inequality~\cite[Proposition~A.1]{fou26Rayleigh} to the cumulant of the exclusion
\begin{equation*}
\left\lvert \phi_{\aleph_{|\rho|}, p+q} \right\rvert \leq  \sum_{T \in \mT_{\aleph_{|\rho|} \cup \lbr p+q \rbr}} \prod_{\{x,y\} \in E_G} \ind_{x\sim y}.
\end{equation*}
Like in the proof of the bound on these cumulants~\cite[Proposition~A.1]{fou26Rayleigh}, we use the integration variables $\uz_{p+q}^*$ and $(\overline{x}_{\rho_k})_k$, the latter moving rigidly the clusters $ (\ux_{\Psi_k})_k$. Integrating over successive leaves of the tree, removing the edge $\{i,j\}$ leads to a factor $|\Psi_{\rho_i}|\cdot |\Psi_{\rho_j}| C_d \varepsilon^d$, depending on the number of particles in each cluster. For this reason, this time we need to discriminate according to the degrees $d_1, \dots, d_{|\rho|+p+q}$ of the vertices $\Psi_{\rho_1}, \dots, \Psi_{\rho_{|\rho|}}, x_1^*, \dots, x_{p+q}^*$. Since the	number of trees on $\lbr 1, m \rbr$ with prescribed degrees $d_1, \dots, d_m$ is equal to 
\begin{equation*}
\frac{(m-2)!}{\prod_{i=1}^{m} (d_i - 1)!},
\end{equation*}
we get (also integrating the density~$M_\beta$ over $\uv_{p+q}^*$)
\begin{align*}
& \int \left\lvert f^{\varepsilon,\rho}_{n}(0, \uz^{[0]}_{\Psi_{\rho_1}}, \dots, \uz^{[0]}_{\Psi_{\rho_{|\rho|}}}) \right\rvert \prod_{k=1}^{|\rho|} \dd \overline{x}_{\rho_k}\\
& \hspace{25mm} \leq C_0^N e^{-\beta  \|\uv_N\|^2} \sum_{p, q \geqslant 0} \frac{(C_0\lambda)^{p} \mu^{q}}{p!q!} \varepsilon^{d(|\rho|+p+q-1)} \sum_{d_1, \dots, d_{|\rho|+p+q}}  \frac{(|\rho|+p+q-2)!}{\prod_{i=1}^{|\rho|+p+q} (d_i - 1)!} \prod_{i=1}^{|\rho|} |\Psi_{\rho_i}|^{d_i}.
\end{align*}
Now, we observe that
\begin{equation*}
\frac{(|\rho|+p+q-2)!}{p!q!} = \binom{p+q}{p} \binom{p+q+|\rho|-2}{p+q} (|\rho|-2)! \leq \binom{p+q}{p} 2^{p+q+|\rho|-2} (|\rho|-2)!\ .
\end{equation*}
We write on the other hand
\begin{align*}
\prod_{i=1}^{|\rho|} \sum_{d_i \geqslant 1}  \frac{|\Psi_{\rho_i}|^{d_i}}{ (d_i - 1)!} \leq \prod_{i=1}^{|\rho|} |\Psi_{\rho_i}| \exp( |\Psi_{\rho_i}|) \leq e^{2N},
\end{align*}
and similarly 
\begin{align*}
\prod_{i=1}^{p+q} \sum_{d_{|\rho|+i} \geqslant 1}  \frac{1}{ (d_{|\rho|+i} - 1)!}  \leq e^{p+q}.
\end{align*}
We eventually get, since $|\rho| \leq N$ and $C_0 \lambda \leq \mu$,
\begin{align*}
 \int \left\lvert f^{\varepsilon,\rho}_{n}(0, \uz^{[0]}_{\Psi_{\rho_1}}, \dots, \uz^{[0]}_{\Psi_{\rho_{|\rho|}}}) \right\rvert \prod_{k=1}^{|\rho|} \dd \overline{x}_{\rho_k} & \leq (2e^2 C_0)^N e^{-\beta \|\uv_N\|^2} (|\rho|-2)! \varepsilon^{d(|\rho|-1)}  \sum_{p, q \geqslant 0} \binom{p+q}{p} (2 \mu \varepsilon^d)^{p+q}  \\
 & \leq (2e^2 C_0)^{N} e^{-\beta \|\uv_N\|^2}  (|\rho|-2)! \varepsilon^{d(|\rho|-1)} \sum_{r \geqslant 0} (4 \varepsilon )^{r},
\end{align*}
eventually leading, for an absolute constant $C$, to
\begin{align} \label{ineq:initial_cum}
 \int \left\lvert f^{\varepsilon,\rho}_{n}(0, \uz^{[0]}_{\Psi_{\rho_1}}, \dots, \uz^{[0]}_{\Psi_{\rho_{|\rho|}}}) \right\rvert \prod_{k=1}^{|\rho|} \dd \overline{x}_{\rho_k} & \leq (CC_0)^{N} e^{-\beta \|\uv_N\|^2}  (|\rho|-2)! \varepsilon^{d(|\rho|-1)}.
\end{align}

\noindent
\textbf{Second step: bounding the collision kernels.} The observable~$H$ being bounded, we can dominate it roughly. Moreover, let us observe that the factors bounding the initial cumulants~\eqref{ineq:initial_cum} decompose on the clusters as
\begin{equation} \label{eq:veloc_decay}
(CC_0)^N e^{-\beta \|\uv_N \|^2} = \prod_{j=1}^r (CC_0)^{|\Psi_{\rho_j}|} e^{-\beta \|\uv_{\Psi_{\rho_j}} \|^2},
\end{equation}
so that for a generic cluster~$R = \rho_j$, we harness the associated velocity decay to study
\begin{align*} 
S_{R} \doteq \sum_{k \geqslant 0} \sum_{\ul_{k} \in \Lambda_k} p_\mu^{|\ul_k|} \sum_{\Upsilon \in \mT_{R, k}^{\dd}} \int \dd \uv_{r+k} \ind_{T^\dd_{\rho_j} = \Upsilon}  \prod_{e \in E_\Upsilon} s_{e} \langle v_{m}({\tau_e}^+ )  - v_{m'}({\tau_e}^+),  \omega_e \rangle_+ \dd \omega_e \dd \tau_e e^{-\beta\|\uv_{r+k}\|^2}.
\end{align*}
With our notation, $r+k = |\Psi_R|$.
Since we just bounded the dependency on the initial value, the trajectories do not depend anymore on the particles' tags, allowing us to bound roughly the sum over $\ul_k$, using that $p_\mu \leq 1$.
We will use the velocity decay to bound the velocities appearing in the collision kernels. We bound each edge separately, starting from the furthest from time~$t$, so that the other collision kernels will only contain velocities that do not depend on the considered edge. 

For each edge corresponding to a clustering condition, using the Cauchy-Schwartz inequality, we bound all the possibilities for this edge by
\begin{equation} \label{ineq:sum_edge_1}
2 \sum_{m,m' = 1}^{r+k} (|v_{m}({\tau^\star_e}^+ )| + | v_{m'}({\tau^\star_e}^+)| ) \leq 4 (r+k) \sqrt{(r+k) \|\uv_{r+k} \|^2},
\end{equation}
where the factor~2 stems from the choice of the edge's decoration. For each edge associated to a particle's adjunction, since one of the vertices corresponds to the particle that is added, we only have to sum over the choice of the second particle as
\begin{equation} \label{ineq:sum_edge_2}
2 \sum_{m = 1}^{r+k} (|v_{m}({\tau^\star_e}^+ )| + | v_{m'}({\tau^\star_e}^+)| ) \leq 2  \sqrt{(r+k) \|\uv_{r+k} \|^2}.
\end{equation}
Determining whether the edge is a clustering condition or a particle's adjunction corresponds to an additional factor~$2^{r+k}$. In the end, harnessing the velocity decay~\eqref{eq:veloc_decay} and integrating the ordered times, we get 
\begin{align}
S_{R} & \leq \sum_{k \geqslant 0} 8^k  \int \dd \uv_{r+k} \dd \uom_{E_\Upsilon} \dd \utau_{E_\Upsilon}  (r+k)^{r-1} \left(\sqrt{(r+k) \|\uv_{r+k} \|^2} \right)^{r+k-1} e^{-\|\uv_{r+k}\|^2} \nonumber \\
 &\leq  \sum_{k \geqslant 0} 8^k \frac{ t^{r+k-1}}{(r+k-1)!} (r+k)^{r-1} \sqrt{r+k}^{r+k-1} \left( C_d \sqrt{r+k} \right)^{r+k-1} \nonumber \\
 &\leq  \sum_{k \geqslant 0} (\hat{C}_d t)^{r+k-1} (r+k)^{r-1} e^{r+k}, \label{lgn:dom_SR}
 \end{align}
using in the end that $(r+k)^{r+k-1} \leq (r+k-1)! e^{r+k}$.

\noindent
\textbf{Last step: combinatorial manipulations.}
We now gather it all to bound the cumulants using their expression~\eqref{eq:cumulant_exp_int} from Proposition~\ref{prop:dyn_trees}. 
Using the exchangeability of particles to reduce the partition~$\rho \in \mP_n$ to the cardinal of its subsets, and denoting~$C_H \doteq \|H \|_\infty$ and $\hat{C} \doteq CC_0$, we write for a constant~$C_d$ depending only on the dimension that
\begin{align*}
 \left\lvert \int f^\varepsilon_n[H]^{[\tree]}(t) \right\rvert 
& \leq \sum_{i = 1}^n  \sum_{r_1 + \dots + r_i = n} \frac{(i-2)!\varepsilon^{d(i-1)} n!}{i! r_1! \dots r_i!}  \frac{C_H^n}{\mu^{n - i}} \prod_{j=1}^i \sum_{k_j \geqslant 0} (C_d t)^{r_j+k_j-1} (r_j+k_j)^{r_j-1} \\
&\leq \frac{C_H^n n!}{\mu^{n - 1}} \sum_{i = 1}^n  \varepsilon^{i-1}    \prod_{j=1}^i \sum_{r_j \geqslant 0} \sum_{k_j \geqslant 0} (C_d t)^{r_j+k_j-1} e^{k_j + r_j} 
\end{align*} 
in the scaling $\mu \varepsilon^d = \varepsilon$, and harnessing the denominators~$r_j !$ to control the terms~$(r_j + k_j)^{r_j-1}$. Now, for~$t$ small enough the series are convergent (and the term $r_j = k_j = 0$ actually does not appear, as it corresponds to the empty tree), so that we get
\begin{align*}
 \left\lvert \int f^\varepsilon_n[H]^{[\tree]}(t) \right\rvert 
& \leq \frac{C_H^n n!}{\mu^{n - 1}} \sum_{i = 1}^n  \varepsilon^{i-1} \cdot  4^i \\
& \leq 4 \frac{C_H^n n!}{\mu^{n - 1}} (1 + 8\varepsilon),
\end{align*} 
which concludes the proof. The second point~\eqref{ineq:In_1_n} of the proposition is proved observing that the term $8\varepsilon$ corresponds to the sum for $i \geq 2$, i.e. to the sum over the non-trivial partitions. \CQFD
\end{prof}

\subsection{Limit cumulants}  \label{sec:limit_cumulants}
Since the cumulants~$f_n^\varepsilon$ decrease very quickly as $n$ increases (Proposition~\ref{prop:cum_bound}), we only need the convergence of the first cumulant to compute the large deviation principle, and the convergence of the fluctuation field. Indeed, in the rescaling~$(\lambda^{n-1} f_n^\varepsilon)$ that appears in the statistical study of the tagged particles, the cumulants vanish for $n\geqslant 2$, as soon as $\textstyle \frac{\lambda}{\mu}$ goes to zero.
Nevertheless, for completeness, we also compute their limit in the rescaling~$(\mu^{n-1} f_n^\varepsilon)$. The corresponding formula could be used to capture the fine scales of the dynamics in a further study.

\subsubsection{Rare encounters of tagged particles} \label{sec:rare_encounters2}
First of all, we will show that the encounters between tagged particles are very rare. Indeed, when we have bounded roughly the sum over the possible tags~$\ul_k \in \Lambda_k$ by~$2^k$, one could have looked at the influence of the encounters with at least one tagged particle, which corresponds to the sum
\begin{align*}
\sum_{\substack{\ul_k \in \Lambda_k \\ \ul_k \neq \underline{0}_k }} p_\mu^{|\ul_k|} & = \sum_{i=1}^k	\binom{k}{i} p_\mu^i\\
& \leq p_\mu 2^k.
\end{align*}
This leads to the exact same bound as the previous one, with an additional factor~$p_\mu \ll 1$. In the end, we get the estimate presented in the proposition below.
\begin{prop}[Rare tagged encounters] \label{prop:rare_encounters2}
Denoting 
\begin{equation*}
\int f_n^\varepsilon[H]^{[\tree], \underline{0}}(t)
\end{equation*}
the integrated cumulants where the domain of integration of $\dd \nu$ is restricted to $\{\ul_k^* = \underline{0}_k\}$, one has in the usual mixed scaling~\eqref{eq:scaling}, for any $t$ small enough and an absolute constant $C > 0$, that
\begin{equation*}
\left\lvert \int f^\varepsilon_n[H]^{[\tree]}(t) - \int f^\varepsilon_n[H]^{[\tree], \underline{0}}(t) \right\rvert \leq \frac{C^n n!}{\mu^{n-1}} \| H \|_\infty^n\cdot p_\mu.
\end{equation*}
\end{prop}

\subsubsection{Discarding cycles} \label{sec:disc_cycles}

In the expansion of the integrated cumulants in terms of dynamics trees (Section~\ref{sec:dyn_trees}), we emphasized that non-clustering encounters may happen, either stemming from a non-clustering overlap, or from a non-clustering collision. We state here that the dynamics presenting such cycles are negligible in the said expansion, in particular proving the smallness of the cycle part of the expansion (Section~\ref{sec:DiscCycl}). Using a more precise computation than the previous method~\cite{2023longcor}, we achieve an optimal bounding factor in~$\varepsilon$ instead of $\varepsilon |\log \varepsilon|$.
\begin{prop}[Cycles are rare in the dynamics]	\label{prop:discard_cycles}
We have the following estimate on the expansion proved in Proposition~\ref{prop:dyn_trees} for the integrated cumulants, under the constraint that a cycle happens in the dynamics:
\begin{align*}
& \Bigg\lvert \sum_{\rho \in \mP_n} \int f^{\varepsilon,\rho}_{n}(0)\dd \uv_{n} \prod_{j=1}^{|\rho|} \frac{\dd \overline{x}_{\rho_j}}{\mu^{|\rho_j|-1}} \sum_{k_j \geqslant 0}  \sum_{\Upsilon \in \mT_{\rho_j, k_j}^{\dd}} H^{\otimes \rho_j} \ind_{T^\dd_{\rho_j} = \Upsilon}  \nonumber \\ 
& \hspace{18mm}  \times \dd \uv^*_{k_j} \dd \uom_{E_\Upsilon} \dd \utau_{E_\Upsilon}  \prod_{e \in E_\Upsilon} s_{e} \langle v_{m}({\tau_e}^+ )  - v_{m'}({\tau_e}^+),  \omega_e \rangle_+ \times \ind^\between_{\rho_j} \Bigg \rvert \leq \frac{C^n n!}{\mu^{n-1}} \| H \|_\infty^n\cdot \varepsilon.
\end{align*}
The cycle condition imposes strong geometric constraints that result in a additional factor~$\varepsilon$ in the usual bounds for the integrated cumulants. The dynamics implying cycles are thus negligible in the expansion above.
\end{prop}
The proof of this proposition is given in Section~\ref{app:geom_recoll}, based on geometric estimates.

\subsubsection{Convergence of the integrated cumulants} \label{sec:conv_cumulants}

\begin{figure}[h!] 
\centering
\begin{tikzpicture}
\begin{scope}[xshift = 0cm]
\draw[black!75, very thick ] (1, 0) -- (1, -1.5); 
\draw[red] (1, 0) circle (1 mm) node[above = 1mm]{$1$};
\draw[black!75, very thick ] (1, -1.5) -- (3, -6); 
\draw[red] (3, -6) circle (1 mm) node[below = 1mm]{$1$};
\draw[black!75, very thick ] (1, -1.5) -- (-1, -6); 
\draw[blue] (-1, -6) circle (1 mm) node[below = 1mm]{$2$};
\draw[black!75, very thick ] (0, -3.75) -- (1, -6); 
\draw[blue] (1, -6) circle (1 mm) node[below = 1mm]{$3$};
\draw[black, very thick, dotted ] (0.8, -1.5) -- (-1, -1.5) node[left = 2.5mm]{$\tau_{1}$}; 
\draw[black, very thick, dotted ] (-0.2, -3.75) -- (-1, -3.75) node[left = 2.5mm]{$\tau_{2} $};  
\end{scope}
\begin{scope}[xshift = 7cm]
\draw[black!75, very thick ] (1, 0) -- (1, -1.5); 
\draw[red] (1, 0) circle (1 mm) node[above = 1mm]{$1$};
\draw[black!75, very thick ] (1, -1.5) -- (4, -6); 
\draw[red] (4, -6) circle (1 mm) node[below = 1mm]{$1$};
\draw[black!75, very thick ] (1, -1.5) -- (-2, -6); 
\draw[blue] (-2, -6) circle (1 mm) node[below = 1mm]{$2$};
\draw[black!75, very thick ] (2, -3) -- (0, -6); 
\draw[blue] (0, -6) circle (1 mm) node[below = 1mm]{$3$};
\draw[black!75, very thick ] (2.5, -3.75) -- (1, -6); 
\draw[blue] (1, -6) circle (1 mm) node[below = 1mm]{$4$};
\draw[black, very thick, |-> ] (4.5, -6) node[right = 1mm]{$t=0$} -- (4.5, 0.5); 
\draw[black, very thick, dotted ] (1.2, -1.5) -- (4.3, -1.5) node[right = 2.5mm]{$\tau_{1} = \tau_{\{1,2\}}$}; 
\draw[black, very thick, dotted ] (2.2, -3) -- (4.3, -3) node[right = 2.5mm]{$\tau_{2} = \tau_{\{1,3\}}$}; 
\draw[black, very thick, dotted ] (2.7, -3.75) -- (4.3, -3.75) node[right = 2.5mm]{$\tau_{3} = \tau_{\{1,4\}}$}; 
\end{scope}
\end{tikzpicture}
\caption{Left: example of non-contributing dynamics tree in the expansion~\eqref{lgn:cumexp} of the first correlation function. Right: only contributing tree for 3~added particles.}
\label{fig:limtree}
\end{figure}
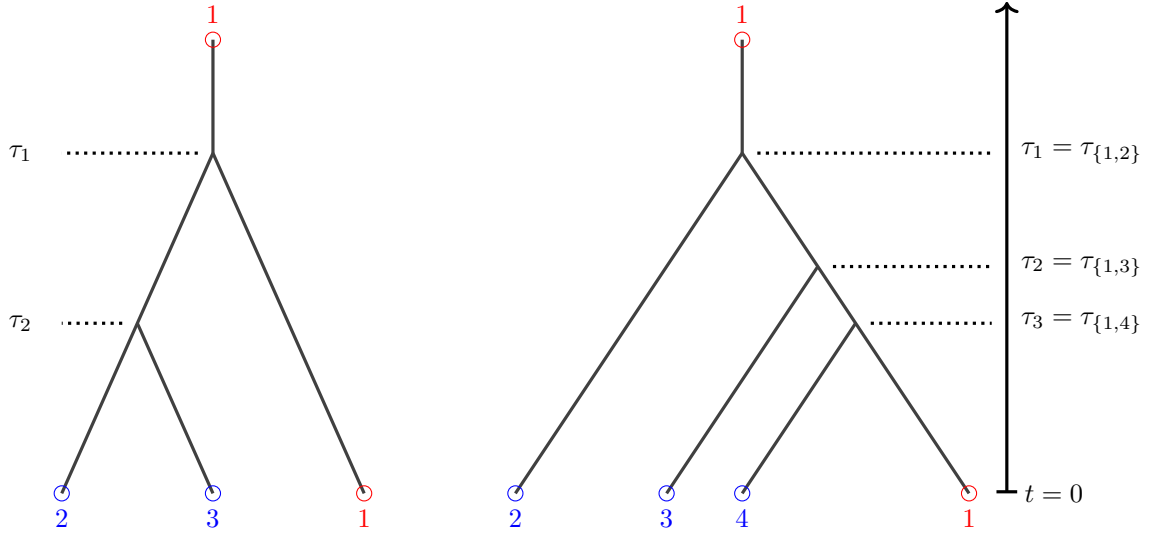

Now that we stated that the cycles and the other tagged particles asymptotically have a negligible impact on the dynamics of $n$~fixed particles, we go back to the expression~\eqref{def:In_rho_H} given in~Proposition~\ref{prop:dyn_trees} for the integral of cumulants, to determine their limit (Theorem~\ref{theo:cum_conv} below). Indeed, by Propositions~\ref{prop:rare_encounters2} and~\ref{prop:discard_cycles} above, up to an error of order~$(\varepsilon + p_\mu)$, the pseudo-trajectories~$\uz_{n+k}^{[0,t]}$ only include clustering encounters, and none of the added particles is tagged. Hence, the velocities of these pseudo-trajectories only depend on the dynamics tree~$T^\dd_{\lbr n \rbr, k}$ and on the dynamics parameters in the integral, so that they are equal to the velocities of the limit pseudo-trajectories~$\uze_{n+k}^{[0,t]}$ (see Section~\ref{sec:pseudo-traj}) . This way, similarly as in the companion paper~\cite{fou26Rayleigh}, the positions of both pseudo-trajectories are $(k\varepsilon)$-close, and in particular they are identical for the $n$~final particles. Since the only tagged particles are among these $n$~final particles, and since the equilibrium~$M_\beta$ only depends on the velocities, we get
\begin{equation*}
[\varphi_0^{\otimes \ul_n} M_\beta^{\otimes n+k}](\uze_{n+k}^{[0]}) = [\varphi_0^{\otimes \ul_n} M_\beta^{\otimes n+k}](\uz_{n+k}^{[0]}).
\end{equation*}
The observable~$H^{\otimes n}$ also only depends on the $n$~studied particles, so that one can bound
\begin{align*}
& \left\lvert H^{\otimes n}(\uze_{n}^{[0,t]})  [\varphi_0^{\otimes \ul_n} M_\beta^{\otimes n+k}](\uze_{n+k}^{[0]}) - H^{\otimes n}(\uz_{n}^{[0,t]})  F_{n+k}^\varepsilon(0,\uz_{n+k}^{[0]})  \right\rvert \\[0.7em]
& \hspace{10mm} = \left\lvert H^{\otimes n}(\uz_{n}^{[0,t]}) \right\rvert \cdot \left\lvert  [\varphi_0^{\otimes \ul_n} M_\beta^{\otimes n+k}](\uz_{n+k}^{[0]}) - F_{n+k}^\varepsilon(0,\uz_{n+k}^{[0]})   \right\rvert \\[0.7em]
& \hspace{20mm} \leq \left\lvert H^{\otimes n}(\uz_{n}^{[0,t]}) \right\rvert \cdot C_0^{n+k} e^{-\frac{\beta}{4}  \|\uv_{n+k}\|^2} \varepsilon,
\end{align*}
by the initial error estimate~\cite[Proposition~2.1]{fou26Rayleigh}. 
%% §§ numéro 2.1 ??
Eventually, recalling the result~\eqref{ineq:In_1_n} of Proposition~\ref{prop:cum_bound}, the leading term  in the cumulant expansion corresponds to the case of a single cluster $\rho_1 = \lbr 1, n \rbr$, once again with an error of order~$\varepsilon$. We end up with the following limit formula for the integrated cumulants.
\begin{theo}[Convergence of the cumulants] \label{theo:cum_conv}
For~$t > 0$ small enough, the integrated cumulants $\textstyle \int f^\varepsilon_n[H](t, \ul_n) $ converge as $\varepsilon$ goes to~0 in the scaling~\eqref{eq:scaling} towards the following limit formula
\begin{align}
& \int f_n[H](t, \ul_n)  \label{lgn:cumexp} \\
& \hspace{8mm} \doteq \sum_{k \geqslant 0} \sum_{\Upsilon \in \mT_{\lbr n\rbr, k}^{\dd}} \int \frac{\dd \overline{x} \dd \uv_{n+k}}{\mu^{n-1}}   \dd \uom_{E_\Upsilon} \dd \utau_{E_\Upsilon}  \prod_{e \in E_\Upsilon} s_{e} \langle v_{m}^{[{\tau_e}^+ ]}  - v_{m'}^{[{\tau_e}^+]},  \omega_e \rangle_+ H^{\otimes n}  [\varphi_0^{\otimes \ul_n} M_\beta^{\otimes n+k}](\uze_{n+k}^{[0]}), \nonumber
\end{align}
where none of the added particles in the pseudo-trajectories are tagged, and with the quantitative bound
\begin{align*}
 \left\lvert \int f_n[H](t, \ul_n) - \int f^\varepsilon_n[H](t, \ul_n) \right\rvert \leq \frac{C^n n!}{\mu^{n-1}} \| H \|_\infty^n\cdot \left( \varepsilon + p_\mu \right).
\end{align*}
\end{theo}
Let us observe finally that some of the trees in the formula above do not contribute to the sum. Indeed, imagine that in a tree two added particles (hence non-tagged ones), both meet as their first (forwards) encounter at a time~$\tau_e$. Then, changing the sign of the encounter leads to a weight
\begin{equation*} 
- s_{e} \langle v_{m}^{[{\tau_e}^+]}  - v_{m'}^{[{\tau_e}^+]},  \omega_e \rangle_+   M_\beta\left({v_m^{[0]}}'\right) M_\beta\left( {v_{m'}^{[0]}}'\right) = - s_{e} \langle v_{m}^{[{\tau_e}^+]}  - v_{m'}^{[{\tau_e}^+]},  \omega_e \rangle_+   M_\beta\left({v_m^{[0]}}\right) M_\beta\left( {v_{m'}^{[0]}}\right),
\end{equation*}
using the equilibrium structure.
Note that since they are added, they do not contribute to the weight~$\varphi_0^{\otimes \ul_n} H^{\otimes n}$, so that the contributions of this tree and of its counterpart with a changed sign cancel out. In particular, for the first limit cumulant~$f_1 = F_1$, as pictured in Figure~\ref{fig:limtree}, the only shape of tree that contributes is the linear one. Hence, it is enough to know the signs associated to each encounter, and one may write
\begin{align}
& \int F_1[H](t, \ell) \label{eq:F1_lim}\\
& \hspace{8mm} = \sum_{k \geqslant 0} \sum_{\us_{k}} \int \dd z_1 \dd \uv_{k}^*   \dd \uom_{k} \dd \utau_{k} \prod_{i=1}^{k} s_{i} \langle v_{1}^{[{\tau_i}^+ ]}  - v_{1+i}^{[{\tau_i}^+]},  \omega_i \rangle_+ H(\zeta_1^{[0,t]})  [\varphi_0^{\otimes \ell} M_\beta^{\otimes 1+k}](\uze_{1+k}^{[0]}). \nonumber
\end{align}
Fortunately, this expansion coincides with the limit pseudo-trajectory formulation~\eqref{eq:pseudo_traj_form} of the first correlation function, solution at the limit to the linear Rayleigh--Boltzmann equation~\eqref{eq:phi}.

\subsection{Convergence of the cumulant generating function} \label{sec:cumulant_generating_function_conv}

We study here the convergence of the cumulant generating function in the case~\eqref{def:cumulant_generating_function_marq} of an observable~$H = \widehat{H} \ind_{\ell=1}$ weighting only the tagged particles, generalized to observables depending on the pseudo-trajectory~$\uz^{[0,t]}$ in the whole time interval~$[0,t]$ as in~\eqref{def:cumulant_generating_function3}. It writes
\begin{align*}
\mathfrak{G}_\varepsilon^{[0,t]}[H]& =  \sum_{p \geqslant 1} \frac{\lambda^p}{p!} \int f^\varepsilon_{p}\left[e^H - 1 \right](t, \underline{1}_p).
\end{align*}
For $p\geq 2$, using the bounds on the cumulants (Proposition~\ref{prop:cum_bound}), one has
\begin{align*}
\sum_{p \geqslant 2} \frac{\lambda^p}{p!} \left\lvert \int f^\varepsilon_{p}\left[e^H - 1\right](t, \underline{1}_p) \right\rvert & \leq \lambda \sum_{p \geqslant 2} \left( \frac{\lambda}{\mu} \right)^{p-1} \left( C \| e^H - 1 \|_\infty \right)^p \\
& \leq  \lambda \times 2 \frac{\lambda}{\mu}  \left( C \| e^H - 1 \|_\infty \right)^2
\end{align*}
as soon as $ 2 \lambda C \| e^H - 1 \|_\infty \leq \mu$ in the scaling~\eqref{eq:scaling}. Hence, rescaling properly the generating function, we get the following proposition thanks to the convergence of the first cumulant (Theorem~\ref{theo:cum_conv}).
\begin{prop}[Convergence of the cumulant generating function] \label{prop:cum_gen_fun_conv} For any time~$t$ smaller than the short time~$T_\beta$ considered in Proposition~\ref{prop:cum_bound}, still in the mixed scaling~\eqref{eq:scaling}, the rescaled cumulant generating function converges as
\begin{equation} \label{eq:conv_G}
 \frac{1}{\lambda} \mathfrak{G}_\varepsilon^{[0,t]}[H] \xrightarrow[\varepsilon \to 0]{}  \int F_1[e^H - 1](t, 1) = \int F_1[e^H](t, 1) - 1. 
\end{equation}
\end{prop}
All the cumulants of order greater than~2 vanish in this scaling, at the opposite of the non-linear version~\cite{2023grandev} in which they contribute to add non-linearity in the limiting generating function.  Next section is dedicated to the equations driving this limit.

\subsection{Hamilton--Jacobi equations for the limit cumulant generating function} \label{sec:HJ}
For observables that behave according to the structure of the trajectories, i.e. of the form 
\begin{equation*}
H(z^{[0,t]}) = g(t, z^{[t]}) - \int_0^t [\partial_s - v \cdot \nabla_x] g(s,z(s)) \dd s
\end{equation*}
as introduced in Section~\ref{sec:stat_refinements}, we will study the limit cumulant generating function (see Section~\ref{sec:cumulant_generating_function_conv}), denoted
\begin{equation}  \label{def:Itg}
\mI(t,g) \doteq \int \tilde{F}_1\left[\exp\left( g_t - \int_0^t (\partial_s - v \cdot \nabla_x) g_s \right)\right](t).
\end{equation}
It is relevant to consider observables in the functional space
\begin{equation} \label{def:gset} 
g \in \mathbb{B}_{T, \beta} \doteq \left\{ h \ \Bigl\lvert \ \exists C > 0: \ \forall t \leqslant T, \ \left\|(h(t, v) - \frac{\beta}{4}|v|^2)_+ \right\|_{\Lp^\infty(\mD)} + \|(\partial_t - v\cdot \nabla_x)h(t)\|_{\Lp^\infty(\mD)} < C \right\}.
\end{equation}
Indeed, the relevant variable that will appear as boundary condition of a Boltzmann equation will be~$e^g(t)$ (see the proof of Proposition~\ref{prop:identifHJ}).
Note that the time of validity of Theorem~\ref{theo:LD} depends on this norm~$C$ chosen for $g \in \mathbb{B}_{T,\beta}$. One may choose to extend this class of functions to observables growing as small inverse Gaussians, using part of the exponential decay of the cumulants to compensate this growth (see~\cite[Chapter~7]{2023grandev}), yet that would mean losing control on the solutions to the Hamilton-Jacobi system below (Proposition~\ref{prop:HJSyst}).

We will use the formula~\eqref{eq:F1_lim} above to find the Hamilton--Jacobi equation it satisfies. Observe that the term~$k = 0$ in this expansion corresponds to the initial value~$\mI(0,g)$. For any~$k \geqslant 1$, we will split the dynamics tree in two, to isolate the influence of particle~$2$ at time~$\tau_1$. To this extent, let us define
$\widetilde{\us}_{k-1} \doteq (s_2, \dots, s_k)$, and similarly $\widetilde{\uv}_{k-1}, \widetilde{\uom}_{k-1}, \widetilde{\utau}_{k-1}$, and $\widetilde{\uze}_{k} \doteq (\zeta_1, \zeta_3, \dots, \zeta_{1+k})$.
We split according to the sign~$s_1$, which labels the first collision (between~$v_1$ and~$v_2 = v^*_1$), writing
\begin{align}
& \mI(t,g) - \mI(0,g) \label{lgn:Itg1}\\
& \hspace{10mm} =  \int \dd z_1 \dd v_2 \dd \omega_1 \dd \tau_{1} \langle v_{1}  - v_{2},  \omega_1 \rangle_+ M_\beta(v_2') \sum_{k \geqslant 1} \sum_{\widetilde{\us}} \int \dd \widetilde{\uv}_{} \dd \widetilde{\uom}_{} \dd \widetilde{\utau}_{} \prod_{i=2}^{k} s_{i} \langle v_{1}^{[{\tau_i}^+ ]}  - v_{1+i}^{[{\tau_i}^+]},  \omega_i \rangle_+  \nonumber \\
& \hspace{35mm}  \times \exp\left({g(t,\zeta_1^{[t]}) - \int_0^t (\partial_s - v \cdot \nabla_x)g(s,\zeta_1^{[s]}) }\right)  [\varphi_0 M_\beta^{\otimes k}](\widetilde{\uze}_{k}^{[0]}) \nonumber \\
& \hspace{10mm} - \int \dd z_1 \dd v_2 \dd \omega_1 \dd \tau_{1} \langle v_{1}  - v_{2},  \omega_1 \rangle_+ M_\beta(v_2) \sum_{k \geqslant 1} \sum_{\widetilde{\us}} \int \dd \widetilde{\uv}_{} \dd \widetilde{\uom}_{} \dd \widetilde{\utau}_{} \prod_{i=2}^{k} s_{i} \langle v_{1}^{[{\tau_i}^+ ]}  - v_{1+i}^{[{\tau_i}^+]},  \omega_i \rangle_+ \nonumber\\
& \hspace{35mm} \times \exp\left({g(t,\zeta_1^{[t]}) - \int_0^t (\partial_s - v \cdot \nabla_x)g(s,\zeta_1^{[s]}) }\right)  [\varphi_0 M_\beta^{\otimes k}](\widetilde{\uze}_{k}^{[0]}), \nonumber 
\end{align}
where the pseudo-trajectory~$\widetilde{\uze}_{k}^{[0,t]}$ corresponds to the scattering case~$s_1 = 1$ in the first term, and to the overlapping case~$s_1 = -1$ in the second one.
Now, by the fundamental theorem of calculus for the transport equation, since there is no collision in the pseudo-trajectory on~$[\tau_1, t]$, one has
\begin{equation*}
g(t,\zeta_1(t)) - \int_0^t (\partial_s - v \cdot \nabla_x) g(s,\zeta_1(s)) = g(\tau_1,\zeta_1(\tau_1^+)) - \int_0^{\tau_1} (\partial_s - v \cdot \nabla_x) g(s,\zeta_1(s)).
\end{equation*}
In the overlapping case, we have $g(\tau_1,\zeta_1(\tau_1^+)) = g(\tau_1,x_1^{[\tau_1]}, v_1) = g(\tau_1,\zeta_1(\tau_1^-))$, yet there is a discontinuity in the scattering case. Indeed, in the scattered case, writing
\begin{equation*}
g(\tau_1,\zeta_1(\tau_1^+)) = g(\tau_1,x_1^{[\tau_1]}, v_1) - g(\tau_1,x_1^{[\tau_1]}, v_1') + g(\tau_1,x_1^{[\tau_1]}, v_1'),
\end{equation*}
and recalling that the mapping $(v_1', v_2') \mapsto (v_1, v_2)$ is of Jacobian~1, the first term in the expansion~\eqref{lgn:Itg1} above writes exactly as the second one, with an additional factor
\begin{equation*}
\exp\left( g(\tau_1,x_1^{[\tau_1]}, v_1') - g(\tau_1,x_1^{[\tau_1]}, v_1) \right).
\end{equation*}
This way, the expansion~\eqref{lgn:Itg1} of the cumulant generating function corresponds to an integral over time of the same expansion at time~$\tau_1$, with an additional weight depending on~$z_1(\tau_1)$: this is precisely the partial derivative of the functional~$\mI(t,g)$ with respect to~$g(t)$ as an independent variable, in the direction of the said weight (the functional~$\mI(t,g)$ is clearly analytical in this variable~$g(t)$ by the expansion that we used above). Hence, one can write
\begin{align*}
& \mI(t,g) - \mI(0,g) =  \int_0^t \dd \tau \frac{\partial \mI}{\partial g(t)}(\tau,g)\left[\int \dd z_1 \dd v_2 \dd \omega \langle v_{1}  - v_{2},  \omega \rangle_+ M_\beta(v_2) (e^{g(\tau,z_1') - g(\tau,z_1)} - 1) \right].
\end{align*}
One may identify the partial derivative with the function~$z \mapsto \partial_{g(t)} \mI(t, g, z)$ such that
\begin{equation*}
\frac{\partial \mI}{\partial g(t)}(t,g)\left[h \right] = \int_{\mD} \partial_{g(t)} \mI(t, g, z) h(z) \dd z.
\end{equation*}
The function $z \mapsto \partial_{g(t)} \mI(t, g, z)$ is shown~\cite{2023grandev} to be continuous in~$x$, with values in the space of measures weighted by the inverse of the Maxwellian~$M_\beta^{-1}$. Hence, we are left with  
\begin{align} \label{eq:Hamilton_mild}
\mI(t,g) = \mI(0,g) + \int_0^t \dd \tau \int \dd z_1 \dd v_2 \dd \omega \langle v_{1}  - v_{2},  \omega \rangle_+ M_\beta(v_2) \partial_{g(t)} \mI(\tau, g, z_1) (e^{g(\tau,z_1') - g(\tau,z_1)} - 1),
\end{align}
yielding the following proposition.
\begin{prop}[Hamilton--Jacobi system for the limit cumulant generating function] \label{prop:HJSyst} 
Introducing the Hamiltonian
\begin{equation*}
 \mH\left( q , p \right) \doteq \int \dd z_1 \dd v_2 \dd \omega \langle v_{1}  - v_{2},  \omega \rangle_+ M_\beta(v_2) q(z_1) (e^{p(x_1, v_1) - p(x_1, v_1')} - 1),
\end{equation*}
for any time $t < T_\beta$ (see Propositions~\ref{prop:cum_bound} and~\ref{prop:cum_gen_fun_conv}), the limit cumulant generating function is a mild solution of the equation
\begin{equation*} 
\partial_t \mI(t,g) = \mH\left( \partial_{g(t)} \mI(t,g), g(t)\right),
\end{equation*}
in the sense of~\eqref{eq:Hamilton_mild}.
\end{prop}
At fixed~$t > 0$ and~$g \in \mathbb{B}_{t, \beta}$ (defined above in~\eqref{def:gset}), this Hamiltonian equation incites to introduce the following Hamilton--Jacobi system
\begin{equation} \label{eq:HJsystem}
\left\{ \begin{array}{lll}
(\partial_s - v\cdot \nabla_x) q^{[t]} = \frac{\partial \mH}{\partial p} (q^{[t]}, p^{[t]})  & , & q^{[t]}(0) = M \varphi_0 \exp(p^{[t]}(0)), \\
(\partial_s - v\cdot \nabla_x)(p^{[t]} - g) = -\frac{\partial \mH}{\partial q}(q^{[t]}, p^{[t]}) & , & p^{[t]}(t) = g(t),
\end{array} \right.
\end{equation}
where the unknowns~$(q^{[t]}, p^{[t]})$ are meant to find the minimizing values of~$\left(\partial_{g(t)}\mI(t,g) , g(t)\right)$ for this Hamiltonian. Next proposition is dedicated to proving that the mild Hamiltonian solution 
\begin{align} \label{eq:Ichapo}
\hat{\mI}(t,g) \doteq \mI(0,g) + \int_0^t \dd s \int_{\mD} q^{[t]}(s) (\partial_s - v\cdot \nabla_x)(p^{[t]}(s) - g(s))  + \int_0^t  \mH(q^{[t]}(s), p^{[t]}(s)) \dd s
\end{align}
is well-defined, and to identify it with the functional~$\mI(t,g)$.

\begin{prop}[Identification of the Hamiltonian solutions] \label{prop:identifHJ} For any time~$t>0$ and any observable~$g \in \mathbb{B}_{t, \beta}$, the Hamilton--Jacobi system~\eqref{eq:HJsystem} admits a unique global solution~$(q^{[t]}, p^{[t]})$ such that $(q^{[t]}e^{-p^{[t]}}, e^{p^{[t]}}) \in \Lp^\infty\left([0,t], \mF_{1, \beta/4} \right)^2$, and the functional~$\hat{\mI}(t,g)$ defined as~\eqref{eq:Ichapo} from this solution coincides on~$[0,t]$ with our functional
\begin{equation*} 
\mI(t,g) = \hat{\mI}(t,g).
\end{equation*}
\end{prop}

\begin{prof}
To prove the well-posedness of the equation, we compute the change of unknowns
\begin{equation*} 
\left\{ \begin{array}{l}
\gamma(s) = e^{g(s)} \\
\theta(s) = (\partial_s - v\cdot \nabla_x) g(s),
\end{array} \right.
\end{equation*}
and we consider the functional~$\mJ(t, \theta, \gamma) \doteq \mI(t, g)$.
For this functional, the formula for the partial derivative with respect with~$\gamma$ is even simpler: where the derivation of the term~$\gamma = e^g$ with respect to~$e^g$ let it invariant, the derivation with respect with~$\gamma$ makes it disappear. Thus, in the previous formulas, we have to add it back in the weight corresponding to the direction of the derivation. The associated Hamiltonian is hence
\begin{equation} \label{def:hatH}
 \hat{\mH}\left( \chi , \eta \right) \doteq \int \dd z_1 \dd v_2 \dd \omega \langle v_{1}  - v_{2},  \omega \rangle_+ M_\beta(v_2) \chi(z_1) (\eta(z_1') - \eta(z_1)),
\end{equation}
for the variables~$\chi \doteq q e^{-p}$ and~$\eta \doteq e^p$, and the Hamilton--Jacobi system writes as
\begin{equation*}
\left\{ \begin{array}{lll}
(\partial_s - v\cdot \nabla_x) \chi + \theta \chi  = \frac{\partial \hat{\mH}}{\partial \eta} (\chi, \eta)  & , & \chi(0) = M_\beta \varphi_0, \\
(\partial_s - v\cdot \nabla_x)\eta - \theta \eta  = -\frac{\partial \hat{\mH}}{\partial \chi}(\chi, \eta) & , & \eta(t) = \gamma(t).
\end{array} \right.
\end{equation*}
This system is equivalent to the following linear Boltzmann--Hamilton--Jacobi system
\begin{equation} \label{eq:modBoltz}
\left\{ \begin{array}{lll}
(\partial_s - v\cdot \nabla_x) \chi = - \theta \chi  + \int \dd v_2 \dd \omega \langle v  - v_{2},  \omega \rangle_+  \Bigl(M_\beta(v_2')\chi(z') - M_\beta(v_2)\chi(z)\Bigr) \\[1em]
(\partial_s - v\cdot \nabla_x)\eta = + \theta \eta  - \int \dd v_2 \dd \omega \langle v  - v_{2}, \omega \rangle_+ M_\beta(v_2) \Bigl(\eta(z') - \eta(z)\Bigr).
\end{array} \right.
\end{equation}
 Contrary to~\cite{2023grandev}, these equations are decoupled, since the coupling was stemming from the non-linearity. Moreover, they are not symmetrical: these linear equations are very close to the Rayleigh--Boltzmann equation~\eqref{eq:phi}, yet the first one has the same collision kernel as the equation on~$M_\beta \varphi$, and the second one the same as the equation on~$\varphi$. Appendix~\ref{app:RB} is dedicated to showing that in our functional setting, this system admits a unique global positive solution~$(\chi, \eta) \in \Lp^\infty\left([0,t], \ \mF_{1, \beta/4} \times \mF_{1, -\beta/4} \right)$, proving the first part of the proposition.

To perform the identification, we start to show it on small times. We prove that
\begin{equation*} 
\mJ(t, \theta, \gamma) = \hat{\mJ}(t, \theta, \gamma) \ \Bigl( \doteq \hat{\mI}(t,g) \Bigr),
\end{equation*}
where the second term is defined through the same change of unknowns. 
Both functionals are mild solutions to the same Hamilton--Jacobi equation (see for example~\cite{2023grandev} for the algebraic details). It is hence enough to show that there exists a unique solution to this mild equation in a regularity space common to both functionals.

Denoting $\mB_{R,\beta}$ the $(\mF_{1, - \beta/4})$-ball of radius~$R>0$, let us define the norm
\begin{equation*} 
\triple{ \mJ(t)}_{R, \beta} \doteq \sup_{\substack{\|\theta\|_\infty \leqslant R \\ \gamma \in \mB_{R, \beta}}} \lvert \mJ(t, \theta, \gamma) \rvert,
\end{equation*}
and let us assume the regularity condition, for every~$0 < \beta < \beta'$,
\begin{equation} \label{eq:regularityJ} 
\forall G \in \mC^0(\mD), \left[ \forall z \in \mD, |G(z)| \leq R C_\beta |v| e^{\frac{\beta}{4}|v|^2}  \right] \Rightarrow  \left\lvert  \int \partial_{\gamma} \mJ(t, \theta, \gamma) \ G \right\rvert \leq \frac{R \tilde{C}_\beta}{\sqrt{\beta' - \beta}} \triple{ \mJ(t)}_{R, \beta'},
\end{equation}
for some explicit constants~$C_\beta, \tilde{C}_\beta$ depending continuously on~$\beta$.
Then, considering the Hamilton--Jacobi equation~\eqref{eq:Hamilton_mild} and the modified Hamiltonian~$\hat{\mH}$ defined in~\eqref{def:hatH}, we have for any functionals~$\mJ, \mJ'$ satisfying the regularity asumption~\eqref{eq:regularityJ} above, and~$\gamma$ in the ball~$\mB_{R, \beta}$, that
\begin{align} \label{ineq:HJmoinsJ'}
\left\lvert \hat{\mH}\bigl(\partial_{\gamma}\mJ, \gamma\bigr) - \hat{\mH}\bigl(\partial_{\gamma}\mJ', \gamma\bigr) \right\rvert 
&  \leq \left\lvert \int \dd z \dd v_2 \dd \omega \langle v  - v_{2},  \omega \rangle_+ M_\beta(v_2) \Bigl[ \partial_{\gamma}\mJ - \partial_{\gamma}\mJ' \Bigr](z) \Bigl|\gamma(z') - \gamma(z)\Bigr|\right\rvert.
\end{align}
Using the fact that~$\gamma\in \mB_{R, \beta}$, we have on the one hand for any~$z \in \mD$,
\begin{align*}
\left\lvert \int \dd v_2 \dd \omega \langle v  - v_{2},  \omega \rangle_+ M_\beta(v_2) | \gamma(z)| \right\rvert
& \leq R e^{\frac{\beta}{4}|v|^2} \nu_\beta(v) \\
& \leq R e^{\frac{\beta}{4}|v|^2} C_\beta |v|,
\end{align*}
with the notation~\eqref{def:nubeta} for the loss factor~$\nu_\beta(v)$. On the other hand, using a variant of Lemma~\ref{lemm:opkern2} to rewrite the collision operator in terms of its integral kernel~\cite[Lemma~3.3.1]{fouthese}, we observe that
\begin{align*}
\left\lvert \int \dd v_2 \dd \omega \langle v  - v_{2},  \omega \rangle_+ M_\beta(v_2) | \gamma(z')| \right\rvert
& \leq R \int \dd v_2 \dd \omega \langle v  - v_{2},  \omega \rangle_+ M_\beta(v_2) e^{\frac{\beta}{4} |v'|^2} \\
& \leq R C_\beta \int \frac{\dd \eta\ e^{\frac{\beta}{4}|\eta|^2}}{|\eta - v|^{d-2}}  e^{-\beta \left[\frac{|\eta - v|^2}{8} + \frac{|\eta|^2 - |v|^2}{4} + \frac{\bigl( |\eta|^2- |v|^2\bigr)^2 }{8 |\eta - v|^2 } \right] }  \\
& \leq R C_\beta e^{\frac{\beta}{4}|v|^2} \int \frac{\dd \eta}{|\eta - v|^{d-2}}  e^{-\beta \left[\frac{|\eta - v|^2}{8} + \frac{\bigl( |\eta|^2- |v|^2\bigr)^2 }{8 |\eta - v|^2 } \right] },
\end{align*}
so that we recognize precisely the collision kernel of the modified collision operator~\eqref{def:Kgain2}. Thanks to its bound shown in Appendix~\ref{app:RB}, Lemma~\ref{lemm:opkern2}, we get
\begin{align*}
\left\lvert \int \dd v_2 \dd \omega \langle v  - v_{2},  \omega \rangle_+ M_\beta(v_2) | \gamma(z')| \right\rvert
& \leq \overline{C}_\beta R \frac{ e^{\frac{\beta}{4}|v|^2}}{1 + |v|}.
\end{align*}
In the end, the functional in~$z$ integrated against~$\Bigl[ \partial_{\gamma}\mJ - \partial_{\gamma}\mJ' \Bigr]$ in~\eqref{ineq:HJmoinsJ'} satisfies the assumptions of the regularity condition~\eqref{eq:regularityJ}, so that 
\begin{align*}
\left\lvert \hat{\mH}\bigl(\partial_{\gamma}\mJ, \gamma\bigr) - \hat{\mH}\bigl(\partial_{\gamma}\mJ', \gamma\bigr) \right\rvert 
&  \leq \frac{\tilde{C}_\beta}{\sqrt{\beta' - \beta}} \triple{ \mJ(t)}_{R, \beta'}.
\end{align*}
This inequality is precisely the hypothesis of the abstract Cauchy--Kovalevskaya theorem stated in~\cite[Appendix~A]{2023grandev}, for equations of the form~\eqref{eq:Hamilton_mild} that we consider. According to this theorem, there exists a short time on which there exists a unique solution to the mild equation on~$\mJ$, with the regularity~\eqref{eq:regularityJ} that we imposed. To show that both functionals~$\mJ$ and~$\hat{\mJ}$ are identical on short times, it is now enough to show that they share this regularity assumption. It is trivial in the case of~$\hat{\mJ}$ since by construction we have
\begin{equation*} 
\partial_{\gamma}\hat{\mJ}(t, \theta, \gamma) = \chi(t).
\end{equation*}
For the other one, we use its analyticity in~$\gamma$ to write for any~$\lambda \in \R$ that
\begin{equation*} 
\int \partial_{\gamma}\mJ(t, \theta, \gamma)(z) G(z) \dd z = \frac{1}{2 \pi\lambda} \int_0^{2 \pi} \mJ(t, \theta, \gamma + e^{i\theta} \lambda G) e^{-i\theta} \dd \theta.
\end{equation*}
Since we assumed that for any~$z\in \mD$,
\begin{align*}
G(z)&  \leq R C_\beta |v|e^{\frac{\beta}{4}|v|^2} \\
& \leq R \frac{ \tilde{C}_\beta}{\sqrt{\beta' - \beta}} e^{\beta' |v|^2}, 
\end{align*}
if we choose
\begin{equation*} 
 \lambda = \frac{\sqrt{\beta' - \beta}}{R \tilde{C}_\beta}
\end{equation*}
we obtain
\begin{equation*} 
\gamma + e^{i\theta} \lambda G \in \mB_{R, \beta'}.
\end{equation*}
Eventually, thanks to the analytical formula, we retrieve exactly the wanted regularity domination~\eqref{eq:regularityJ}.
Now that the identification is established on short times, on any time interval~$[0,t]$, one has uniform bounds on~$\hat{\mJ}$ thanks to the long-time existence result on the linear Boltzmann--Hamilton--Jacobi system. Hence, after the first small time of identification, since $\mJ = \hat{\mJ}$, both functionals still belong to the same functional framework, allowing to extend the proof of the uniqueness until any fixed large time~$t$, concluding the proof. 
\CQFD
\end{prof}

\section{Geometric control of dynamics cycles} \label{app:geom_recoll}

We prove here Proposition~\ref{prop:discard_cycles}, adding the analysis of geometric cycle constraints to the bounds proved in Section~\ref{sec:bounds}. Thanks to a fine computation to handle the appearing singularities (Section~\ref{sec:singularities}), we achieve a full convergence rate in a full power of~$\varepsilon$, better than the previous existing one~\cite{2023longcor}.

\subsection{Parametrizing the cycle}

First of all, it will be useful in the following to parametrize the encounters in terms of the angle~$\sigma$ instead of~$\omega$ (see Appendix~\ref{app:scattering} for the associated change of variables).
For a general cluster~$R = \rho_j$, at fixed number of added particles~$k_j$ and fixed labels, we want to show that the cycle condition in the following integral
\begin{equation} \label{cycle_integral}
\sum_{\Upsilon \in \mT_{R, k}^{\dd}} \int \ind_{T^\dd_{R} = \Upsilon} \dd \uv_{r+k}   \dd \usig_{E_\Upsilon} \dd \utau_{E_\Upsilon}  \prod_{e \in E_\Upsilon} s_{e} | v_{m}({\tau_e}^+ )  - v_{m'}({\tau_e}^+)| e^{- \beta \|\uv_{r+k}\|^2}\times \ind^\between_{R}
\end{equation}
leads to similar bounds as without the cycle constraint (see~\eqref{lgn:dom_SR}), with an additional small factor~$\varepsilon$. This way, we will be able to use the same computation as in Section~\ref{sec:bounds}, when bounding the integrated cumulants, gaining smallness from the cycle.

We have denoted~$\uv_{r+k} = (\uv_R, \uv_k^*)$, and we define $\bV = \max( 1, \| \uv_{r+k} \|_2 )$ to control the total energy. For simplicity, we will not precise in all the following, every time we use an adverb of time, that \emph{time is looked at backwards}.

Let us hence suppose that two particles, say $i$ and $j$, create the first cycle at a time $t_\cy \in [\tau_{c +1}, \tau_c]$, through a non-clustering collision or overlap. The cycling condition imposes strong geometric constraints, providing smallness when integrating over well-chosen parameters in the dynamics, restrained to small geometric volumes. Nevertheless, to integrate over a collision angle in~\eqref{cycle_integral}, we need to make sure that the velocities appearing in the product over the edges do not depend on this angle: we will hence start identifying the corresponding edges, and sum over them to make the associated velocities disappear like in~\eqref{ineq:sum_edge_1} and~\eqref{ineq:sum_edge_2}.

\begin{figure}[h!]\centering 
\begin{tikzpicture}
\begin{scope}[scale = 0.9]
\draw[black, thick, -{Latex[length=4mm, fill = white]}] (1,0) -- (3*49/50,4*49/50);
\draw[black] (3/2,4/2) node[left= 1mm]{$u_j$};
\draw[black, thick, fill = black] (3,4) circle (1 mm);
\draw[black, thick, fill = black] (3.2,4) circle (1 mm);
\draw[black, thick, fill = black] (4.2,2) circle (1 mm);
\draw[black, thick, fill = black] (4.4,2) circle (1 mm);
\draw[black, thick] (3.2,0) -- (3.2 + 1,2); % ui
\draw[black, thick] (3.2 + 1/2,2/2) node[left]{$u_i$};
\draw[black, thick, -{Latex[length=4mm, fill = white]}] (4.2,2) -- (4.2 - 1*21.36/22.36,2 + 2*21.36/22.36); %        ui'
\draw[black] (4.2 - 1/2,2 + 2/2)node[right]{$u_i'$};
\draw[black, thick, -{Latex[length=4mm, fill = white]}] (4.4,2) -- (5.7,3.4) node[right]{$u_k'$};
\draw[black, thick] (5.4,0) -- (5.4 - 1,2); % uk
\draw[black](5.4 - 1/4,2/4) node[right]{$u_k$}; 
\draw[black, very thin, loosely dashed] (1,0) -- (7 , 0) node[right]{$\tau_0$};
\draw[black, very thin, loosely dashed] (4.4,2) -- (7, 2) node[right]{$\tau_\dd$};
\draw[black, very thin, loosely dashed] (3.2,4) -- (7 , 4) node[right]{$t_\cy$};
\draw[black, thick, fill = black!10, opacity = 0.6, rounded corners=8pt] (0.7, 0) rectangle (3.5, -1);
\draw[black, thick, fill = black!10, opacity = 0.6] (2.1, -0.5) node{$C_i$};
\draw[black, thick, fill = black!10, opacity = 0.6] (5.4, -0.5) circle (5mm) node{$C_k$};
\draw[black] (3.5,-1)  node[below]{ Deflection case~1 ($j$ stems from~$C_i$)};
\end{scope}
%%%%%%%%%%%%%%%%%%%%%%%%%%%%%%%%%%%%%%%%%%%%%%%%%%%%%%%%%%%%%%
%                                                            %
%%%%%%%%%%%%%%%%%%%%%%%%%%%%%%%%%%%%%%%%%%%%%%%%%%%%%%%%%%%%%%
\begin{scope}[xshift = 8cm, scale = 0.9]
\draw[black, thick, -{Latex[length=4mm, fill = white]}] (5.5,0) -- (5.5 - 3.5*49/50,4*49.5/50);
\draw[black] (6.5 - 4.5*0.3,4*0.3) node[right= 1mm]{$u_j$};
\draw[black, thick, fill = black] (2,4) circle (1 mm);
\draw[black, thick, fill = black] (1.8,4) circle (1 mm);
\draw[black, thick, fill = black] (1.5,2) circle (1 mm);
\draw[black, thick, fill = black] (1.7,2) circle (1 mm);
\draw[black, thick] (0.5,0) -- (0.5 + 1,2); % ui
\draw[black, thick] (0.5 + 1/2,2/2) node[left]{$u_i$};
\draw[black, thick, -{Latex[length=4mm, fill = white]}] (1.5,2) -- (2 - 0.2*21.36/22.36,2 + 2*21.36/22.36); %        ui'
\draw[black] (1.2 - 1/2, 2 + 2/2)node[right]{$u_i'$};
\draw[black, thick, -{Latex[length=4mm, fill = white]}] (1.7,2) -- (2.3,3) node[below = 3mm]{$u_k'$};
\draw[black, thick] (2.7,0) -- (2.7 - 1,2); % uk
\draw[black](2.7 - 1/4,2/4) node[right]{$u_k$}; 
\draw[black, very thin, loosely dashed] (0.5,0) -- (7 , 0) node[right]{$\tau_0$};
\draw[black, very thin, loosely dashed] (1.5,2) -- (7, 2) node[right]{$\tau_\dd$};
\draw[black, very thin, loosely dashed] (2,4) -- (7 , 4) node[right]{$t_\cy$};
\draw[black, thick, fill = black!10, opacity = 0.6] (0.5, -0.5) circle (5mm) node{$C_i$};
\draw[black, thick, fill = black!10, opacity = 0.6, rounded corners=8pt] (2.4, 0) rectangle (5.8, -1);
\draw[black, thick, fill = black!10, opacity = 0.6] (4.1, -0.5) node{$C_k$};
\draw[black] (3.5,-1)  node[below]{ Deflection case~2 ($j$ stems from~$C_k$)};
\draw[black, thin, loosely dotted, |->] (8.2,-1)  node[below]{\small $t$} -- (8.2,4.5) node[above]{\small $0$} ;
\end{scope}
%%%%%%%%%%%%%%%%%%%%%%%%%%%%%%%%%%%%%%%%%%%%%%%%%%%%%%%%%%%%%%
%                                                            %
%%%%%%%%%%%%%%%%%%%%%%%%%%%%%%%%%%%%%%%%%%%%%%%%%%%%%%%%%%%%%%
\begin{scope}[xshift = 0cm, yshift = -6.5cm, scale = 0.9]
\draw[black, thick, -{Latex[length=4mm, fill = white]}] (0,0) -- (3*49/50,4*49/50); %uj
\draw[black] (3/2,4/2) node[left= 1mm]{$u_j$};
\draw[black, thick, fill = black] (3,4) circle (1 mm);
\draw[black, thick, fill = black] (3.2,4) circle (1 mm);
\draw[black, thick] (4.2,0) -- (4.2 + 1/2,1); % ui
\draw[black, thick] (4.2 + 1/3,2/3) node[left]{$u_i$};
\draw[black, thick, fill = black] (4.7,1) circle (1 mm);
\draw[black, thick, fill = black] (4.9,1) circle (1 mm);
\draw[black, thick, -{Latex[length=4mm, fill = white]}] (4.7,1) -- (4.2 - 1*21.36/22.36,2 + 2*21.36/22.36); %         ui'
\draw[black] (4.2,2)node[right]{$u_i'$};
\draw[black, thick, -{Latex[length=4mm, fill = white]}] (4.9,1) -- (6,2.3) node[below = 2.5mm]{$u_k'$};
\draw[black, thick] (5.4,0) -- (5.4 - 1/2,1); % uk
\draw[black](5.4 - 1/4,2/4) node[right]{$u_k$}; 
\draw[black, thick, -{Latex[length=4mm, fill = white]}] (1,0) -- (1 + 4*55.56/56.56,4*55.56/56.56); %ul
%\draw[black, thick,densely dotted, opacity = 0.6] (0,0) -- (0.5, -0.5);
\draw[black, thick] (1+1,1) node[right = 1.5mm]{$u_l$};
\draw[black, very thin, loosely dashed] (0,0) -- (7 , 0) node[right]{$\tau_0$};
\draw[black, very thin, loosely dashed] (4.9,1) -- (7, 1) node[right]{$\tau_\dd$};
\draw[black, very thin, loosely dashed] (3.2,4) -- (7 , 4) node[right]{$t_\cy$};
\draw[black, very thin, loosely dashed] (3.8,2.8) -- (7 , 2.8) node[right]{$\tau_\ov$};
%\draw[black, thick, fill = black] (5,4) circle (1 mm);
%\draw[black, thick, fill = black!10, opacity = 0.6] (0.5, -0.5) circle (7.07mm) node{$C_l$};
\draw[black, thick, fill = black!10, opacity = 0.6, rounded corners=8pt] (-0.3, 0) rectangle (1.3, -1);
\draw[black, thick, fill = black!10, opacity = 0.6] (0.5, -0.5) node{$C_l$};
\draw[black] (3.5,-1)  node[below]{ Deflection case~3 ($j$ stems from~$C_l$)};
\end{scope}
%%%%%%%%%%%%%%%%%%%%%%%%%%%%%%%%%%%%%%%%%%%%%%%%%%%%%%%%%%%%%%
%                                                            %
%%%%%%%%%%%%%%%%%%%%%%%%%%%%%%%%%%%%%%%%%%%%%%%%%%%%%%%%%%%%%%
\begin{scope}[xshift = 8cm, yshift = -6.5cm, scale = 0.9]
\draw[black, thick, -{Latex[length=4mm, fill = white]}] (0,0) -- (3*49/50,4*49/50);
\draw[black] (3/2,4/2) node[left= 1mm]{$v_j$};
\draw[black, thick, fill = black] (3,4) circle (1 mm);
\draw[black, thick, fill = black] (3.2,4) circle (1 mm);
\draw[black, thick, -{Latex[length=4mm, fill = white]}] (5.2,0) -- (4.2 - 1*21.36/22.36,2 + 2*21.36/22.36); %       ui'
\draw[black] (4.8,1)node[right]{$v_i$};
\draw[black, very thin, loosely dashed] (0,0) -- (7 , 0) node[right]{$\tau_0$};
\draw[black, very thin, loosely dashed] (3.2,4) -- (7 , 4) node[right]{$t_\cy$};
\draw[black, very thin, loosely dashed] (3.8,2.8) -- (7 , 2.8) node[right]{$\tau_\ov$};
\draw[black, thick, -{Latex[length=4mm, fill = white]}] (1,0) -- (1 + 4*55.56/56.56,4*55.56/56.56);
\draw[black] (1 + 4/4,4/4) node[right = 2mm]{$u_k$};
%\draw[black, thick,densely dotted, opacity = 0.6] (0,0) -- (0.5, -0.5);
%\draw[black, thick, fill = black] (5,4) circle (1 mm);
\draw[black, thick, fill = black!10, opacity = 0.6] (0.5, -0.5) circle (7.07mm) node{$C_1$};
\draw[black, thick, fill = black!10, opacity = 0.6] (5.7, -0.5) circle (7.07mm) node{$C_2$};
\draw[black] (3.5,-1)  node[below]{ Overlap case};
\draw[black, thin, loosely dotted, |->] (8.2,-1)  node[below]{\small $t$} -- (8.2,4.5) node[above]{\small $0$} ;
\end{scope}
\end{tikzpicture}
\caption{The three deflection cases, and the overlap case} \label{fig:arroz} 
\end{figure}
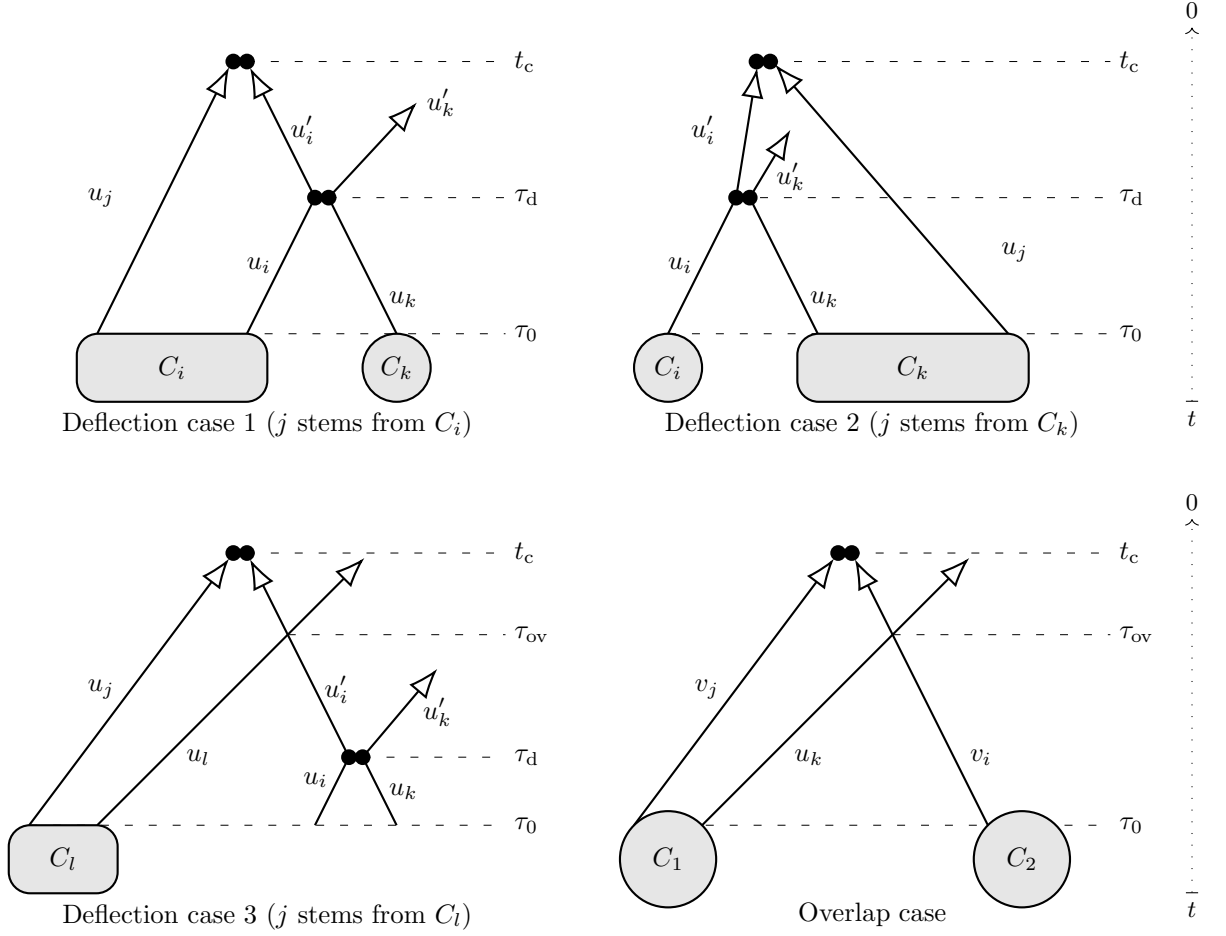

The parametrization of the cycle will depend on whether one of both particles~$i$ or~$j$ has undergone a deflection between $t$ and~$\tau_c$ or not. We hereafter define the relevant interactions to parametrize the cycle.
Should it be the case, let us denote~$k$ the particle that deflected~$i$ at the closest time~$\tau_\dd \geq \tau_c$, and if such a deflection never happened let us denote~$k$ the particle that overlapped~$i$ at a closest time~$\tau_\ov \geq \tau_c$. One can arbitrarily choose~$i$, among~$i$ and~$j$, to be the closest to have been deflected, or to have been overlapped if none of both has ever been deflected. 
Eventually, something that might happen is that the first connection between~$i$ and~$j$ (before their cycle) stems from~$i$ overlapping a fourth particle (call it~$l$), after having encountered~$k$ (see~$i,j,k,l$ in~Fig.~\ref{fig:arroz}). 

The choice of this configuration $(i,j,k,l) \in \lbr 1, r+k \rbr^4 $ and of $(\tau_c, \tau_\dd, \tau_\ov) \subset \{\tau_p \}_{p \leq r+k-1}$ corresponds to a combinatorial factor of order $(r+k)^7$.

Since the first cycle happens between~$\tau_c$ and~$\tau_{c+1}$, before the time~$\tau_c$ every encounter is clustering, registered in the tree~$\Upsilon$, and every velocity in the dynamics thus only depends on the signs and angles associated to the edges of~$\Upsilon$. This way, we may sum over these edges like in Section~\ref{sec:bounds}, as soon as we do it in an order respecting the interdependence among the encounters.

The set~$K \subset E_\Upsilon$ of edges possibly impacted by the collision between~$i$ and~$k$ at time~$\tau_\dd$ is defined as follows: 
we start with all the edges between~$\tau_{c+1}$ and~$0$, since the cycle happens before (as we look at time backwards, time~0 is \emph{after} the cycle occurs). Then, we can define the connected component of~$k$ when removing the edge~$\tau_\dd$ between~$i$ and~$k$, component of which we add to~$K$ all the edges after the deflection time~$\tau_\dd$ (see the definition of the dynamics tree in Fig.~\ref{fig:dyn_tree}). Eventually, we add to~$K$ all the potential overlaps undergone by~$i$ after~$\tau_\dd$, apart from the one with~$l$ if appropriate, since we need it to close the cycle. We hence condition on this set~$K$, which corresponds to a factor~$2^c$ once~$c$ is fixed, and we can sum beforehand over all these edges as in~\eqref{ineq:sum_edge_1} and~\eqref{ineq:sum_edge_2}. To be finally able to integrate over the deflection angle~$\sigma_\dd$ associated with the time~$\tau_\dd$ without impacting velocities in the product over the edges, it only remains to dominate roughly the one associated with the overlap time~$\tau_\ov$ in the deflection case, as follows:
\begin{equation*}
| v_{i}({\tau_\ov}^+ )  - v_{l}({\tau_\ov}^+)| \leq 2 \bV.
\end{equation*}
This factor will be found back in the final estimation~\eqref{lgn:final_est_cycle}. We will now treat both cases separately: as said before, if one of both particles has been deflected before $\tau_c$, we will use the angle~$\sigma_\dd$ of the corresponding collision to parametrize the cycle condition. Else, we will use the distance between the aggregates involved in the closest overlap.	
Both ways will raise some singularities in velocities, that we will control afterwards, in Section~\ref{sec:singularities}.

\subsection{Deflected case }

 If one of the particles has been deflected before $\tau_c$, we have fixed the closest deflection, involving $i$ and $k$ at time $\tau_{\mathrm{d}}$. 
Taking the encounter time just before $\tau_{\mathrm{d}}$ as a reference (call it $\tau_0$), we shall denote
$u_i, u_j, u_k$ the velocities of particles $i,j,k$ between $\tau_0$ and $\tau_{\mathrm{d}}$, and $u_i', u_j', u_k'$ these velocities between $\tau_{\mathrm{d}}$ and $t_\cy$ (cf.~Fig.~\ref{fig:arroz}). We first suppose that $k\neq j$. Then, there are three possibilites for the origin of particle $j$, considering that at time~$t_\cy$ it creates a cycle with $i$: it can stem from the connected component~$C_i$ of $i$ at time $\tau_0$, from the connected component~$C_k$ of $k$ at time $\tau_0$, or from the connected component~$C_l$ of a particle~$l$ overlapping~$i$ at a time~$\tau_{\ov}$, between~$\tau_\dd$ and $\tau_c$ (Fig.~\ref{fig:arroz}). The dependency of $x_i - x_j$ along the setting at~$\tau_0$ is different in each case.

\begin{itemize}[label = $\triangleright$]
\item If $j$ is in the connected component of~$i$ (Fig.~\ref{fig:arroz}, top left), this is the simplest case and
\begin{equation*}
[x_i - x_j](t_\cy) = [x_i -  x_j](\tau_0) + (\tau_\dd-\tau_0)[u_i - u_j] + (t_\cy - \tau_\dd)[u_i' - u_j].
\end{equation*}
\item If $j$ is in the connected component of~$k$ (Fig.~\ref{fig:arroz}, top right), then 
\begin{align*}
 [x_i - x_j] (t_\cy) & = [x_i -  x_k + x_k -  x_j](\tau_\dd) + (t_\cy - \tau_\dd)[u_i' - u_j] \\ 
& =\varepsilon \omega_\dd + [x_k - x_j](\tau_0) + (\tau_\dd-\tau_0)[u_k - u_j] + (t_\cy - \tau_\dd)[u_i' - u_j]. \nonumber
\end{align*}
\item Eventually, if $j$ is in the connected component of~$l$ (Fig.~\ref{fig:arroz}, bottom left), then 
\begin{align*}
[x_i - x_j](t_\cy) & = [x_i -  x_l + x_l -  x_j](\tau_\ov) + (t_\cy - \tau_\ov)[u_i' - u_j] \\
& = \varepsilon \omega_\ov + [x_l -  x_j](\tau_0) + (\tau_\ov - \tau_0)[u_l - u_j] +  (t_\cy - \tau_\ov)[u_i' - u_j]. \nonumber
\end{align*}
Actually, it may happen that $u_l$ is deflected between $\tau_0$ and $\tau_\ov$, yet it only requires to change the term $(\tau_\ov - \tau_0)[u_l - u_j]$, splitting it according to each segment between two deflections.
\end{itemize}
In each case, using the cycle condition $[x_i - x_j](t_\cy) = \varepsilon \omega_\cy + \zeta$, for a $\zeta \in \Z^d$ due to the periodicity of the domain and satisfying $|\zeta| \leq t \bV$, one can write
\begin{equation} \label{eq:cone}
(t_\cy - \tau_*)[u_i' - u_j] = \Delta x(\tau_*) + \varepsilon \omega_*,
\end{equation}
where $\tau_* \in \{ \tau_\dd, \tau_\ov \}$, $\Delta x(\tau_*)$ is independent of $\omega_\dd$ (see the explicit formula~(\ref{eq:deltx})), and $\omega_* \in \omega_\cy - \{0, \omega_\dd, \omega_\ov  \} $ satisfies $|\omega_*| \leq 2$.
This relationship~(\ref{eq:cone}) defines a cone to which $u_i' - u_j$ must belong to achieve the cycle, whose height depends on the relative position of~$i$ and~$j$ at time $\tau_*$: if $\Delta x(\tau_*)$ is large enough compared to $\varepsilon$, we will be able to control the height of this cone, and hence its angle, which constrains~$\sigma_\dd$. Otherwise, it will provide a strong condition on the encounter times. First, if 
\begin{equation} \label{ineq:deltax}
|\Delta x (\tau_*)| \geq 4 \varepsilon,
\end{equation}
then
\begin{equation*}
\Bigl |(t_\cy - \tau_*) [u_i' - u_j]\Bigr | \geq \frac{|\Delta x (\tau_*)|}{2},
\end{equation*}
so that
\begin{equation} \label{ineq:tauc}
\frac{1}{|t_\cy - \tau_*|} \leq \frac{4 \mathbf{V}}{ |\Delta x(\tau_*)|}.
\end{equation}
Now, the value of $\Delta x(\tau_*)$ is respectively given by
\begin{equation} \label{eq:deltx} \def\arraystretch{1.5}
\Delta x(\tau_*) = \left\{ \begin{array}{ll}
{[} x_j - x_i](\tau_0) + (\tau_\dd-\tau_0)[u_j - u_i] & \hspace{9mm} \sif j \in C_i, \\
{[} x_k - x_j](\tau_0) + (\tau_\dd - \tau_0)[u_k - u_j] & \hspace{9mm} \sif j \in C_k, \\
{[} x_l - x_j](\tau_0) + (\tau_\ov - \tau_0)[u_l - u_j] &\hspace{9mm} \sif j \in C_l, \\
\end{array} \right.
\end{equation}
which is affine in $\tau_* \in \{ \tau_\dd, \tau_\ov\}$, so that denoting $(\tau_* - \tau_0^*)\Delta u $ the projection of $\Delta x(\tau_*)$ on $\Delta u \in \{u_j - u_i, u_k-u_j, u_l - u_j \}$ according to the considered case, one has $|\Delta x(\tau_*)|	 \geq |(\tau_* - \tau_0^*)\Delta u|$.
Thus, using~(\ref{eq:cone}) and~(\ref{ineq:tauc}), we know that
\begin{equation*}
u_i' - u_j = \frac{ \Delta x(\tau_*) + \varepsilon \omega_{*}}{|t_\cy - \tau_{*}|}
\end{equation*}
belongs to a cylinder of width at most \[ \frac{4 \varepsilon \mathbf{V}}{|\Delta x(\tau_*)|} \leq \frac{4 \varepsilon \mathbf{V}}{|(\tau_* - \tau_0^*) \Delta u|}, \] (which contains the previous cone of controlled height, cf Fig.~\ref{fig:geom}, \textbf{a.}). Hence, $u_i'$ also belongs to such a cylinder, since $u_j$ remains unchanged. Now, the deflection condition states that \[|u_i' - (u_i + u_k)/2| = |u_i - u_k|/2 , \] so that~$u_i' - (u_i + u_k)/2$, whose direction is given by~$\sigma_\dd$ (see Fig.~\ref{fig:omesig}), also has to belong to a sphere of radius $|u_i - u_k|/2$. The intersection of this sphere with the previous cylinder is included in the union of two cones, of which solid angle is maximum when the cylinder is tangent to the sphere (cf. Fig.~\ref{fig:geom},~\textbf{b.}), namely 
\begin{equation*}
C_d \min\left[ 1, \left(\frac{ \varepsilon \mathbf{V}}{|(\tau_* - \tau_0^*) \Delta u |\cdot |u_i - u_k|}\right)^{d-1}\right]. 
\end{equation*}
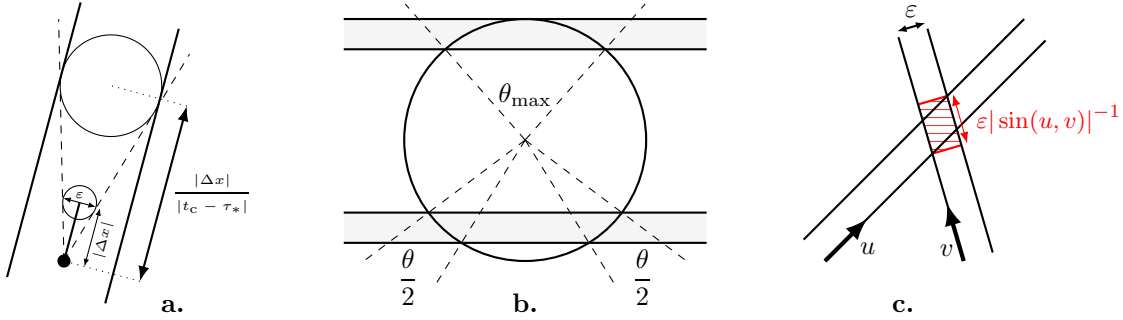
\begin{figure}[h!]\centering 
\begin{tikzpicture}
%%%%%%%%%%%%%%%%%%%%%%%%%%%%%%%%%%%%%%%%%%%%%%%%%%%%%%%%%%%%%%
%                     Pris dans la gouttière                 %
%%%%%%%%%%%%%%%%%%%%%%%%%%%%%%%%%%%%%%%%%%%%%%%%%%%%%%%%%%%%%%
\begin{scope}[xshift = 0cm, rotate = -15, scale=0.8]
\draw[black, thick] (0,0) -- (0, 1);
\draw[black, thin] (0,1) node[above = -0.7mm]{\tiny $\varepsilon$} circle (0.28);
\draw[black, thick] (-0.28*3,-0.5) -- (-0.28*3, 4.2);
\draw[black, thick] (0.28*3,-0.5) -- (0.28*3, 4.2); %(0.28*3/0.96, 4.2);
\draw[black, thin, {Latex[length=1mm, fill = black]}-{Latex[length=1mm, fill = black]}] (-0.28,1) -- (0.28, 1);
\draw[black, thin, {Latex[length=1mm, fill = black]}-{Latex[length=1mm, fill = black]}] (0.35,0) node[below right = -0.5mm, rotate = 75]{\tiny $|\Delta x|$} -- (0.35, 1);
\draw[black, thin] (0,3) circle (0.28*3) ;
\draw[black, thin,dotted] (1.3,0) -- (0,0) ;
\draw[black, thin,dotted] (1.3,3) -- (0,3) ;
\draw[black] (1.3,1.5) node[right]{\tiny $\frac{|\Delta x|}{|t_\cy - \tau_*|}$};
\draw[black, thick,{Latex[length=2mm, fill = black]}-{Latex[length=2mm, fill = black]}] (1.3,0) -- (1.3, 3) ;
\draw[black, thin, dashed] (0,0) -- (0.28*4, 0.96*4);
\draw[black, thin, dashed] (0,0) -- (-0.28*4, 0.96*4);
\draw[black, fill = black] (0,0) circle (1mm);
\draw[black] (2, -0.5) node[above]{\textbf{a.}};
\end{scope}
\begin{scope}[xshift = 4.5cm, rotate = 0, scale=0.8]           %%   Le soleil est un cône
\draw[white, ultra thin, fill = gray!8] (-1, 0.3) rectangle (5, 0.8);
\draw[white, ultra thin, fill = gray!8] (-1, 3.5) rectangle (5, 4);
\draw[black, thick] (2,2) circle (2cm);
\draw[black, thick] (-1, 0.3) -- (5, 0.3);
\draw[black, thick] (-1, 0.8) -- (5, 0.8);
\draw[black, thick] (-1, 4) -- (5, 4);
\draw[black, thick] (-1, 3.5) -- (5, 3.5);
\draw[black, thin, dashed] (2, 2) -- (2+3*0.6614, 2+3*0.75);
\draw[black, thin, dashed] (2, 2) -- (2-3*0.6614, 2+3*0.75);
\draw[black, thin] (2, 2.4) node[above]{$\theta_{\max}$};
\draw[black, thin, dashed] (2, 2) -- (2-3.3*0.8, 2-3.3*0.6);
\draw[black, thin, dashed] (2, 2) -- (2-2.8*0.527, 2-2.8*0.85);
\draw[black, thin, dashed] (2.2-2.8*0.65, 2.3-2.8*0.7) node[anchor = north east]{$\frac{\theta}{2}$};
\draw[black, thin, dashed] (2, 2) -- (2+3.3*0.8, 2-3.3*0.6);
\draw[black, thin, dashed] (2, 2) -- (2+2.8*0.527, 2-2.8*0.85);
\draw[black, thin, dashed] (1.8+2.8*0.65, 2.3-2.8*0.7) node[anchor = north west]{$\frac{\theta}{2}$};
\draw[black] (2, -1) node[above]{\textbf{b.}};
\end{scope}
%%                      Rayons croisés
\begin{scope}[xshift = 9.5cm, rotate = 0, scale=0.8]
\draw[red, thick] (3 - 0.28*1.849, 0.96*1.849) -- (3 - 0.28*1.849+ 0.5/0.96 - 0.28*0.5/0.96*0.28, 0.96*1.849+0.28*0.5/0.96*0.96); 
\draw[red, thick] (0.707*3.8594, 0.707*3.8594) -- (0.707*3.8594 - 0.5/0.96 + 0.28*0.5/0.96*0.28, 0.707*3.8594 - 0.28*0.5/0.96*0.96); 
\draw[red, pattern=horizontal lines, pattern color = red!80] (3 - 0.28*1.849, 0.96*1.849) -- (3 - 0.28*1.849+ 0.5/0.96 - 0.28*0.5/0.96*0.28, 0.96*1.849+0.28*0.5/0.96*0.96) -- (0.707*3.8594, 0.707*3.8594) -- (0.707*3.8594 - 0.5/0.96 + 0.28*0.5/0.96*0.28, 0.707*3.8594 - 0.28*0.5/0.96*0.96)  -- cycle;
\draw[black, thick] (0.25/0.707,0.25/0.707) -- (4, 4); %1 gauche
\draw[black, thick] (0.5/0.707,0) -- (4+0.25/0.707, 4-0.25/0.707); % 1 droite
\draw[black, ultra thick, -{Latex[length=3mm, fill = black]}] (0.5/0.707,0) -- (0.5/0.707 + 0.707, 0.707) node[below = 1.7mm]{$ u $}; % 1 droite vecteur
\draw[black, thick] (3,0) -- (3-0.28*4+0.28*0.5/0.96*0.28, 4*0.96-0.28*0.5/0.96*0.96); % 2 gauche
\draw[black, ultra thick, -{Latex[length=3mm, fill = black]}] (3,0) -- (3 - 0.28, 0.96) node[below = 4mm]{$ v $}; % 2 gauche vecteur
\draw[black, thick] (3 + 0.5/0.96 - 0.28*0.5/0.96*0.28, 0.28*0.5/0.96*0.96) -- (3+0.5/0.96-0.28*4, 4*0.96); % 2 droite
\draw[black, thick, {Latex[length=1.3mm, fill = black]}-{Latex[length=1.3mm, fill = black]}] (3-0.28*4+0.28*0.5/0.96*0.28-0.28*0.1, 4*0.96-0.28*0.5/0.96*0.96+0.96*0.1) -- (3+0.5/0.96-0.28*4-0.28*0.1, 4*0.96+0.96*0.1); % ecart 2
\draw[black, thick] (0.5*3-0.5*0.28*4+0.5*0.28*0.5/0.96*0.28-0.5*0.28*0.1 + 0.5*3+0.5*0.5/0.96-0.5*0.28*4-0.5*0.28*0.1, 0.5*4*0.96-0.5*0.28*0.5/0.96*0.96+0.5*0.96*0.1+ 0.5*4*0.96+0.5*0.96*0.1) node[anchor = south]{$\varepsilon$}; 
\draw[red, {Latex[length=1.5mm, fill = red]}-{Latex[length=1.5mm, fill = red]}] (3 - 0.28*1.849+ 0.5/0.96 - 0.28*0.5/0.96*0.28 + 0.096, 0.96*1.849+0.28*0.5/0.96*0.96+0.028) node[anchor = south west]{\small $ \varepsilon  |\sin(u, v)|^{-1}$} -- (0.707*3.8594+0.096, 0.707*3.8594+0.028)  ;
\draw[black] (2, -1) node[above]{\textbf{c.}};
\end{scope}
\end{tikzpicture}
\caption{Geometrical estimates} \label{fig:geom}
\end{figure}
Eventually, denoting $\ind_{ i \between j}$ the condition that $i$ and~$j$ create the first cycle within the tree~$\Upsilon$, we can bound the following integral using that $d \geq 3$,
\begin{align} \label{ineq:defcase}
\int \ind_{i\between j} \dd \sigma_{\dd} \dd \tau_* & \leq C_d \sum_{|\zeta | \leqslant t \mathbf{V}} \int \min\left[ 1, \left(\frac{ \varepsilon \mathbf{V}}{|(\tau_* - \tau_0^*) \Delta u |\cdot |u_i - u_k|}\right)^{d-1} \right] \dd \tau_* \\
& \leq C_d  (t\mathbf{V})^d \left(  \frac{\varepsilon \mathbf{V}}{|\Delta u|\cdot |u_i - u_k| } + (d-2) \left[\frac{\varepsilon \mathbf{V}}{|\Delta u|\cdot |u_i - u_k| }\right]^{(d-1) - (d-2)} \right) \nonumber \\
& \leq \tilde{C}_d \frac{t^d \bV^{d+1} \cdot \varepsilon }{|\Delta u|\cdot |u_i - u_k|} . \label{ineq:singdef}
\end{align}
Finally, if the condition~(\ref{ineq:deltax}) is not satisfied, then the condition $|\tau_* - \tau_0^*| \leq 2\varepsilon/|\Delta u|$  provides an even better bound integrating over $\tau_*$, with the same singularity $|\Delta u|^{-1}$.

The case $k = j$ is treated in a similar manner, without singularity since the disjunction on the height of the cone is not necessary: the cycle must occur because of a periodic shift $\zeta \neq 0_{\Z^d}$, necessarily providing a height larger than 1, which yields a negligible term (of order~$\varepsilon^{d-1}$, see for instance~\cite[proof of Proposition~B.2]{2023longcor}).

\subsection{Non-deflected case}

If none of both particles~$i$ and~$j$ has ever been deflected, we have considered the closest overlap, involving~$i$ and~$k$ at time~$\tau_\ov$. Taking the encounter time just before~$\tau_\ov$ as a reference, call it~$\tau_0$, we will use the same notation as previously, knowing that none of the velocities changes between~$\tau_0$ and~$\tau_\ov$ and that the velocities of~$i$ and~$j$ are given by their initial velocities~$v_i$ and~$v_j$ (see~Fig.~\ref{fig:arroz}, bottom right).
\noindent
Then, the overlap at time $\tau_\ov$ and the cycle creation at time $t_\cy$ give respectively the following conditions
\begin{equation} \label{eq:ovcond} \left\{
\begin{array}{l}
x_k(\tau_0) - x_i(\tau_0) + (\tau_\ov - \tau_0) (u_k - v_i) = \zeta_\ov + \varepsilon \omega_\ov \\
x_j(\tau_0) - x_i(\tau_0) + (t_\cy - \tau_0) (v_j - v_i) = \zeta_\cy + \varepsilon \omega_\cy,
\end{array}\right.
\end{equation}
for some~$(\zeta_\ov, \zeta_\cy, \omega_\ov, \omega_\cy) \in (\Z^d)^2 \times (\Sf^{d-1})^2$.
Since before~$\tau_\cy$ all encounters are clustering, removing the two encounters~$t_\cy$ and~$\tau_\ov$ divides the past encounters that involve~$i$ and~$j$ into two connected components: $C_1$ containing~$j$ and~$k$ and $C_2$ containing~$i$.
By construction, when $\vec{x}_{\ov}$ varies alone, both dynamics components~$C_1$ and~$C_2$ move rigidly one with respect to the other, and the distance $x_k(\tau_0) - x_j(\tau_0)$ remains fixed, so that
\begin{equation*}
\vec{x}_{\ov} =  x_k(\tau_0) - x_i(\tau_0) = [x_k(\tau_0) - x_j(\tau_0)] + [x_j(\tau_0) - x_i(\tau_0)] 
\end{equation*}
has to belong~(\ref{eq:ovcond}) to two cylinders of axes $(u_k - v_i)$ and $(v_j - v_i)$ and of common width $\varepsilon$.

\paragraph{If $j \neq k$,} the volume of the intersection of these cylinders is of order $\varepsilon^d |\sin(u_k - v_i, v_j - v_i)|^{-1}$ (cf. Fig.~\ref{fig:geom},~\textbf{c.}), so that with the same notation as in the previous deflected case~(\ref{ineq:defcase}), recalling the changes of variables~\eqref{eq:jacob} and~\eqref{lgn:change_var_sigom}, we have
\begin{align} \label{ineq:singove}
\int \ind_{i \between j} \dd \sigma_\ov \dd \tau_\ov = \frac{\mu}{|v_i - u_k| } \int \ind_{i \between j} \dd \vec{x}_\ov & \leq C_d  \frac{(t\bV)^{2d} \cdot \varepsilon}{|\sin(u_k - v_i, v_j - v_i)|\cdot |v_i - u_k| }.
\end{align}

\paragraph{If $j = k$,} since none of the particles has been deflected, this means that the cycle is due to the periodic conditions with $\zeta_\ov \neq \zeta_\cy$. Hence, the difference of both lines of the system~(\ref{eq:ovcond}) provides
\begin{equation*}
( \tau_\ov - t_\cy) (u_j - u_i) = (\zeta_\ov - \zeta_\cy) + \varepsilon (\omega_\ov - \omega_\cy),
\end{equation*}
so that for $\varepsilon$ small enough $(u_j - u_i)$ must belong to a cone of opening at most $3\varepsilon$, and integrating as before over $v_i$ and $v_j$ provides a negligible term (of order $\varepsilon^{d-1}$, see once again~\cite[proof of Proposition~B.2]{2023longcor} for instance).

\subsection{Handling the singularities}  \label{sec:singularities}

Now, we are left with some singularities to handle in our bounds~\eqref{ineq:singdef} and~\eqref{ineq:singove}. To do so, we will either: use some of the relative velocities appearing in the product over the edges to cancel them, use a previous deflection to integrate them, or integrate them using the initial velocities~$\uv_{r+k}$ in the case where there is no such deflection. First, the singularity $|u_i - u_k|^{-1}$ in~(\ref{ineq:singdef}) and~(\ref{ineq:singove}) appears in the product over the edges in the integral~(\ref{cycle_integral}) that we want to bound, so that it cancels out. The remaining singularities above hence consist in 
\begin{align*}
\frac{1}{|u_j - u_i|} + \frac{1}{|u_j - u_k|} + \frac{1}{|u_j - u_l|} + \frac{1}{ |\sin(u_k - v_i, v_j - v_i)|  }.
 \end{align*} 
For one of the three first singularities, of the form, $|u_j - u_n|^{-1}$, we will discriminate on the history of~$j$ and~$n$.

\paragraph{If~$j$ or~$n$} has been deflected by a particle~$m$ at a closest time~$\tau_{\tilde{\dd}}$, we will integrate over this collision's parameters. If this deflection is between~$j$ and~$n$, by the scattering identity $|u_j - u_n| = |u_j^\star - u_n^\star | $ between post and pre-collisional velocities, the singularity cancels out in the product over the edges. Otherwise, let us say by symmetry that~$j$ is colliding~$m \neq n$, so that in particular $u_n$ does not depend on~$\sigma_{\tilde{\dd}}$. Using the deflection equation~(\ref{eq:sigma}), we write
\begin{align*}
\int \frac{\dd \sigma_{\tilde{\dd}} \dd \tau_{\tilde{\dd}}}{|u_j - u_n|}  & = \int \dd \sigma_{\tilde{\dd}} \dd t_{\tilde{\dd}} \left| \frac{u_j^\star + u_m^\star}{2} - u_n + \frac{|u_j^\star - u_m^\star|}{2} \sigma_{\tilde{\dd}} \right|^{-1}\\
& = \frac{2}{|u_j^\star - u_m^\star|} \int \dd \sigma_{\tilde{\dd}} \dd \tau_{\tilde{\dd}}  \left| \frac{u_j^\star + u_m^\star - 2 u_n}{|u_j^\star - u_m^\star|}  + \sigma_{\tilde{\dd}} \right|^{-1}.
\end{align*}
We are brought back to studying an integral of the following form, computed in hyperspherical coordinates (cf.~Section~\ref{app:scattering}) for $w = (1,0,\dots, 0)$, its maximum value being reached for any $|w| = 1$, 
\begin{equation*}
\int \frac{\dd \sigma}{|w + \sigma|} = C_d \int_0^\pi \frac{\sin^{d-2} \theta }{\sqrt{1 - \cos \theta}}\dd \theta = \tilde{C}_d.
\end{equation*}
The remaining singularity $\lvert u_j^\star - u_m^\star\rvert^{-1}$, due to the collision, now appears in the product over the edges in~\eqref{cycle_integral} and hence get cancelled like~$\lvert u_i - u_k\rvert^{-1}$.

Nevertheless, to be able to integrate over these parameters~$(\sigma_{\tilde{d}}, \tau_{\tilde{d}})$, we first need to dispose of the velocities that are impacted by them: as we did for $K$ in the beginning of this appendix, we fix the corresponding set~$J$ of edges, and sum over them beforehand (along with the choice of~$K$, the choice of these disjoint sets corresponds to a factor~$3^\cy$).

\paragraph{If neither~$j$ nor~$n$} has ever been deflected before~$\tau_\dd$, we integrate directly the singularity over the velocity~$v_j$ at time~0, using part of the exponential decay~$e^{-\frac{\beta}{2}|v_j|^2}$, as the singularity is locally integrable in dimension~$d > 1$.

\paragraph{The sine singularity} only appears in the non-deflected case, and so is integrated in the same way over $v_i$ and $v_j$, using for example once again hyperspherical coordinates
\begin{align*}
\int \frac{e^{-{v_j}^2 - {v_i}^2}}{|\sin(u_k - v_i, v_j - v_i)|} \dd v_j \dd v_i & = \int \left(\int \frac{e^{-(v + v_i)^2 - {v_i}^2}}{|\sin(u_k - v_i, v)|} \dd v \right)\dd v_i\\
& \leq \int \left(C_d \int r^{d-1} e^{-r^2 + 2 |v_i| r }\frac{ \sin^{d-2} \theta}{|\sin \theta |}\dd \theta \dd r \right) e^{-{v_i}^2} \dd v_i \leq \tilde{C}_d.
\end{align*}

\subsection{Conclusion}
In the end, we have obtained enough smallness and there is no more singularity, so that we can eventually sum over all the remaining edges, and conclude the domination exactly as in Section~\ref{sec:bounds}, yielding the same bound as~\eqref{lgn:dom_SR}. The choice of the disjoint sets $K, J \subset E_\Upsilon$ gives a factor $3^c$, which becomes $3^{r+k}/2$ with the sum over $c$, so that we get the following final bound (the constant~$C$ stemming from the different cases)
\begin{align} 
& \sum_{\Upsilon \in \mT_{R, k}^{\dd}} \int \ind_{T^\dd_{R} = \Upsilon} \dd \uv_{r+k}   \dd \usig_{E_\Upsilon} \dd \utau_{E_\Upsilon}  \prod_{e \in E_\Upsilon} s_{e} | v_{m}({\tau_e}^+ )  - v_{m'}({\tau_e}^+)| e^{- \beta \|\uv_{r+k}\|^2}\times \ind^\between_{R} \label{lgn:final_est_cycle} \\
& \hspace{1cm} \leq C (r+k)^7 3^{r+k} \frac{t^{r+k-1}}{(r+k-1)! } \int  \dd \uv_{r+k} e^{-\beta \bV^2} (r+k)^{r-1} \left(\sqrt{(r+k)\bV^2}\right)^{r+k-1} \cdot \bV^{2d+1} \varepsilon \nonumber \\[1em]
& \hspace{1cm} \leq \varepsilon\cdot (\tilde{C} t)^{r+k-1} (r+k)^{r-1}  . \nonumber
\end{align} 
This bound is similar to~\eqref{lgn:dom_SR}, possibly for different constants, with an additional factor~$\varepsilon$. Hence, the same computation as in Section~\ref{sec:bounds} leads to Proposition~\ref{prop:discard_cycles}. 
Note that more precisely, $\tau_\dd$ or $\tau_\ov$ has been integrated before, and $v_i$ and $v_j$ also might have been integrated before with the singularities, yet it only changes constants.

To sum up the strategy, we fixed a configuration of cycle $[i,j,k,l, (\tau_c, \tau_\dd, \tau_\ov), J, K]$, summed over all the edges in~$K$ whose associated velocities would have been impacted by the collision angle~$\sigma_\dd$, so that we could integrate over the collision parameters~$\sigma_\dd$, and~$\tau_\dd$ or~$\tau_\ov$, or~$\vec{x}_\ov$. This made appear some singularities that we have handled by summing first over the edges of~$J$, to free some collision parameters or velocities, eventually integrating over them.
\CQFD

\section{Proof of the fluctuation and large deviation results} \label{sec:fluct_LD_proof}

We eventually prove in this section Theorems~\ref{theo:fluct} and~\ref{theo:LD}. We  start presenting tightness properties, so that we will be able to prove these theorems in a weak form, to extend them afterwards to stronger topologies.

\subsection{Tightness}  \label{sec:tightness}

We expose here tightness results for the empirical measure and the fluctuation field, useful to extend weak results, in the sense of observables, to results in the strong Skorokhod topology on the set of trajectories~$\mathrm{Traj}([0,t], \mM(\mD))$. 

The results of this section are given without complete proofs, which can be found in the paper dealing with the fluctuations and large deviations  of the general symmetric hard-sphere dynamics~\cite{2023grandev}. Indeed, these proofs rely on the bounds on the cumulants, that we proved in Section~\ref{sec:bounds} in our mixed model, and are otherwise identical. 
\begin{prop} \label{prop:tight_LD} There exists a distance~$d$, based on normalized bounded observables, and associated to the strong Skorokhod topology on $\mathrm{Traj}([0,t], \mM(\mD))$, such that
\begin{equation} \label{eq:tight_LD}
\lim_{A \to +\infty} \lim_{\lambda \to + \infty } \frac{1}{\lambda} \log \Pro \left[ \sup_{s \in [0,t]} d(0, \tilde{\pi}^\varepsilon_s) \geq A \right] = - \infty, 
\end{equation}
and for any $\eta > 0$,
\begin{equation} \label{eq:tight_LD2}
\lim_{\delta \to +\infty} \lim_{\lambda \to + \infty } \frac{1}{\lambda} \log \Pro \left[ \sup_{|s-s'| < \delta} d( \tilde{\pi}^\varepsilon_s , \tilde{\pi}^\varepsilon_{s'} ) > \eta \right] = - \infty.  
\end{equation}
\end{prop}
The proof of this proposition follows from the estimate on the cumulants given in Proposition~\ref{prop:cum_bound}, with the same proof as in~\cite[Proposition~7.3.2]{2023grandev}.

\begin{prop}[Tightness of the fluctuation field] \label{prop:fluct_tight}
There exists a distance~$\tilde{d}$, based on some normalized bounded observables, and associated to the strong Skorokhod topology on $\mathrm{Traj}([0,t], \mM(\mD))$, such that
\begin{equation*} 
\lim_{A \to +\infty} \lim_{\mu \to + \infty } \Pro \left[ \sup_{s \in [0,t]} \tilde{d}(0,\zeta_s^\varepsilon) \geq A \right] = 0 
\end{equation*}
and for any $\eta > 0$,
\begin{equation} \label{eq:tightness_fluctuation_field}
\lim_{\delta \to 0} \lim_{\mu \to + \infty }\Pro \left[ \sup_{|s-s'| < \delta} \tilde{d}( \zeta^\varepsilon_s , \zeta^\varepsilon_{s'}) > \eta \right] = 0.  
\end{equation}
\end{prop}
The proof of this proposition is once again identical to the one in~\cite[Proposition~6.2.3]{2023grandev}, using the bounds proved in Proposition~\ref{prop:cum_bound}. 
Note that in the integrated form, dominating roughly the bounded observables, it is easy to remark that the bounding estimates that are true in the non-linear symmetric case remain true in our linear tagged model, since by symmetry one can rewrite the sum over the tags~$\ul_p$ as
\begin{align*} 
\sum_{\ul_p \in \Lambda_p } \lambda^{|\ul_{p}|} \mu^{p - |\ul_{p}|} \int M_\beta^{\otimes p} \varphi_0^{\otimes \ul_{p}} \ind_{\mX_{p}^\varepsilon} \dd \uz_{p} = \mu^p \int M_\beta^{\otimes p} \left( \frac{\lambda}{\mu} \varphi_0 + 1\right)^{\otimes p} \ind_{\mX_{p}^\varepsilon} \dd \uz_{p},
\end{align*}
which formally corresponds to a bounded symmetric initial data~$\Bigl[M_\beta ( p_\mu \varphi_0 + 1 )\Bigr]$.

\subsection{Convergence of the fluctuation field} \label{sec:fluct_proof}

Since the fluctuation field~$(\zeta_t^\varepsilon)$ defined in~\eqref{def:fluctuation_field} is tight by Proposition~\ref{prop:fluct_tight}, identifying its limit moments is enough to characterize its limit, as they decrease fast enough (see the method of moments in~\cite[Theorem~30.1]{bill1979}). 
Hence, to identify the limit fluctuation field with a Gaussian process, and to find its covariance, we consider the following sampled observable 
\begin{equation*}
H(z^{[0,t]}, \ell) = \ind_{\ell = 1} \sum_{j=1}^J \psi_j(z^{[\theta_j]}).
\end{equation*} 
For simplicity, we denote~$\tilde{f}_n(t) \doteq f_n(t, \ul_n = \underline{1}_n)$ and $\tilde{F}_1(t) \doteq F_1(t, \ell_1 = 1)$ the cumulants associated to the tagged particles only.
Recall that the random set~$\mS_\lambda$ contains the labels of tagged particles. With the same notation as in Sections~\ref{sec:stat_descr} and~\ref{sec:stat_refinements}, we have
\begin{equation*}
\zeta^\varepsilon_t[H] = \frac{1}{\sqrt{\lambda}} \sum_{j=1}^J \left[ \sum_{i \in \mS_\lambda} \psi_j\left(Z_{\varepsilon,i}^{[\theta_j]}\right) - \lambda \int \tilde{F}_1^\varepsilon(\theta_j) \psi_j  \right].
\end{equation*}
To characterize the law of this fluctuation field, we look at its Fourier transform, using the generalized cumulant generating function~\eqref{def:cumulant_generating_function3} to write
\begin{equation} \label{eq:fourier_fluct}
\E\left[ e^{i \zeta^\varepsilon_t[H]} \right] =  \exp\left(\mathfrak{G}_\varepsilon^{[0,t]}\left[\frac{iH}{\sqrt{\lambda}}\right]\right)  \exp\left(-i\sqrt{\lambda} \sum_{j=1}^J \int \tilde{F}_1^\varepsilon(\theta_j) \psi_j  \right).
\end{equation}
Expanding the cumulant generating function~\eqref{def:cumulant_generating_function3} yields

\begin{align*}
\mathfrak{G}_\varepsilon^{[0,t]}\left[\frac{i H}{\sqrt{\lambda}}\right] & = \sum_{p \geqslant 1} \frac{\lambda^p}{p!} \int \tilde{f}_p^\varepsilon\left[e^{\frac{i}{\sqrt{\lambda}}H} - 1\right] \\
& = \sum_{p \geqslant 1} \frac{\lambda^p}{p!} \int \tilde{f}_p^\varepsilon\left[ \frac{i}{\sqrt{\lambda}}H -\frac{1}{2\lambda} H^{ 2 } + O\left( \frac{\|H\|}{\sqrt{\lambda}} \right)^3 \right].
\end{align*}
Using the estimates on the cumulants stated in Proposition~\ref{prop:cum_bound}, and the convexity of $x \mapsto (1+x)^p$, we bound for any $p \geqslant 2$
\begin{align*}
\frac{\lambda^p}{p!} \left\lvert \int \tilde{f}_p^\varepsilon\left[ \frac{i}{\sqrt{\lambda}} H + O\left( \frac{\|H\|}{\sqrt{\lambda}} \right)^2  \right] \right\rvert \leq \frac{\lambda^p C^p}{\mu^{p-1}}\cdot \left[ \frac{\|H \|^p}{\lambda^{\frac{p}{2}}} + 2^p O\left(\frac{\|H \|}{\sqrt{\lambda}}\right)^{p+1} \right],
\end{align*}
so that 
\begin{equation*}
\sum_{p \geqslant 2} \frac{\lambda^p}{p!} \left\lvert \int \tilde{f}_p^\varepsilon\left[e^{\frac{i}{\sqrt{\lambda}}H} - 1 \right] \right\rvert \leq 2 (C\|H\|)^2 \frac{\lambda}{\mu}.
\end{equation*}
Then the cumulant generating function writes
\begin{align*}
\mathfrak{G}_\varepsilon^{[0,t]}\left[\frac{iH}{\sqrt{\lambda}}\right] & = \int \tilde{F}_1^\varepsilon \left[ i\sqrt{\lambda} H - \frac{H^{2}}{2} + O\left( \frac{\|H\|}{\sqrt{\lambda}} \right)^3 \right] + 2 (C\|H\|)^2 \frac{\lambda}{\mu} 
\end{align*}
where the first term simplifies in the Fourier transform~\eqref{eq:fourier_fluct}. Eventually, in the mixed scaling~\eqref{eq:scaling} this provides the following convergence
\begin{equation*}
\E\left[ e^{i \zeta^\varepsilon_t[H]} \right] \xrightarrow[\varepsilon \to 0]{} \exp\left( - \frac{1}{2} \sum_{i,j=1}^J \int \tilde{F}_1\left[\psi_i(z^{[\theta_i]}) \psi_j(z^{[\theta_j]})  \right] \right)
\end{equation*}
thanks to the convergence of the integrated cumulants shown in Theorem~\ref{theo:cum_conv}. By the method of moments mentioned above, and thanks to the tightness~\eqref{eq:tightness_fluctuation_field} of the fluctuation field,  $\zeta^\varepsilon_t$ converges to a Gaussian process~$\zeta_t$, with time-observable covariance
\begin{equation*}
\E\Bigl[ \zeta_s[g] \zeta_t[h] \Bigr] = \int \tilde{F}_1\left[g(z^{[s]})h(z^{[t]})\right].
\end{equation*}
At each time, it is hence characterized by its covariance~\eqref{eq:cov_fluct} with respect to observables. This concludes the proof of Theorem~\ref{theo:fluct}, recognizing the expansion~\eqref{eq:F1_lim} of the solution $M_\beta\varphi$ to the Rayleigh--Boltzmann equation. \CQFD

\subsection{Large deviations of the empirical measure} \label{sec:LD_proof}

We explain here how the convergence of the cumulant generating function (Sections~\ref{sec:cumulant_generating_function_conv} and~\ref{sec:HJ}) leads to the large deviation principle for the empirical measure (Theorem~\ref{theo:LD}), following the method of~\cite{2023grandev}.

\paragraph{Upper bound}

Let~$\bF$ be a closed set for the Skorokhod topology. In particular, $\bF$ is also closed for the weaker topology given by the opens of the form
\begin{equation*}
 \bO_{h, \delta}(\mathfrak{v}) \doteq \Bigl\{ \fm \in \mathrm{Traj}([0,t], \mM(\mD))  \ , \ \left\lvert \{h, \fm - \mathfrak{v} \} \right\rvert < \delta \Bigr \}, 
\end{equation*}
for observables $h \in \mC^\infty_c$ and measures $\mathfrak{v} \in \mathrm{Traj}([0,t], \mM(\mD))$, where $\delta > 0$ is fixed (we recall the definition~\eqref{def:crochet} of the transport filtered mean).
For any open set of this form, one can write
\begin{align*}
\Pro \bigl(\tilde{\pi}^\varepsilon_t \in \bO_{h, \delta}(\mathfrak{v})\bigr) 
& \leq \E\left[ \exp\Bigl( \lambda \{ h, \tilde{\pi}^\varepsilon_t \} - \lambda \{ h, \mathfrak{v} \} + \lambda \delta \Bigr) \right]. 
\end{align*}
Now, since by definition
\begin{align*}
\{ h, \tilde{\pi}^\varepsilon_t \} = \frac{1}{\lambda} \sum_{i \in \mS_\lambda} \left[ h\left(t, Z_i^{[t]}\right) + \int_0^t (\partial_s + v\cdot \nabla_x)h\left(s, Z_i^{[s]}\right) \right],
\end{align*}
and harnessing the definition~\eqref{def:cumulant_generating_function3} of the cumulant generating function, we get
\begin{align*}
\Pro \bigl(\tilde{\pi}^\varepsilon_t \in \bO_{h, \delta}(\mathfrak{v})\bigr) 
& \leq \exp\left( \mathfrak{G}_\varepsilon^{[0,t]}\left[ h_t + \int (\partial_s + v\cdot \nabla_x) h_s \right] - \lambda \{ h, \mathfrak{v} \} + \lambda \delta \right).
\end{align*}
Taking the limit superior, thanks to the convergence~\eqref{eq:conv_G} of the generating function, and using the notation~$\mI(t,g)$ for the limit cumulant generating function~\eqref{def:Itg}, one has
\begin{align*}
\limsup_{\varepsilon \to 0} \frac{1}{\lambda} \log \Pro \bigl(\tilde{\pi}^\varepsilon_t \in \bO_{h, \delta}(\mathfrak{v})\bigr) 
& \leq \mI(t,h)-1 - \{ h, \mathfrak{v} \} + \delta.
\end{align*}
By definition of the Legendre transform~\eqref{def:Legendre}, 
for our fixed $\delta > 0$ there exists an observable~$g \in \mathbb{B}_{t, \beta}$ such that 
\begin{equation*}
\{ g, \mathfrak{v} \} - \mI(t, g) + 1> \mathbf{\Lambda}(t, \mathfrak{v}) - \delta,
\end{equation*}
so that
\begin{align*}
\limsup_{\varepsilon \to 0} \frac{1}{\lambda} \log \Pro \bigl(\tilde{\pi}^\varepsilon_t \in \bO_{g, \delta}(\mathfrak{v})\bigr) 
& \leq -\mathbf{\Lambda}(t, \mathfrak{v}) + 2 \delta.
\end{align*}
 By the tightness result stated in Proposition~\ref{prop:tight_LD},  
$\bF$ is in fact compact, and one can extract a finite covering of open sets~$\bF \subset \cup_{i=1}^k \bO_{g_i, \delta}(\mathfrak{v}_i)$, so that 
\begin{align*}
\limsup_{\varepsilon \to 0} \frac{1}{\lambda} \log \Pro \bigl(\tilde{\pi}^\varepsilon_t \in \bF\bigr) \leq - \inf_{i\leqslant k } \mathbf{\Lambda}(t, \mathfrak{v}_i) + 2 \delta \leq - \inf_{\mathfrak{v} \in \bF}  \mathbf{\Lambda}(t, \mathfrak{v}) + 2 \delta.
\end{align*}
The result~\eqref{ineq:upperbound} follows, considering that~$\delta  > 0$ may be chosen arbitrarily small. 

\paragraph{Lower bound}
The lower bound follows from more elaborate methods, that works only for measures~$\mathfrak{v}$ such that the Legendre transform~$\mathbf{\Lambda}(t, \mathfrak{v})$ is reached for an observable~$h \in \mathbb{B}_{t, \beta}$, which is the case when~$\mathfrak{v}$ is a strong solution of the biased Boltzmann--Hamilton--Jacobi equation~\eqref{eq:HJsystem}
\begin{equation*} 
(\partial_s - v\cdot \nabla_x) \mathfrak{v} = \int \dd v_\cy \dd \omega \langle v  - v_\cy,  \omega \rangle_+   M_\beta(v_\cy) \left(\mathfrak{v}(v')  e^{p(z) -  p(z')} - \mathfrak{v}(v) e^{p(z') - p(z)} \right),
\end{equation*}
for some~$p \in\mathbb{B}_{t, \beta}$, whence the restriction on the infimum defining the lower bound. The algebraic details are exaclty the same as in~\cite[Chapter~7]{2023grandev} and proceed of industrious topological considerations, so that we do not replicate them here.  \CQFD

\section*{\Huge Appendix}
\setcounter{section}{0}
\renewcommand{\thesection}{\Alph{section}}
\renewcommand{\theHsection}{appendixsection.\Alph{section}}

\section{Scattering parametrization} \label{app:scattering}

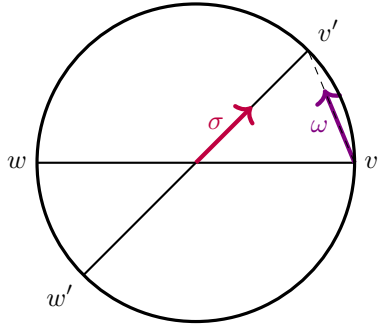
\begin{figure}[h!] %  ht: here OR top 
\centering 
\begin{tikzpicture}
\begin{scope}[scale = 0.7]
\draw[black,  thick, -] ( -3, 0) node[left]{$w$} -- (3,0) node[right]{$v$};
\draw[black,  thick, -] ( -1.414*3/2,-1.414*3/2) node[anchor = north east]{$w'$} -- ( 1.414*3/2,1.414*3/2) node[anchor = south west]{${v}'$};
\draw[purple, ultra thick, ->] ( 0, 0) -- ( 1.414*3/4,1.414*3/4) ;
\draw[violet, ultra thick, ->] ( 3, 0) -- ( 3 - 0.383*3/2,0.924*3/2);
\draw[purple, ultra thick] ( 1.414*3/8*1.4,1.414*3/8*1.4) node[left]{$\sigma$};
\draw[violet, ultra thick]  ( 3 - 0.383*3/4,0.924*3/4) node[left]{$\omega$};
\draw[black, thin, dashed]  ( 3 ,0) -- ( 1.414*3/2,1.414*3/2);
\draw[black, very thick,] (0,0) circle (3 cm);
\end{scope}
\end{tikzpicture}
\caption{Deflection parameters}
\label{fig:omesig}
\end{figure}

This appendix is dedicated to the parametrizations of the scattering. Indeed, the scattering system~\eqref{eq:coll}, denoting $(v, w) \doteq (v_i, v_j)$, is equivalent to 
\begin{equation} \label{eq:sigma}
\left\{  \def\arraystretch{2}  \begin{array}{l}
v' = \frac{v+w}{2} + \sigma \frac{|v-w|}{2} \\
w' = \frac{v+w}{2} - \sigma \frac{|v-w|}{2},
\end{array} \right.
\end{equation}
with~$\sigma = \frac{{v_i}' - {v_j}'}{|{v_i}' - {v_j}'|}$ denoting the deviation's direction with respect to the mean velocity. Since $(v,w,v',w')$ are in the same plane, $\sigma$ and $\omega$ also belong to this same plane (see Fig.~\ref{fig:omesig}). Thus, we might reduce their dependency to the 1-dimensional circle, where the following holds
\begin{equation*}
\sigma = \pi - 2 \omega,
\end{equation*}
$\sigma$ going through the circle, and $\omega$ through the half-circle. 

\paragraph{Change of variable $\omega \mapsto \sigma$.}

Let us denote $u = w-v$ the precollisional relative velocity. We want to change the variable~$\omega \in \Sf^{d-1}_+$ into~$\sigma \in \Sf^{d-1}$ in the following integral, written in hyperspherical coordinates
\begin{align*}
\int_{\Sf^{d-1}_+} f(\omega)\langle \omega, u \rangle \dd \omega & = \int f(\omega) |u| \cos(\theta_{d-2}) \sin^{d-2} \theta_1 \sin^{d-3} \theta_2 \dots \sin^2 \theta_{d-3} \sin \theta_{d-2}\ \dd \underline{\theta}_{d-1},
\end{align*}
where the $d-2$~first angles~$(\theta_i)_{i \leqslant d-2}$ range~$[0,\pi]$, and $\theta_{d-1}$ ranges~$[0, 2 \pi]$. Thus, choosing the angle $\theta \doteq \widehat{(u, \omega)}$ as $\theta_{d-2}$, the half-sphere that $\omega$ ranges corresponds to $\theta \in [0, \pi/2]$. In Fig.~\ref{fig:omesig}, one can see that $\widehat{(u, \sigma)} = \pi - 2 \widehat{(u, \omega)}$, so that we will compute the natural change of variable $\theta = (\pi - \phi)/2$. Denoting $\check{\theta} = (\theta_1, \dots, \theta_{d-3}, \theta_{d-1})$, we get by trigonometry formulas that
\begin{align}
\int_{\Sf^{d-1}_+} f(\omega)\langle \omega, u \rangle \dd \omega & = \int_0^\pi \int f(\omega) \frac{|u|}{2} \sin\left(\frac{\phi}{2}\right) \cos\left(\frac{\phi}{2}\right) \sin^{d-2} \theta_1 \sin^{d-3} \theta_2 \dots \sin^2 \theta_{d-3}\   \dd \check{\theta} \dd \phi\nonumber \\
& = \int_0^\pi \int f(\omega) \frac{|u|}{4} \sin(\phi) \sin^{d-2} \theta_1 \sin^{d-3} \theta_2 \dots \sin^2 \theta_{d-3}\  \dd \check{\theta} \dd \phi \nonumber\\
& = \frac{|u|}{4} \int_{\Sf^{d-1}} f(\omega)  \dd \sigma. \nonumber%\label{lgn:change_var_sigom}
\end{align} 
It is useful in the computations of this paper to be able to switch from a variable to another, with the following formula.
\begin{lemm}[Change of scattering angle] For any measurable function~$f$, one has
\begin{align}
\int_{\Sf^{d-1}_+} f(\omega)\langle \omega, u \rangle \dd \omega = \frac{|u|}{4} \int_{\Sf^{d-1}} f(\omega)  \dd \sigma. \label{lgn:change_var_sigom}
\end{align} 
\end{lemm}

\section{Linear Boltzmann equations} \label{app:RB} 

In this section, we prove the well-posedness of the system of linear Boltzmann-like equations~\eqref{eq:modBoltz} that corresponds to the Hamilton--Jacobi system presented in Section~\ref{sec:HJ}. Before stating this result in Proposition~\ref{prop:BHJ}, we start with a general study of the linear Rayleigh--Boltzmann equation.
 
For an inverse temperature~$\beta > 0$ and $k\in \N^*$, we define the weighted functional space $\mF_{k, \beta}$ of measurable functions defined almost everywhere on the domain $\mD^k$ such that
\begin{equation} \label{def:foncspace}
\|f_k\|_{k, \beta} \doteq \supess_{\uz_k \in \mD^k} \bigl\lvert f_k(\uz_k) \exp( \beta \|\uv_k\|^2) \bigr \rvert < \infty,
\end{equation}
hence decreasing at least as the Gaussian equilibrium~$M_{2\beta}^{\otimes k}$ in velocities.

\subsection{Linear Rayleigh--Boltzmann equation}

We proceed to a similar decomposition of this equation as in~\cite{1996Glassey}, and compute the integrals in hyperspherical coordinates. Through the change of unknown~$R \doteq M_\beta^\frac{1}{2} \varphi$, the linear Rayleigh--Boltzmann equation~\eqref{eq:phi} might indeed be written
\begin{equation} \label{eq:R}
\partial_t R + v \cdot \nabla_x R = \hat{\mK} R(v) - \nu_\beta(v) R(v),
\end{equation}
where the gain operator is
\begin{equation} \label{def:Kgain2}
\hat{\mK} R(v) \doteq \int \langle v_c - v, \omega \rangle_+ M_\beta(v_c) M_\beta^{\frac{1}{2}}(v) M_\beta^{-\frac{1}{2}}(v') R(v') \dd \omega \dd v_c.
\end{equation}
and the loss factor is
\begin{equation} \label{def:nubeta}
\nu_\beta(v) \doteq \int \langle v_c - v, \omega \rangle_+ M_\beta(v_c) \dd \omega \dd v_c.
\end{equation}
The loss factor being bounded below in~$v$, it will provide some decay of the solution, while the gain operator will rather make it grow. We will prove in the following lemma that, in this choice of unknown, it is nonetheless bounded.
\begin{lemm}[Integral kernel of the gain part] \label{lemm:opkern2} The gain operator~$\hat{\mK}$ defined above~\eqref{def:Kgain2} has the following kernel integral formula
\begin{equation*} 
\hat{\mK} R(v) = \sqrt{ \frac{\beta}{2\pi}} \int \frac{R(\eta)}{|\eta - v|^{d-2}}  \exp\left(-\beta \left[\frac{|\eta - v|^2}{8} + \frac{\bigl( |\eta|^2- |v|^2\bigr)^2 }{8 |\eta - v|^2 } \right]\right) \dd \eta,
\end{equation*}
and is bounded on~$\Lp^\infty(\mD)$ as a consequence of the following estimate on its integral kernel
\begin{equation*} 
\sqrt{ \frac{\beta}{2\pi}} \int \frac{R(\eta)}{|\eta - v|^{d-2}}  \exp\left(-\beta \left[\frac{|\eta - v|^2}{8} + \frac{\bigl( |\eta|^2- |v|^2\bigr)^2 }{8 |\eta - v|^2 } \right]\right) \dd \eta \leq \frac{17|\Sf^{d-1}|}{\beta |v|} \cdot
\end{equation*}
\end{lemm}
\begin{prof}
The proof of the kernel formula can be found in detail in~\cite{1996Glassey}, or in a more concise way in our PhD thesis~\cite[Lemma~3.3.2]{fouthese}. We denote~$u \doteq v - \eta$, so that
\begin{equation*} 
|\eta|^2 - |v|^2 = |v-u|^2 - |v|^2 = |u|^2 - 2 u \cdot v.
\end{equation*}
Using the hyperspherical coordinates (see Appendix~\ref{app:scattering}) with respect to~$v$, denoting~$\dd J(\underline{\theta})$ the Jacobian corresponding to the measure on the unit sphere, with~$\theta = \theta_{d-2} \doteq \widehat{(u, v)}$, we thus have
\begin{align*}
\int \frac{\dd \eta}{|\eta - v|^{d-2}} e^{-\beta \left[\frac{|\eta - v|^2}{8} + \frac{\bigl( |\eta|^2- |v|^2\bigr)^2 }{8 |\eta - v|^2 } \right]} & = \int \dd r \dd J(\underline{\theta}_{d-1}) r e^{- \frac{\beta}{8} \left[r^2 + ( r - 2|v| \cos \theta)^2  \right]}.
\end{align*}
Now, one can integrate
\begin{align*}
\int_0^\infty r e^{- \frac{\beta}{8} \left[r^2 + ( r - 2|v| \cos \theta)^2 \right]} \dd r &= \left[ \frac{-2}{\beta} e^{- \frac{\beta}{8} \left[r^2 + ( r - 2|v| \cos \theta)^2 \right]} \right]_0^\infty   + |v|\cos\theta \int_0^\infty r e^{- \frac{\beta}{8} \left[r^2 + ( r - 2|v| \cos \theta)^2 \right]} \dd r.
\end{align*}
In the second integral, part of the weight concentrates around~0, and the rest around~$2|v| \cos \theta$, so it can be split into
\begin{align*}
\int_0^{|v|\cos \theta} e^{- \frac{\beta}{8} \left[r^2 + ( r - 2|v| \cos \theta)^2 \right]} \dd r & \leq \int_{|v|\cos \theta}^\infty e^{- \frac{\beta}{8} u^2} \dd u \\
& \leq \frac{8}{\beta |v| \cos \theta} e^{- \frac{\beta}{8} |v|^2 \cos^2 \theta}
\end{align*}
(bounding by standard estimates on the Gaussian cumulative distribution function), and similarly
\begin{align*}
\int_{|v|\cos \theta}^\infty e^{- \frac{\beta}{8} \left[r^2 + ( r - 2|v| \cos \theta)^2 \right]} \dd r & \leq  \frac{8}{\beta |v| \cos \theta} e^{- \frac{\beta}{8} |v|^2 \cos^2 \theta}.
\end{align*}
In the end, we get 
\begin{align*}
\int_0^\infty r e^{- \frac{\beta}{8} \left[r^2 + ( r - 2|v| \cos \theta)^2 \right]} \dd r & \leq  \frac{2}{\beta} e^{- \frac{\beta}{8} \left[ 2|v| \cos \theta \right]^2} + \frac{16}{\beta } e^{- \frac{\beta}{8} |v|^2 \cos^2 \theta}.
\end{align*}
It remains to integrate over~$\theta$ as 
\begin{align*}
\int_0^\pi \sin \theta e^{- \frac{\beta}{8} |v|^2 \cos^2 \theta } \dd \theta & = \int_{-1}^1 e^{- \frac{\beta}{8} |v|^2 x^2 } \dd x \leq \frac{1}{|v|} \sqrt{\frac{8\pi}{\beta}} \cdot 
\end{align*}
The contribution of the integral over~$\theta$ to the volume of the sphere is~$2$, so that the integral over the other angles yields~$|\Sf^{d-1}|/2$. In the end, 
\begin{align*}
 \int_0^\infty r e^{- \frac{\beta}{8} \left[r^2 + ( r - 2|v| \cos \theta)^2 \right]} \dd r & \leq \frac{|\Sf^{d-1}|}{2} \left( \frac{2}{\beta}\cdot \frac{1}{|v|} \sqrt{\frac{2\pi}{\beta}} + \frac{16}{\beta}\cdot \frac{1}{|v|} \sqrt{\frac{8\pi}{\beta}} \right),
\end{align*}
which concludes the proof when multiplying by the factor~$\sqrt{\beta (2\pi)^{-1}}$. The kernel is also bounded for small values of~$|v|$ in dimensions~$d \geqslant 2$, so that the operator is bounded from~$\Lp^\infty(\mD)$ to itself.
\CQFD
\end{prof} 

\subsection{Boltzmann--Hamilton--Jacobi system}

We study here the Boltzmann--Hamilton--Jacobi system~\eqref{eq:modBoltz} appearing in the description of the large deviations of the empirical measure, for long times. In the quadratic case, it has been studied by Chenjiayue Qi for long times in the case of small initial data~\cite{Qi2024}. Our situation here is completely linear.

To gain boundedness of the gain operator appearing in the Rayleigh--Boltzmann equation, we computed a change of unknown that distributed differently the exponential weights. However, doing so, we caused the new gain operator to be greater than the loss factor for small velocities, so that the latter do not compensate totally the gain operator. Because of this, they will not immediately provide a maximum principle. 

In the specific case of the system that we study here, we have no hope of a maximum principle, because of the gain terms~$-\theta_t \chi_t$ and~$\theta_t\eta_t$ that depend on the observable~$\theta_t$. We recall the system
\begin{equation} \label{eq:modBoltz2}
\left\{ \begin{array}{lll}
(\partial_s - v\cdot \nabla_x) \chi = - \theta \chi  + \int \dd v_2 \dd \omega \langle v  - v_{2},  \omega \rangle_+  \Bigl(M_\beta(v_2')\chi(z') - M_\beta(v_2)\chi(z)\Bigr) \\[1em]
(\partial_s - v\cdot \nabla_x)\eta = + \theta \eta  - \int \dd v_2 \dd \omega \langle v  - v_{2}, \omega \rangle_+ M_\beta(v_2) \Bigl(\eta(z') - \eta(z)\Bigr)
\end{array} \right.
\end{equation}
of unknowns $(\chi, \eta)$, with the boundary conditions
\begin{equation*} 
\left\{ \begin{array}{l} \chi(0) = M_\beta \varphi_0 \\ \eta(t) = \gamma(t) \in \mF_{1, -\beta/4}.
\end{array} \right.
\end{equation*}
Recall the definition~\eqref{def:foncspace} of the functional space~$\mF_{1, \frac{\beta}{2}}$, embedded with the $\Lp^\infty(\mD)$-norm with a weight~$M_\beta^{-1}$. It thus corresponds to functions decaying faster than the~$\beta-$Gaussian in velocities.
\begin{prop}[Well-posedness of the Boltzmann--Hamilton--Jacobi system] \label{prop:BHJ} For any time~$t > 0$, for a bounded gain weight~$\theta \in \Lp^\infty(\mD)$, and boundary conditions
\begin{equation*} 
\left\{ \begin{array}{l} \chi(0) = M_\beta \varphi_0 \in \mF_{1, \beta/4} \\ 
\eta(t) = \gamma(t) \in \mF_{1, -\beta/4},
\end{array} \right.
\end{equation*}
the system~\eqref{eq:modBoltz2} has a unique global positive solution~$(\chi, \eta) \in \Lp^\infty\left([0,t], \ \mF_{1, \beta/4} \times \mF_{1, -\beta/4} \right)$.
\end{prop}

\begin{prof} We will compute a different change of unknown to go from the operators of the Boltzmann--Hamilton--Jacobi system~\eqref{eq:modBoltz2} to the bounded operator~
$\hat{\mK}$ defined in~\eqref{def:Kgain2}. Considering~$\chi(0) \in \mF_{1, \beta/4}$, we will look at $R(0) \doteq M_\beta^{-\frac{1}{2}} \chi(0) \in \Lp^\infty(\mD)$. The existence of a global solution for the first equation of the system~\eqref{eq:modBoltz2} above, on~$\chi$ in~$\mF_{1, \beta/4}$, is thus equivalent to the global existence in~$\Lp^\infty(\mD)$ of the modified Rayleigh--Boltzmann equation
\begin{equation*}
\partial_t R + v \cdot \nabla_x R = \hat{\mK} R(v) - \theta R  - \nu_\beta(v) R(v),
\end{equation*}
with the modified gain operator~$\hat{\mK}$, which is shown to be bounded in~$\Lp^\infty(\mD)$ (Lemma~\ref{lemm:opkern2}). Hence, the operator $R \mapsto \hat{\mK} R - \theta R$ is also bounded for~$\theta \in \Lp^\infty(\mD)$. Classical results in kinetic theory~\cite{Transport} states that there exists a global unique positive solution to this equation, with the bound, for any time~$t \geqslant 0$,
\begin{equation*} 
\|R(t)\|_{\Lp^\infty} \leq \|R(0)\|_{\Lp^\infty} e^{Ct},
\end{equation*}
where the constant~$C$ depends on the operators' kernels as
\begin{equation*} 
C = \sup_{z \in \mD} \left[ \int_{\R^d}  \frac{\sqrt{ \beta(2\pi)^{-1}}}{|u - v|^{d-2}} \cdot  e^{-\beta \left[\frac{|u - v|^2}{8} + \frac{\bigl( |u|^2- |v|^2\bigr)^2 }{8 |u - v|^2 } \right]} \dd u \ - \theta(z) - \nu_\beta(v)\right].
\end{equation*}
This concludes the proof, returning to the variable~$\chi \in \mF_{1, \beta/4}$. The same computation holds backwards in time for~$\eta$, observing that~$M_\beta\eta$ satisfies a linear Boltzmann-like equation with the same collision operator as~$\chi$. Note that even in the case~$\theta = 0$, one can check that the constant~$C$ is positive because of small velocities. In this case nevertheless, one can still prove a maximum principle using more elaborate methods implying a bootstrap argument to fall back in spaces in which a complete spectral study is possible (see for example~\cite[Section~3.1]{2020isa}, or more recently~\cite{Qi2024}). \CQFD 
\end{prof}

%%%--%%%--%%%--%%%--%%%-- Bibliographie --%%%--%%%--%%%--%%%
\vspace{1cm}

\bibliographystyle{abbrv}  % abbrv  acm   alpha (noms complets)    apalike   ieeetr ...
%\bibliography{../../arefs}
\bibliography{arefs}

\end{document}